\documentclass[twoside, 12pt, MSc]{muthesis_2020}
\usepackage[english]{babel}
\usepackage{bookmark,booktabs}
\usepackage{color,bold-extra,mathrsfs,float}
\usepackage{comment,graphics,aliascnt}
\usepackage[dvipsnames]{xcolor}
\usepackage{hyperref}
\usepackage[hyperpageref]{backref}
\hypersetup{
    hypertexnames=false, 
    colorlinks,
    linkcolor={blue},
    citecolor={blue},
    urlcolor={blue},
    pdfauthor={Zhengbo Zhou},
    pdfcreator={Zhengbo Zhou via GNU Emacs and AUCTeX}
}
\renewcommand*{\backref}[1]{}
\renewcommand*{\backrefalt}[4]{%
\ifcase #1 %
\or
   (Cited on p.~#2)%
\else
   (Cited on pp.~#2)%
\fi
}

\usepackage{enumerate,amsmath,amsfonts,amssymb,comment,mathtools}
\usepackage{amsthm}

\usepackage{listings}
\definecolor{codegreen}{rgb}{0,0.6,0}
\definecolor{codegray}{rgb}{0.5,0.5,0.5}
\definecolor{codepurple}{rgb}{0.58,0,0.82}
\definecolor{mygreen}{RGB}{28,172,0}
\definecolor{mylilas}{RGB}{170,55,241}
\definecolor{backcolour}{rgb}{0.95,0.95,1.92}
\lstdefinestyle{mystyle}{
  language=matlab,
  keywordstyle = \linespread{1}\ttfamily\color{blue}\footnotesize,
  basicstyle =   \linespread{1}\ttfamily\footnotesize,
  commentstyle=  \linespread{1}\ttfamily\color{teal}\upshape\footnotesize,
  numberstyle=   \tiny\color{blue},
  stringstyle =  \linespread{1}\color{purple}\footnotesize,
  breakatwhitespace=false,
  breaklines=true,
  captionpos=b,
  keepspaces=true,
  numbers=left,
  numbersep=5pt,
  showspaces=false,
  showstringspaces=false,
  showtabs=false,
  tabsize=4,
  aboveskip=\medskipamount,
  columns=fixed,
  fontadjust=true,
  basewidth=0.5em,
  upquote=true,
}
\def\inline{\lstinline[basicstyle=\upshape\ttfamily]}

\DeclarePairedDelimiter\abs{\lvert}{\rvert}

\usepackage[most]{tcolorbox}
\tcbset{colback=green!40!white, arc=0pt,outer arc=0pt, boxrule=0pt, frame empty, breakable,
        highlight math style= {enhanced, 
            colframe=red,colback=red!10!white,boxsep=0pt}
        }
\newtcbox{\mybox}[1][red]
  {on line, arc = 0pt, outer arc = 0pt,
    colback = #1!10!white, colframe = #1!50!black,
    boxsep = 0pt, left = 1pt, right = 1pt, top = 2pt, bottom = 2pt,
    boxrule = 1pt, bottomrule = 1pt, toprule = 1pt}

\newtheorem{theorem}{Theorem}[chapter]
\newtheorem{lemma}[theorem]{Lemma}
\newtheorem{proposition}[theorem]{Proposition}

\newtheorem{corollary}[theorem]{Corollary}
\newtheorem{algorithm}{Algorithm}

\theoremstyle{definition}
\newtheorem{definition}[theorem]{Definition}
\newtheorem{remark}[theorem]{Remark}
\newtheorem{example}[theorem]{Example}

\numberwithin{equation}{chapter}

\newcommand{\R}{\mathbb{R}} 
\newcommand{\C}{\mathbb{C}} 

\newcommand{\n}{^n}
\newcommand{\Rn}{\R^n}
\newcommand{\mn}{^{m\times n}}
\newcommand{\nn}{^{n\times n}}

\newcommand{\tp}{^T} 
\newcommand{\ctp}{^*}
\newcommand{\inv}{^{-1}}
\DeclareMathOperator{\diag}{diag}
\DeclareMathOperator{\rank}{rank}
\DeclareMathOperator{\tr}{Tr}

\newcommand{\mat}{MATLAB}

\newcommand{\wt}{\widetilde}

\def\eu{\mathrm{e}} 
\def\im{\mathrm{i}} 

\renewcommand{\l}{\lambda}

\def\ycite[#1#2#3#4#5]#6{\cite[$\mit{#1#2#3#4}$#5]{#6}}
\newcommand{\iter}[1]{^{(#1)}} 
\newcommand{\norm}[1]{\|{#1}\|_2} 
\newcommand{\sign}[1]{\mathrm{sign}\left({#1}\right)}
\newcommand{\off}{\textrm{off}} 

\newcounter{mylineno}
\makeatletter
\let\oldtabcr\@tabcr
\def\nonumberbreak{\oldtabcr\hspace{3.5pt}}
\def\mynewline{\refstepcounter{mylineno}%
  \llap{\footnotesize\arabic{mylineno}\hspace{5pt}}%
}

\gdef\@tabcr{\@stopline \@ifstar{\penalty%
    \@M \@xtabcr}\@xtabcr\mynewline}

\newenvironment{code}{%
  \mathcode`\:="603A  
  
  \setcounter{mylineno}{0}
  \par
  \upshape
  \begin{list} 
  {} {\leftmargin = 1cm}
  \item[]
  \begin{tabbing}

  \hspace*{.3in} \= \hspace*{.3in} \=
  \hspace*{.3in} \= \hspace*{.3in} \= \kill
  \mynewline
}{\end{tabbing}\end{list}}
\makeatother
\mathcode`@="8000
{\catcode`\@=\active\gdef@{\mkern1mu}}

\newcommand{\ift}{\text{ if }}
\newcommand{\gnorm}[1]{\|{#1}\|}
\renewcommand{\diag}{\mathrm{diag}}
\newcommand{\eref}[1]{\eqref{#1}}

\def\mat{MATLAB}

\begin{document}
\pagenumbering{roman}

\title{A Mixed Precision Eigensolver Based on the Jacobi Algorithm}
\author{Zhengbo Zhou}
\school{Mathematics}
\faculty{Science and Engineering}
\submitdate{Submit date: 10 September 2022}
\def\wordcount{12038}

\beforeabstract
The classic method for computing the spectral decomposition of a real symmetric
matrix, the Jacobi algorithm, can be accelerated by using mixed precision
arithmetic. The Jacobi algorithm is aiming to reduce the off-diagonal
entries iteratively using Givens rotations. We investigate how to use the 
low precision to speed up this algorithm based on the approximate 
spectral decomposition in low precision.  

We first study two different index choosing techniques, classical and
cyclic-by-row, for the Jacobi algorithm. Numerical testing suggests that
cyclic-by-row is more efficient. Then we discuss two different methods of
orthogonalizing an almost orthogonal matrix: the QR factorization and the
polar decomposition. For polar decomposition, we speed up the Newton
iteration by using the one-step Schulz iteration. Based on numerical
testing, using the polar decomposition approach (Newton--Schulz iteration)
is not only faster but also more accurate than using the QR factorization.  

A mixed precision algorithm for computing the spectral decomposition of a real
symmetric matrix at double precision is provided. In doing so we compute
the approximate eigenvector matrix $Q_\ell$ of $A$ in single precision
using \texttt{eig} and \texttt{single} in MATLAB. We then use the
Newton--Schulz iteration to orthogonalize the eigenvector matrix $Q_\ell$
into an orthogonal matrix $Q_d$ in double precision. Finally, we apply the
cyclic-by-row Jacobi algorithm on $Q_d^TAQ_d$ and obtain the
spectral decomposition of $A$. At this stage, we will see, from the testings,
the cyclic-by-row Jacobi algorithm only need less than 10 iterations to
converge by utilizing the quadratic convergence. The new mixed precision
algorithm requires roughly 30\% of the time used by the Jacobi algorithm on
its own. 
\afterabstract

\prefacesection{Acknowledgements}
It is my pleasure to acknowledge the advices from my supervisors Prof.
Nicholas J. Higham and Prof. Fran\c coise Tisseur. I wouldn't start my
research in numerical linear algebra without their guidance and this thesis
cannot be completed without their unlimited patience and insightful
feedback.  

I would also like to thank all my friends for providing emotional support.
In addition, two special thanks to Zeyu Wu, for his suggestion on coding
aspects, and Anheng Xi, for the translation of several German articles from
Jacobi and Sch{\"o}nhage.   

Finally, I thank my parents, Guangyuan Zhou and Ling Li, for their
unconditional love and for giving me the chance to study abroad. By
standing upon the shoulders of my parents, I have seen further than them.

\afterpreface

\cleardoublepage
\pagenumbering{arabic}

\chapter{Introduction}

Modern hardware increasingly supports low precision floating-point
arithmetic. Using low precision, we have opportunities to accelerate linear
algebra computations. A symmetric eigenproblem is solving
\begin{equation}\notag
  Ax = \lambda x,
\end{equation}
with $x\neq 0$ and the matrix $A\in\R\nn$ is symmetric. Equivalently, the
problem is aiming to compute the spectral decomposition
\begin{equation}\notag
  A = Q\varLambda Q\tp,\quad A = A\tp \in \R\nn,\quad Q\tp Q = QQ\tp = I,
\end{equation}
and $\varLambda$ is the diagonal matrix with the eigenvalues of $A$ along
the diagonal.

Traditional ways of solving such a problem include the Jacobi algorithm,
the power method and the QR algorithm~\ycite[2000,
Section~5--7]{GvdV2000Eigenvaluecomputation20th}. This thesis is
concerned with the 
Jacobi algorithm proposed by Jacobi in 1846 which is an iterative method
that reduces the off-diagonal entries at each step. In addition, we
investigate how the usage of low precision arithmetic can speed up the
Jacobi algorithm.

In 2000, Drmač and Veselić~\ycite[2000]{DV2000Approximateeigenvectorsas}
provided
a way to speed up the Jacobi algorithm on the real symmetric matrix $A$.
Let $U$ be the approximate eigenvector matrix of $A$, then applying the
Jacobi algorithm on $U\tp A U$ can speed the reduction up.

In Chapter~\ref{chap:norms} we provide essential definitions and examples
including norms, orthogonal matrix and the singular value decomposition. We
introduce two important classes of the orthogonal matrix, the Givens
rotation and the Householder transformation.

The Jacobi algorithm is introduced in Chapter~\ref{chap:jacobi_algorithm}.
We describe the derivation of the Jacobi algorithm and its convergence
rate. In Section~\ref{sec:classical-jacobi} and~\ref{sec:cyclic-jacobi} we
study how to choose the pair of index required by the Jacobi algorithm.

The classical Jacobi algorithm, mentioned in
Section~\ref{sec:classical-jacobi}, chooses the largest off-diagonal entry
as the pair of index and the cyclic-by-row Jacobi algorithm, mentioned in
Section~\ref{sec:cyclic-jacobi}, chooses the indices row by row but
restricts to the superdiagonal part of the matrix. Although the classical
Jacobi algorithm reduces the off-diagonal entries most optimally, the
cyclic-by-row Jacobi algorithm does not require expensive off-diagonal
searches. Finally, we conduct numerical tests to examine the speed and
conclude that we shall stick with the faster cyclic-by-row Jacobi
algorithm.

We provide a procedure to get an approximate eigenvector matrix in
Chapter~\ref{chap:orthogonalisation}. Firstly, compute the
spectral decomposition at low precision $A \approx Q_\ell D_\ell Q_\ell\tp$ 
such that
\begin{equation}\notag
  \gnorm{Q_\ell\tp Q_\ell - I} \lesssim nu_\ell,\quad \gnorm{AQ_\ell -
  Q_\ell D_\ell}\lesssim nu_\ell  \gnorm{A},
\end{equation}
where $u_\ell$ is the unit roundoff in low precision. Then we orthogonalize
$Q_\ell$ to $Q_d$ such that $\gnorm{Q_d\tp Q_d - I} \lesssim nu_d$ where
$u_d$ is the unit roundoff in double precision (``Unit roundoff'' - defined
in Section~\ref{sec:IEEE}). We can use $Q_d$ as the approximate eigenvector
matrix of $A$ and precondition $A$ into $A_{\text{cond}} = Q_d\tp A Q_d$,
then the Jacobi algorithm applied on $A_{\text{cond}}$ should converge in a
few sweeps. In Section~\ref{sec:Householder-QR}
and~\ref{sec:polar-decomposition} we review two methods for orthogonalizing
a matrix: the Householder QR factorization and the polar decomposition.
Particularly, in Section~\ref{sec:polar-decomposition}, we study the Newton
iteration approach and its variant Newton--Schulz iteration to compute the
unitary factor of the polar decomposition. MATLAB codes are presented and
after several numerical testings, among the Householder QR factorization
approach, the Newton iteration approach and the Newton--Schulz iteration
approach, the latter method is the most accurate and efficient one to use.

Finally, we assemble all the ingredients and produce
Algorithm~\ref{alg:jacobi-preconditioned} in
Chapter~\ref{chap:mixed-precision}, a mixed precision Jacobi algorithm. In
Section~\ref{sec:quadratic-conv}, we review some literature on the
quadratic convergence of the cyclic-by-row Jacobi algorithm. We observe
from the numerical testings that for any symmetric matrix $A\in\R\nn$,
after preconditioning, the Jacobi algorithm should converge within 10
iterations. The numerical testings suggest the Jacobi algorithm only needs
two iterations to converge and there is a 75\% reduction of time using
preconditioned Jacobi algorithm compare to the Jacobi algorithm on its own.

\section{Notations}
We use the Householder notations:
\begin{itemize}
\item Use capital letters for matrices.
\item Use lower case letters for vectors.
\item Use lower case Greek letters for scalars.
\end{itemize}

We use $\lambda(A)$ and $\sigma(A)$ to represent the eigenvalues and the
singular values of $A$. We denote $\sigma_1(A)$ or $\sigma_1$ if the matrix
$A$ has been specified as the largest singular value of $A$.

We say that the $Ax = b$ problem is well conditioned if small changes in
the data $A,b$ always make small changes in the solution $x$; otherwise, it
is ill-conditioned. We use $\kappa(A) = \gnorm{A}\gnorm{A\inv}$ as a
measure to visualize the condition of a problem. (Norms will be introduced
in Section~\ref{sec:norm})

We say a matrix $A\in\R\nn$ is almost orthogonal if
$\gnorm{A\ctp A - I} \lesssim nu_\ell$ where $I$ is the $n\times n$
identity matrix. A matrix $A\in\R\nn$ is almost diagonal if
$\gnorm{A - \diag(A)}/\gnorm{A} \lesssim nu_\ell$ and this quantity is zero
if $A$ is diagonal.

\section{IEEE standard}\label{sec:IEEE}
The unit roundoff measures the relative error due to approximation in
floating-point arithmetic. Throughout this thesis, we will only consider
two precisions: single precision and double precision~\ycite[2002,
Section~2.1]{high:ASNA2}. We will denote $u_\ell$ as the unit roundoff at
single precision and $u_d$ as the unit roundoff at double precision. IEEE
standard 754-2019~\ycite[2019, Table~3.5]{IEEE2019} shown in
Table~\ref{tab.unit-roundoff} will be used. For simplicity, any matrix with
subscript $\ell$ will be at single precision and any matrix with subscript
$d$ will be at double precision.
\begin{table}[ht]
\centering
\caption{Unit roundoff defined by IEEE standard 754-2019.}
\label{tab.unit-roundoff}
\begin{tabular}{ll}
  \toprule
  Type & Unit Roundoff \\ \midrule 
  Single precision & $2^{-24} \approx $ 5.96e-8 \\
  Double precision & $2^{-53} \approx$ 1.11e-16 \\
  \bottomrule
\end{tabular}
\end{table}

\section{Reproducible Research}

In this section, we will present our testing environment and system
specifications. With these informations, the reader will have the ability
to reproduce any results we have in this thesis.

Testing environment and system specifications:
\begin{itemize}
\item System OS: Windows 10.
\item MATLAB version: 9.12.0.2009381 (R2022a) Update 4~\cite{MATLAB:2022}.
\item CPU: Intel(R) Core(TM) i5-9600KF CPU @ 3.70GHz.
\item GPU: NVIDIA GeForce RTX 2060.
\item All the testings that require the random number generator will
initiate after the MATLAB command \inline{rng(1,'twister')}\footnote{MATLAB
  built-in function to control the random number generator, see
  \url{https://www.mathworks.com/help/matlab/ref/rng.html}.}.
\end{itemize}


\chapter[Preliminaries]{Matrix Norms, Orthogonal Matrices and Singular
  Value Decomposition}\label{chap:norms} 

In this chapter, we will introduce several concepts that are important for
our analysis. In the first section, vector norms and matrix norms will be
presented. Then we will discuss an important class of matrices, orthogonal
matrices. Finally, a crucial decomposition, singular value decomposition,
will be discussed.  

\section{Norms}\label{sec:norm}

Matrix norms are essential for matrix algorithms. For example, the matrix
norms provide a way of measuring the difference between two matrices and
help us to check if we have the desired solution in finite arithmetic. In
this section, we will do our demonstrations in the real world but the
complex case is similar. Also, we will only introduce the norm for square
matrices, but these theories are all applicable to any general matrices
with careful attention to matching dimensions. 

\subsection{Vector Norms}\label{sec:vector-norm}

Before introducing matrix norm, a brief illustration of vector norm is necessary. 
\begin{definition}
  [vector norm] 
  A vector norm on $\R\n$ is a function $\gnorm : \R\n \to \R$ such that it
  satisfies the following properties 
  \begin{enumerate}
    \item $\gnorm{x} \geq 0$ for all $x\in\R\n$.
    \item $\gnorm{x} = 0$ if and only if $x = 0$.
    \item $\gnorm{\lambda x} = \abs{\lambda} \gnorm{x}$ for all $\lambda \in \R$ and $x \in \Rn$.
    \item $\gnorm{x + y} \leq \gnorm{x} + \gnorm{y}$ for all $x,y\in \Rn$.
  \end{enumerate}
\end{definition}

\begin{example}
  A useful class of vector norm is the $p$-norm 
  \begin{equation}\notag
    \gnorm{x}_p = \left({\sum_{i = 1}^n \abs{x_i}^p}\right)^{1/p},\quad p\geq 1.
  \end{equation}
  One particular example of the $p$-norm is the Euclidean norm or simply $2$-norm defined by taking $p=2$
  \begin{equation}\notag
    \norm{x} =  \left({\sum_{i = 1}^n \abs{x_i}^2}\right)^{1/2} = \sqrt{x\tp x}.
  \end{equation}
  The latter equality is obvious.
\end{example}

\subsection{Matrix Norms}

Let us denote $\R\nn$ as $n\times n$ matrices with real entries.

\begin{definition}
  [matrix norm]
  A matrix norm on $\R\nn$ is a function $\gnorm{\cdot}:\R\nn \to \R$ satisfying the following properties 
  \begin{enumerate}
    \item $\gnorm{A} \geq 0$ for all $A\in\R\nn$.
    \item $\gnorm{A} = 0$ if and only if $A = 0$.
    \item $\gnorm{\lambda A} = \abs{\lambda}\gnorm{A}$ for all $\lambda \in \R$ and $A\in\R\nn$.
    \item $\gnorm{A + B} \leq \gnorm{A} + \gnorm{B}$ for all $A,B \in \R\nn$. 
  \end{enumerate}
\end{definition}

One simple matrix norm is the Frobenius norm 
\begin{equation}
  \gnorm{A}_F = \left(\sum_{i = 1}^n \sum_{j = 1}^n \abs{a_{ij}}^2 \right)^{1/2} = (\tr(A\tp A))^{1/2}.
\end{equation}
We will use the Frobenius norm intensively in the next section on the Jacobi algorithm.

Another important class of matrix norms is called induced norms or subordinate norms. Given a vector norm defined in Section~\ref{sec:vector-norm}, the corresponding induced norm is defined as 
\begin{equation}\notag
  \gnorm{A} = \max_{\gnorm{x} = 1}\gnorm{Ax}.
\end{equation}

\begin{example}
  The matrix 2-norm is induced by the vector 2-norm. From the definition of the induced norm, 
  \begin{equation}\notag
      \norm{A} = \max_{\norm{x} = 1} \norm{Ax} = \sqrt{x_\mu\tp A\tp A x_\mu} = \mu,
  \end{equation}
  for some positive number $\mu$, and $x_\mu$ is the corresponding chosen $x$. Since $\norm{x} = \sqrt{x\tp x}$ for $x\in\Rn$, hence we have $x_\mu\tp A\tp A x_\mu = \mu^2$. Also, by $\norm{x_\mu} = x_\mu\tp x_\mu = 1$, we can pre-multiply $x_\mu$ at both sides and we have  $A\tp A x_\mu = \mu^2 x_\mu$ which implies $\mu$ is the eigenvalue of $A\tp A$. In case we are choosing $x_\mu$ such that $\mu$ is maximum, therefore we have the definition of the matrix 2-norm 
  \begin{equation}\label{eq:example-2.4}
    \norm{A} = \sqrt{\l_{\max}(A\tp A)} = \sigma_1,
  \end{equation} 
  where $\lambda_{\max}(A)$ and $\sigma_1$ denotes the largest eigenvalue and singular value of $A$. 
\end{example}

\begin{definition}
  [consistent norms]
  A norm is consistent if it is submultiplicative
  \begin{equation}
    \gnorm{AB} \leq \gnorm{A} \gnorm{B},
  \end{equation}
  for all $A,B$ such that the product $AB$ defines.
\end{definition}

\begin{example}
  The Frobenius norm and all subordinate (induced) norms are consistent. This is useful for constructing upper bounds.
\end{example}

\section{Orthogonal Matrices}

Recall the definition of the orthogonal matrix,
\begin{definition}
  [orthogonal matrices]
  A matrix $Q\in\R\nn$ is said to be orthogonal if
  \begin{equation}\notag
    Q\tp Q = QQ\tp = I_n,
  \end{equation}
  where $I_n$ denotes the identity matrix of size $n$. If $Q\in\C\nn$, then we call $Q$ unitary rather than orthogonal and we denote its conjugate transpose as $Q\ctp$.
\end{definition}

\begin{theorem}
  The vector 2-norm is invariant under orthogonal transformations.
\end{theorem}

\begin{proof}
  Given an orthogonal matrix $Q\in\R\nn$ and any vector $x\in\Rn$, we have 
  \begin{equation}\notag
    \norm{Qx} = \sqrt{x\tp Q\tp Q x} = \sqrt{x\tp x} = \norm{x},
  \end{equation}
  which proves the theorem.
\end{proof}

\begin{theorem} 
  \label{thm:matrix-invariant-norm}
  The Frobenius norm and matrix 2-norm are the orthogonally invariant norm. Namely, the following norm equality holds for these two norms
  \begin{equation}\notag
    \gnorm{UAV} = \gnorm{A},\quad \text{$A\in\R\nn$, $U$ and $V$ are orthogonal.}
  \end{equation}
\end{theorem}

\begin{proof}
By the definition of the Frobenius norm, we have 
  \begin{equation}
      \label{eq:thm2-2-proof1}
      \begin{aligned}
          \|UAV\|_F^2 &= \tr\big(U AVV\tp A\tp U\tp \big) = \tr\big( U AA\tp U\tp\big).
      \end{aligned}
  \end{equation}
  By the properties $\tr(AB) = \tr(BA)$, \eref{eq:thm2-2-proof1} simplifies to 
  \begin{equation}
      \label{eq:thm2-2-proof2}
      \|U AV\|_F^2 = \tr\big(A\tp U\tp U A\big) = \tr\big(A\tp A\big) = \tr\big(AA\tp\big) = \|A\|_F^2.
  \end{equation}
  Hence we finished the proof for the Frobenius norm.

  By the definition of the matrix 2-norm, we have 
  \begin{equation}\notag
    \norm{UAV}^2 = \l_{\max}(V\tp A\tp U\tp UAV) = \l_{\max}({V\tp A\tp AV}).
  \end{equation}
  Notice that if $V$ is orthogonal, then $A\tp A$ and $V\tp A\tp A V$ share the same eigenvalues. Therefore $\l_{\max}(V\tp A\tp AV) = \l_{\max} (A\tp A)$, therefore we have 
  \begin{equation}\notag
    \norm{UAV}^2 = \lambda_{\max}(A\tp A) = \norm{A}^2.
  \end{equation}
  Hence we proved the theorem.
\end{proof}

\subsection{Givens Rotations}\label{sec:givens-rotation}

\begin{definition}
  [Givens rotation]
  A Givens rotation has the form 
  \begin{equation}\label{eq:givens}
    (G(i, k, \theta))_{pq} =
    \begin{cases}
      1 & \text{if $p,q \notin \{i,k\}$},\\
      c & \text{if $p = q$ and $p \in \{i,k\}$,} \\
      s & \text{if $p = i$ and $q = k$,}\\
      -s & \text{if $p = k$ and $q = i$,}\\
      0 & \text{otherwise.}
    \end{cases}
  \end{equation}
  where $c = \cos(\theta)$ and $s = \sin(\theta)$ for some $\theta$.
\end{definition}

Given a vector $y\in\R^n$ and a Givens rotation $G(i,j,\theta)$, the product $G(i,j,\theta) y$ rotate $y$ through $\theta$ radians clockwise in the $(i,j)$ plane. Hence it is obvious that we can make either $y_i$ or $y_j$ become zero by manipulating $\theta$. This idea will be utilized in Section~\ref{sec:Deriv and Conver} in order to eliminate $(i,j)$th entry of a symmetric matrix.

\begin{proposition}
  Givens rotations are orthogonal.
\end{proposition}
\begin{proof}
  Let $G = G(i,k,\theta)$ be a Givens rotation. By using the trigonometric identity $\cos^2(\theta) + \sin^2(\theta) = 1$, we have the following 
  \begin{equation}\label{eq:1.2}
      G_{ik}=G([i,k],[i,k]) = 
      \begin{bmatrix}
          c & s \\
          -s & c 
      \end{bmatrix}
      \quad \to \quad 
      G_{ik}\tp G_{ik} = 
      \begin{bmatrix}
          1 & 0 \\
          0 & 1
      \end{bmatrix}
      = G_{ik}G_{ik}\tp .
  \end{equation}
  We can calculate the following matrix products $G\tp G$ and $GG\tp$ using \eref{eq:1.2}.
  \begin{equation}\notag
      G\tp G = GG\tp = I.
  \end{equation}
  Hence $G$ is orthogonal.
\end{proof}

Notice that the matrix $G(p,q,\theta)$ can be written as
\begin{equation}\label{eq:1.4}
  \begin{aligned}
      G(p,q,\theta) = [e_1,\dots,e_{p-1},e_{pq},e_{p+1},\dots,e_{q-1},e_{qp},e_{q+1},\dots,e_{n}],
  \end{aligned}
\end{equation}
where $e_i$ is the $i$th column of the $n\times n$ identity matrix if $i\neq p,q$ and $e_{pq}$ and $e_{qp}$ are defined as 
\begin{equation}\label{eq:1.5}
  [e_{pq}]_k = \begin{cases}
      c & \ift k = p \\
      -s & \ift k = q\\
      0 & \ift k \neq p,q
  \end{cases}, \qquad 
  [e_{qp}]_k = \begin{cases}
      s & \ift k = q \\
      c & \ift k = p\\
      0 & \ift k \neq p,q
  \end{cases}.
\end{equation}

\begin{proposition}\label{prop:1.3}
  Suppose $A\in \R\nn$ and $G(p,q,\theta)$ is a Givens rotation. If we define $B$ as $$B = G(p,q,\theta)\tp AG(p,q,\theta),$$ 
  then $B$ and $A$ are the same except in rows and columns $p$ and $q$.
\end{proposition}

\begin{proof} 
  Using the notations in \eref{eq:1.4} and \eref{eq:1.5}, matrix $B$ can be written as
  \begin{equation}\notag
      \begin{aligned}
          B =
          \begin{bmatrix}
              e_{1}\tp A e_{1} & \cdots & e_{1}\tp A e_{p q} & \cdots & e_{1}\tp A e_{q p}  &\cdots & e_{1}\tp A e_{n} \\
              \vdots &  & \vdots &  & \vdots  &  & \vdots \\
              e_{pq}\tp A e_{1} & \cdots & e_{pq}\tp A e_{p q} & \cdots & e_{pq}\tp A e_{q p}  &\cdots & e_{pq}\tp A e_{n} \\
              \vdots &  & \vdots &  & \vdots  &  & \vdots \\
              e_{qp}\tp A e_{1} & \cdots & e_{qp}\tp A e_{p q} & \cdots & e_{qp}\tp A e_{q p}  &\cdots & e_{qp}\tp A e_{n} \\
              \vdots &  & \vdots &  & \vdots  &  & \vdots \\
              e_{n}\tp A e_{1} & \cdots & e_{n}\tp A e_{p q} & \cdots & e_{n}\tp A e_{q p}  &\cdots & e_{n}\tp A e_{n} \\
          \end{bmatrix}.
      \end{aligned}
  \end{equation}
  Hence $B_{ij} = e_i\tp Ae_j = A_{ij}$ if $i,j \neq p,q$ and for $i,j = p$ or $q$, $B_{ij}\neq A_{ij}$ in general.
\end{proof}

\subsection{Householder Transformations}\label{sec:Householder}

Another important class of orthogonal matrix is the Householder transformations. Unlike the Givens rotation which replaces a specific entry into zero, Householder transformation can change multiple entries into zeros. In this section, we will first discuss how this matrix can set entries to zeros and then a MATLAB implementation will be present together with numerical testing.

\begin{definition}
  A Householder transformation (synonyms are Householder matrix and Householder reflector) is an $n\times n$ matrix $P$ of the form
  \begin{equation}\notag
      P = I - 2\frac{vv\tp}{v\tp v},\quad 0\neq v\in\R^n.
  \end{equation}
  The vector $v$ is the Householder vector.
\end{definition}

Householder transformation has useful properties,
\begin{enumerate}
  \item Symmetry:
  \begin{equation}\notag
      P\tp = I\tp - \frac{2}{v\tp v}(vv\tp)\tp = I - \frac{2}{v\tp v}vv\tp = P.
  \end{equation}
  \item Orthogonality:
  \begin{equation}\notag
      \begin{aligned}
          P\tp P = PP\tp = PP & = \left(I - \frac{2}{v\tp v}vv\tp\right)\left(I - \frac{2}{v\tp v}vv\tp\right)\\
          & = I-\frac{4}{v\tp v}vv\tp + \frac{4}{(v\tp v)^2}vv\tp vv\tp\\
          & = I - \frac{4}{v\tp v}vv\tp + \frac{4}{v\tp v}vv\tp = I.
      \end{aligned}
  \end{equation}
\end{enumerate}

\subsubsection{Transforming a Vector}
The Householder transformation can be used to zero components of a vector. Suppose we have a vector $x\in\R^n$ and we want
\begin{equation}\label{eq.2.6}
  y = Px = \left(I - \frac{2}{v\tp v}vv\tp\right) x = x - \frac{2v\tp x}{v\tp v}v
\end{equation}
where $y(2:n) = {0}$. Since $P$ is orthogonal, we require $\norm{y} = \norm{x}$, hence $y = \pm \norm{x}e_1$ where $e_1$ is the first column of the $n\times n$ identity matrix. We can rewrite $v$ as
\begin{equation}\notag
  v = \frac{1}{\alpha} x - \frac{1}{\alpha} y, \quad \alpha = \frac{2v\tp x}{v\tp v}v,
\end{equation}
and since $P$ is independent of the scaling of $v$, hence for $x\neq y$, we can set $\alpha = 1$ and we can choose $v = x -y$. We can check our choice by computing $Px$:
Firstly, we can compute $v\tp x$ and $v\tp v$
\begin{equation}\notag
  v\tp x = (x\tp - y\tp) x=x\tp x - y\tp x,\quad v\tp v = (x\tp - y\tp)(x-y) = 2x\tp x - 2x\tp y.
\end{equation}
Substitute these two components into \eref{eq.2.6}
\begin{equation}\notag
  \begin{aligned}
      Px & = x - \frac{2v\tp x}{v\tp v}v = x - v = x - (x-y) = y.
  \end{aligned}
\end{equation}

The choice of the sign of $y$ depends on the sign of the first entry of $x$. Suppose $x$ is dominant by its first entry, namely $x$ is close to a multiple of $e_1$, then $v = x - \sign{x_1}\norm{x}e_1$ can have a small norm due to cancellation and this may result a large relative error in $2/(v\tp v)$. This relative error can be avoided by setting
\begin{equation}\notag
  v = x + \sign{x_1}\norm{x}e_1
\end{equation}
as suggested in many textbooks, such as \ycite[2002, Section~19.1] {high:ASNA2} and \ycite[2013, Section~5.1.3]{van2013mc}. Then we get
\begin{equation}\label{eq.Px}
  Px = y = x - v = -\sign{x_1}\norm{x}e_1.
\end{equation}

Given $x\neq {0} \in\R^n$, instead of generating the Householder matrix $P$ such that it satisfies \eref{eq.Px}, we should only generate the Householder vector $v$ and the coefficient $\beta = 2/(v\tp v)$. The benefit of doing so can be seen from the Householder QR factorization algorithm in section~\ref{subsec.houseQR}. Based on the previous discussion, we can conclude these into  Algorithm~\ref{alg.house}.

\begin{algorithm}
\label{alg.house}
Given $x\neq{0} \in\R^n$, this algorithm computes $v\in\R^n$ and $\beta\in\R$ such that $P = I_n - \beta vv\tp$ is orthogonal and $Px = -\sign{x_1} \norm{x}e_1$.
\begin{code}
Form $\sigma = x(2:n)\tp x(2:n),\ v = [1,x(2:n)\tp]\tp$ \\
Calculate $\mu = \norm{x} =  \sqrt{x_1^2 + \sigma}$ \\
if \= $x_1 \geq 0$\\
\> $v_1 = x_1 + \mu$\\
else \\
\> $v_1 = x_1 - \mu = -\sigma/(x_1 + \mu)$ \\
end \\
$\beta = 2/(v_1^2 + \sigma)$
\end{code}
\end{algorithm}

Notice that, line 6 rewrites the form of $v_1$ to avoid subtractive cancellation as suggested by \ycite[2013, Algorithm~5.1.1]{van2013mc}
\begin{equation}\notag
  v_1 = x_1 - \mu = \frac{(x_1 - \mu)(x_1 +\mu)}{(x_1 + \mu)} = \frac{x_1^2 - (x_1^2 + \sigma)}{x_1 + \mu} = \frac{-\sigma}{x_1 + \mu}.
\end{equation}

\subsubsection{Implementation and Testing}

During the \mat\ implementation, we do not need to create a new vector $v$. It is sufficient to just overwrite the values of $x$ since the only difference between $v$ and $x$ in Algorithm~\ref{alg.house} is the first entry of these two vectors.

\begin{lstlisting}
function [x,b] = house(x)
n = length(x); sigma = x(2:n)'*x(2:n); 
mu = sqrt(x(1)*x(1) + sigma);
if x(1) >= 0, x(1) = x(1) + mu;
else, x(1) = -sigma/(x(1) + mu); end
b = 2/(x(1)*x(1) + sigma); 
end
\end{lstlisting}

To test the above function, we first generate a vector $x\in\R^5$ using \inline{randn(5,1)}. 
\begin{lstlisting}
x' = 
  -1.3499e+00, 3.0349e+00, 7.2540e-01, -6.3055e-02, 7.1474e-01
\end{lstlisting}
Using the above function and call \inline{[v,b] = house(x)}, we can construct $P$ using the formula $P = I_n - \beta vv\tp$. Then we examine the product $Px$ and we have 
\begin{lstlisting}
(P * x)' = 
  3.4748e+00            0   1.7694e-16  -9.1073e-18   1.1102e-16
\end{lstlisting}
The result satisfies our expectations. 

The Householder transformation can quickly introduce zeros into a vector using orthogonal matrices and this is useful when we study the QR factorization.

\section{Singular Value Decomposition}\label{sec:svd}

\begin{theorem}
    [singular value decomposition]
    \label{thm.svd}
    If $A\in\R\mn$, $m\geq n$, then there exists two orthogonal matrices $U\in\R^{m\times m}$ and $V\in\R\nn$ such that 
    \begin{equation}\label{eq.svd}
        A = U\Sigma V\tp,\quad \Sigma = \diag(\sigma_1,\dots,\sigma_p)\in\R\mn,\quad p = \min\{m,n\},
    \end{equation}
    where $\sigma_1,\dots,\sigma_p$ are all non-negative and arranged in non-ascending order. We denote~\eqref{eq.svd} as the singular value decomposition (SVD) of $A$ and $\sigma_1,\dots,\sigma_p$ are the singular values of $A$.
\end{theorem}

From now, without explicit mention, $A$ will be a square matrix of full rank (our input $Q_\ell$ is always square, full rank matrix). Therefore, if we have the singular value decomposition $A = U\Sigma V\tp$, then $U,V\in\R\nn$ are unitary and $\Sigma = \diag(\sigma_1,\dots,\sigma_n)\in\R\nn$ where $\sigma_1\geq \dots \geq \sigma_n > 0$. 

From the definition, the singular values of $A$ can be computed via computing the eigenvalues of $A\tp A = V\varLambda V\tp$, where $\varLambda = \Sigma^2$. Hence, the singular values of $A$ are the positive roots of the eigenvalues of $A\tp A$. Based on this relationship, we have the following properties.

\begin{corollary}\label{col.2norm}
    If $A\in\R\nn$ of full rank and normal, then $\sigma(A) = \abs{\lambda(A)}$.
\end{corollary}

\begin{proof}
  If $A$ is normal, then $A\tp A = AA\tp$. By spectral theorem, if $A$ is normal, then $A = U\varLambda U\tp$ where $U$ is unitary and $\varLambda$ is a diagonal matrix with eigenvalues of $A$ on the diagonal. Therefore, the singular values of $A$ are 
    \begin{equation}\notag
        \sigma(A) = \sqrt{\lambda(A\tp A)} = \sqrt{\l(U\varLambda ^2U\tp)} = \sqrt{\l(\varLambda^2)} = \sqrt{\l(\varLambda)^2} = \abs{\l(A)}.
    \end{equation}
    Combine this Corollary and \eqref{eq:example-2.4}, we have: if $A$ is normal, then $\norm{A} = \max\{\abs{\l(A)}\}$.
\end{proof}

\section{Test Matices}\label{sec:test-matrices}
Throughout this thesis, we will perform lots of MATLAB testing on our algorithms, hence we want to generate real symmetric matrices $A$ with desired condition number and different eigenvalue distributions.

To accomplish such a goal, we use the MATLAB built-in function \inline{gallery}\footnote{A MATLAB built-in function to generate test matrices. For the full manual, you may refer to \url{https://www.mathworks.com/help/matlab/ref/gallery.html}.}. If we use \inline{gallery('randsvd',n, -kappa)} where \inline{kappa} is a positive integer, then it will generate a symmetric positive definite matrix $A\in\R\nn$ with $\kappa(A) = \mathtt{kappa}$. Our code \inline{my_randsvd} from Appendix~\ref{app:myrandsvd} is the modified version of \inline{gallery('randsvd')} and we deliberately change some eigenvalues to negative such that the output matrix will not necessarily be positive definite. It provides two different eigenvalue distributions.
\begin{enumerate}
  \item Mode \texttt{'geo'}, \inline{A = my_randsvd(n, kappa, 'geo')}: The output matrix $A\in\R\nn$ with $\kappa(A) = \mathtt{kappa}$ and the magnitude of the eigenvalues of $A$ are geometrically distributed.
  \item Mode \texttt{'ari'}, \inline{A = my_randsvd(n, kappa, 'ari')}: The magnitude of the eigenvalues of $A$ are arithmetically distributed.
\end{enumerate}

In Chapter~\ref{chap:mixed-precision}, we will test our algorithm on
matrices with different distributions. For the rest of the thesis, we will
use mode \texttt{'geo'} by default.

\chapter{Jacobi Algorithm}\label{chap:jacobi_algorithm}

The Jacobi algorithm for eigenvalues are first published by Jacobi in 1846~\cite{Jacobi-original-paper-1846} and became widely used in the 1950s after the computer is invented. In this chapter, we will first derive the Jacobi algorithm and discuss its linear convergence. Then the rest of this chapter will focus on two different Jacobi algorithms: Classical and Cyclic-by-row. MATLAB implementation and testing will be given.

\section{Derivation and Convergence}\label{sec:Deriv and Conver}

Given a symmetric matrix $A\in \R\nn$, the idea of the Jacobi algorithm is that at the $k$th step, we use $A^{(k+1)} = Q_{k}\tp A^{(k)}Q_{k}$ to replace $A\iter{k}$, where $Q_{k}$ is orthogonal. We aim to have the property that the off-diagonal entries of $A^{(k+1)}$ are smaller than the off-diagonal entries of $A^{(k)}$. To do this, we use the Givens rotation matrix defined in section~\ref{sec:givens-rotation}.

\begin{definition}
  [$\off$ operator]
  We define the quantity
  \begin{equation}\notag
      \off(A) \coloneqq \sqrt{\|A\|^2_F - \sum_{i=1}^n a^2_{ii}} = \sqrt{\sum_{i=1}^n \sum_{\substack{\text{$j=1$}\\\text{$j\neq i$}}}^n a^2_{ij}},
  \end{equation}
  which is the Frobenius norm of the off-diagonal elements.
\end{definition}

The aim of the Jacobi algorithm on $A$ can be translated to reduce the quantity $\off(A)$. This can be done in 2 steps. 
\begin{enumerate}
  \item Choosing a pair of index $(p,q)$. We assume that $1\leq p < q\leq n$.
  \item Overwriting $A$ by applying $G(p,q,\theta)$ on the matrix $A$ in the way of $$A' = G(p,q,\theta)\tp AG(p,q,\theta),$$ where $\theta$ is chosen such that $A'_{pq}= A'_{qp} = 0$.
\end{enumerate}
Using the Jacobi algorithm on $A$, we produce a sequence of matrices $\big\{A^{(k)}\big\}_{k=0}^\infty$, where 
\begin{equation}\notag
  \begin{cases}
      A^{(0)} = A,  \\
      A^{(k)} = G_{k-1}(p,q,\theta)\tp A^{(k-1)}G_{k-1}(p,q,\theta), & \text{for $k = 1,2,\dots$.}
  \end{cases}
\end{equation}

This set of iterations satisfy

\begin{equation}
  \label{eq:jaco-condition}
  \lim_{k\to\infty}\off(A^{(k)}) = 0.
\end{equation}

In the rest of this section, we will discuss how we can achieve \eqref{eq:jaco-condition} using step 2 and how we can construct such a matrix $G(p,q,\theta)$ given the choice $(p,q)$. Then in the next two sections, we present two different methods of choosing the index $(p,q)$
\begin{itemize}
  \item The Classical Jacobi Algorithm.
  \item The Cyclic-by-row Jacobi Algorithm.
\end{itemize}

Considering the subproblem which only involves four entries of $A$.
\begin{equation}
    \label{eq:2.2}
    \begin{aligned}
        \begin{bmatrix}
            b_{pp} & b_{pq} \\ b_{qp} & b_{qq}
        \end{bmatrix} &= 
        \begin{bmatrix}
            c & s \\ -s & c
        \end{bmatrix}\tp 
        \begin{bmatrix}
            a_{pp} & a_{pq} \\ a_{qp} & a_{qq}
        \end{bmatrix} 
        \begin{bmatrix}
            c & s \\ -s & c
        \end{bmatrix}\\
        & = 
        \begin{bmatrix}
            c^2a_{pp} - cs a_{qp} - cs a_{pq} + s^2a_{qq} &
            cs a_{pp} - s^2 a_{qp} + c^2 a_{pq} - csa_{qq}\\
            c^2a_{pq} + csa_{pp} - csa_{qq} - s^2 a_{pq} &
            csa_{qp} + s^2 a_{pp} + c^2 a_{qq} + csa_{pq}
        \end{bmatrix}.
    \end{aligned}
\end{equation}

Since $A$ is symmetric (symmetry is preserved by orthogonal transformations), we have $a_{pq} = a_{qp}$. Hence we can make further simplification to \eref{eq:2.2}
\begin{equation}
    \label{eq:2.3}
    \begin{bmatrix}
        b_{pp} & b_{pq} \\ b_{qp} & b_{qq}
    \end{bmatrix} =
    \begin{bmatrix}
        c^2 a_{pp} + s^2 a_{qq} -2csa_{pq} & cs(a_{pp} - a_{qq}) + (c^2 -s^2)a_{pq}\\
        cs(a_{pp}-a_{qq}) + (c^2 -s^2)a_{pq} & s^2 a_{pp} + c^2a_{qq} + 2csa_{pq}
    \end{bmatrix}.
\end{equation}
To achieve step 2, we require $b_{pq} = b_{qp} = 0$, which means
\begin{equation}\label{eq:2.4}
    cs(a_{pp}-a_{qq}) + (c^2 -s^2)a_{pq} = 0.
\end{equation}
This equality is essential for us to find the desired Givens rotation and the construction will be discussed in section~\ref{sec:2.3}.

If we take the Frobenius norm on the first line of \eref{eq:2.2} and set $b_{pq}=b_{qp} = 0$, since the Frobenius norm is orthogonally invariant norm (Theorem~\ref{thm:matrix-invariant-norm}), we have 
\begin{equation}\label{eq:2.8}
    b_{pp}^2 + b_{qq}^2 = a_{pp}^2 + a_{qq}^2 + 2a_{pq}^2.
\end{equation}

Using the notations in proposition~\ref{prop:1.3} and the condition $b_{pq} = b_{qp}=0$, we have 
\begin{equation}
    \label{eq:2.5}
    \begin{aligned}
        \off(B)^2 &= \|B\|_F^2 - \sum_{i=1}^n b^2_{ii} = \|A\|_F^2 - \sum_{\substack{i=1\\i\neq p,q}}\big(b^2_{ii}\big) - \big(b^2_{pp} + b^2_{qq}\big)\\
        & = \|A\|_F^2 - \sum^n_{\substack{i=1\\i\neq p,q}}(b^2_{ii}) - a_{pp}^2 - a_{qq}^2 - 2a_{pq}^2 \quad \text{Using \eref{eq:2.8}}\\
        & = \|A\|_F^2 - \sum^n_{\substack{i=1\\i\neq p,q}}(a^2_{ii}) - a_{pp}^2 - a_{qq}^2 - 2a_{pq}^2 \quad \text{By Proposition~\ref{prop:1.3}}\\
        & = \|A\|_F^2 - \sum_{i=1}^n(a^2_{ii}) - 2a_{pq}^2 = \off(A)^2 -2a_{pq}^2.
    \end{aligned}
\end{equation}

Therefore, after $k$ steps and for some choices of the index pair $(p,q)$, we always have 
\begin{equation}\notag
  \off(A\iter{k+1}) = \off(A\iter{k}) - 2a_{pq}^2,
\end{equation}
and this reduction will never terminate until all the off-diagonal entries are zero. However, since each time we introduced a zero, the previous entry we chose will not necessarily remain zero, hence the $\off(A\iter{k})$ will require an infinite number of iterations to converge to zero. In Section~\ref{sec:classical-jacobi}, we proved that even if we choose the largest $a_{pq}$ at each iteration, we still require infinite iteration for $\off(A\iter{k})$ converges to zero.

\subsection{Choices of Givens rotation matrix}\label{sec:2.3}

In this section, we provide a theory to choose $c$ and $s$ for constructing the Givens rotation matrix we need at each iteration. Firstly, if $a_{pq}=0$, there is no need to do further work since $\off(B) = \off(A)$ and we can simply choose $c = 1$ and $s = 0$ which lead to an identity matrix. Otherwise, we look at \eref{eq:2.4}
\begin{equation}
  \label{eq:sec:3.1.1_equality}
    cs(a_{pp}-a_{qq}) + (c^2 -s^2)a_{pq} = 0. \tag{3.6}
\end{equation}
Notice that $c = \cos(\theta)$ and $s = \sin(\theta)$, we have 
\begin{equation}\notag
    \label{2.9}
    \cot(2\theta) = \frac{\cos^2(\theta) - \sin^2(\theta)}{2\cos(\theta)\sin(\theta)} = \frac{c^2 - s^2 }{2cs} = \frac{a_{qq}-a_{pp}}{2a_{pq}} =:\tau.
\end{equation}
By define $t = s/c$, we have transformed \eref{eq:sec:3.1.1_equality} into 
\begin{equation}\label{eq:2.10}
    t^2 + 2\tau t - 1 = 0.
\end{equation}

This quadratic equation can be solved and we would like to choose the root of smaller magnitude.
\begin{equation}\notag
    t_{1,2} = -\tau \pm \sqrt{\tau^2 + 1},\quad \to \quad t_{\mathrm{min}} = 
    \begin{cases}
        -\tau + \sqrt{\tau^2 + 1}, & \ift \tau\geq 0 ; \\
        -\tau - \sqrt{\tau^2 + 1}, & \ift \tau < 0.
    \end{cases}
\end{equation}

Using the fact that $s^2 + c^2 = 1$ and $t_{\mathrm{min}} = s/c$, we can define $c$ and $s$ by 
\begin{equation}\label{eq:2.12}
    s = ct_{\mathrm{min}},\quad c = \frac{1}{\sqrt{1+t_{\mathrm{min}}^2}}.
\end{equation}

\begin{figure}[!tbhp]
\centering
\includegraphics[width=0.5\textwidth]{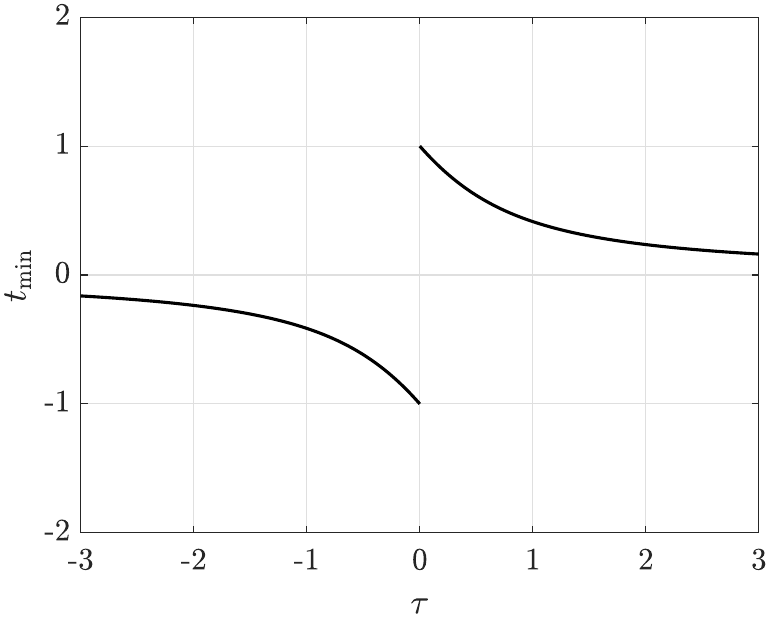}
\caption{Behavior of the solution $t_{\min}$ of \eqref{eq:2.10}.}
\label{fig:1}
\end{figure}

Here we choose the smaller roots (in magnitude) of $t$, so that we ensure $\max(t_{\mathrm{min}}) = 1$ (shown in Figure~\ref{fig:1}), hence the angle of rotation is bounded $|\theta|\leq \pi/4$ and this choice turns out to be important for the quadratic convergence as discussed by Sch{\"o}nhage~\ycite[1961]{1961-first-quadratic-convergence-Schu}.

\subsection{Implementation and Testing}
Section~\ref{sec:2.3} can be summarized in the following algorithm.
\begin{algorithm} \label{alg:jacobi-pair}
Given a symmetric matrix $A\in \R\nn$ and integers $p,q$ such that
$1\leq p < q\leq n$. This algorithm computes a pair $(c,s)$ to construct
$G(p,q,\theta)$ as presented in \eqref{eq:givens}, such that if $B =
G(p,q,\theta)\tp A G(p,q,\theta)$, then $b_{pq} = b_{qp} = 0$.
\begin{code}
if \= $a_{pq} = 0$ \\
\> $c=1,s=0$\\
else\\
\> $\tau = (a_{qq} - a_{pp})/2a_{pq}$\\
\> if \= {$\tau\geq 0$} \\
\> \> $t=-\tau + \sqrt{\tau^2 +1}$ \\
\> else \\
\> \> $t = -\tau - \sqrt{\tau^2 +1}$\\
\> end\\
\> $c = 1/\sqrt{1+t^2}, s = ct$\\
end
\end{code}
\end{algorithm}

When we do the MATLAB implementation, we need to make a minor change to lines 6 and 8. In practice, we need to change these to 
\begin{equation}
  \label{eq:jacobi-pair-preferable-way}
  t = \frac{1}{\tau  \pm  \sqrt{1 + \tau^2}},
\end{equation}
where we choose the positive sign if $\tau \geq 0$ and the minus sign if $\tau <0$. It is easy to prove these two sets of expressions are the same in exact arithmetic, however, they differ in finite precision (See Appendix~\ref{app:algorithm-jacobi-pair}) and \eqref{eq:jacobi-pair-preferable-way} is usually preferable. Then it is straightforward to implement the choice using MATLAB.

\begin{lstlisting}
function [c,s] = jacobi_pair(A,p,q)
if A(p,q) == 0
  c = 1; s = 0;
else
  tau = (A(q,q)-A(p,p))/(2*A(p,q));
  if tau >= 0
    t = 1/(tau + sqrt(1+tau*tau));
  else
    t = 1/(tau - sqrt(1+tau*tau));
  end
  c = 1/sqrt(1+t*t); s = t*c;
end
end
\end{lstlisting}

We can test this function by constructing the matrix $G(p,q,\theta)$ with $(p,q)$ given and examine the $(p,q)$ entry using the following routine,

\begin{lstlisting}
clc; clear; close all;
n = 100; format short e;
A = my_testmatrices(n);
pq = [57,64]; % choose a pair of index
[c,s] = jacobi_pair(A,pq(1),pq(2)); 
% construct the desire Givens rotation matrix
G = eye(n); G([pq(1),pq(2)],[pq(1),pq(2)]) = [c,s;-s,c];
B = G' * A * G;
disp('Before'); disp(A(pq(1),pq(2))); disp(A(pq(2),pq(1)));
disp('After');disp(B(pq(1),pq(2))); disp(B(pq(2),pq(1)));
\end{lstlisting}
This routine should give $b_{pq} = b_{qp} = 0$ where $p = 57$ and $q = 64$.
\begin{lstlisting}
Before
  -9.4140e-01
  -9.4140e-01
After
  -2.2204e-16
  -1.1102e-16
\end{lstlisting}
Before applying the Givens rotations, the entries are $-9.414\times 10^{-1}$ and after such transformation we have $b_{pq}$ and $b_{qp}$ small enough for computers to consider them as zeros. Hence we successfully provide a way of constructing Givens rotations that achieve step 2 as discussed earlier in this section and the remaining part of this chapter is to sophistically choose the index pair $(p,q)$. 

\section{Classical Jacobi Algorithm}\label{sec:classical-jacobi}

From \eref{eq:2.5}, we proved  that $\off(B)^2 = \off(A)^2 - 2a_{pq}^2$, hence it is natural to consider the choice $(p,q)$ such that $a_{pq}^2$ is maximal. This choice was first proposed by Jacobi~\cite{Jacobi-original-paper-1846} in 1846 and the future literature refer to it as the \emph{classical Jacobi algorithm}. We can summarize this choice into Algorithm~\ref{alg:classical Jacobi}.

\begin{algorithm}\label{alg:classical Jacobi}
Given a symmetric matrix $A\in\R\nn$ and a positive tolerance $tol$, this algorithm overwrites $A$ with $V\tp AV$ where $V$ is orthogonal.
\begin{code}
{$V = I_n,\text{ done\_rot = true}$}  \\
while \= {$\text{done\_rot}$}\\
\>  $\text{done\_rot = false}$ \\
\> Choose $(p,q)$ so that $|a_{pq}| = \max_{i\neq j}|a_{ij}|$ \\
\> if \= {$\abs{a_{pq}} > tol\cdot \norm{A}\sqrt{\abs{a_{pp}a_{qq}}}$} \\
\> \> $\text{done\_rot = true}$ \\
\> \> Construct Givens rotation matrix $G$ using
Algorithm~\ref{alg:jacobi-pair}\\
\> \> Update $G\tp A G \to A$ \\
\> \> Update $VG \to V$ \\
\> else \\
\> \> Remove the off-diagonal entries of $A$ \\
\> \> break\\
\> end\\
end
\end{code}
\end{algorithm}

The structure of Algorithm~\ref{alg:classical Jacobi} adapts from~\ycite[2001, Algorithm~2.1]{DHT2001AnalysisCholeskyMethod} and the stopping criterion is adapted from~\ycite[2013, Section~8.5.5]{van2013mc}, \ycite[2000, Theorem~1.1]{DV2000Approximateeigenvectorsas} and \ycite[1992, Section~1]{DV1992Jacobismethodis}.

\subsection{Linear Convergence}

Notice that since $|a_{pq}|$ is chosen to be the largest off-diagonal element, hence we have the inequality
\begin{equation}
    \label{eq:3.1}
    \off(A)^2 \leq N\cdot(2a_{pq}^2),\quad N = \frac{n(n-1)}{2}.
\end{equation}
Using \eref{eq:2.5}, we have 
\begin{equation}
    \label{eq:3.2}
    \begin{aligned}
        \off(B)^2 &= \off(A)^2 - 2a_{pq}^2 \\
        & \leq \off(A)^2 - \frac{1}{N}\off(A)^2 \\
        & = \left(1-\frac{1}{N}\right) \off(A)^2.
    \end{aligned}
\end{equation}

Denote $A^{(k)}$ as the $k$th Jacobi update of $A^{(0)} = A$, then we have the iterative bound 
\begin{equation}
    \label{eq:3.3}
    \off(A^{(k)})^2 \leq \left(1 - \frac{1}{N}\right)^k \off(A^{(0)})^2.
\end{equation}
This implies that the classical Jacobi algorithm converges \emph{linearly} to a diagonal matrix whose entries are the eigenvalues of $A$, and they may appear in any order.

\begin{remark}
    Although the term $\off(A^{(k)})$ decrease at a rate of $(1-1/N)^{1/2}$, if we consider the $N = n(n-1)/2$ steps as one `group' or `sweep' of update, then we actually can expect \emph{asymptotic} quadratic convergence~\ycite[1962]{Wil1962Notequadraticconvergence} where 
    \begin{equation*}
        \off(A^{(K + 1)})\leq c\cdot \off(A^{(K)})^2, \quad \text{for some }c,
    \end{equation*}
    where $K$ is the number of sweeps.
\end{remark}

\subsection{Implementation and Testing}\label{sec:jacobi-imple-and-testing}

Before we implement the Jacobi algorithm, we need two functions (i) the function that calculates the Frobenius norm of the off-diagonal entries and (ii) the function that finds the index of the maximum entry in modulus.

(i) can be done by carefully selecting the diagonal entries and set to zero, and then calculating the Frobenius norm of the new matrix.
\begin{lstlisting}
function offA = off(A)
n = length(A); % get the dimension of the matrix A
A(1:n+1:n*n) = 0; % set the diagonal entries to zero
offA = norm(A,"fro"); % calculate the norm
end
\end{lstlisting}
The third line of the code will eliminate the diagonal entries, then the fourth line will give us the desired output, $\off(A)$.

(ii) can be implemented using the following function:
\begin{lstlisting}
function [p,q] = maxoff(A)
n = length(A); % dimension of matrix A
A(1:n+1:n*n) = 0; A = abs(A); % clear the diagonal entries
[val, idx1] = max(A);
[~, q] = max(val);
p = idx1(q);
end
\end{lstlisting}

At the fourth line, I use \inline{[val, idx1] = max(A)} to find the index, \inline{idx1}, of the maximum entries of each column and I store the values in \inline{val}. Then we can apply \inline{max()} again to find the index of the maximum entry of \inline{val}, denoted as \inline{q}. Finally, we locate the column number of the maximum entry by calling \inline{idx1(q)} as shown in the sixth line.

Making use of these two functions, we can implement the classical Jacobi algorithm in \mat. The function \inline{jacobi_classical} will take two inputs (i) the symmetric matrix $A^{(0)}\in\R\nn$ and (ii) the tolerance $tol$. I choose to output two matrices, $A$ and $V$ such that $\norm{A\iter{0} V - VA} \lesssim nu_d\norm{A\iter{0}}$ and the number of iterations required. 

Notice that if we would like to update $A^{(k)}$ by $G_k\tp A^{(k)}G_k$, since $A^{(k)},G_k\in\R\nn$, the cost would be $\mathcal O(n^3)$ flops. However, we can utilize proposition~\ref{prop:1.3}, namely, we only need to update the $p$ and $q$th rows and columns of $A^{(k)}$, and this version of the update will only take $O(n)$ flops. We can see this reduction in the following code 

\begin{lstlisting}
[c,s] = jacobi_pair(A,p,q); G = [c s; -s c];
A([p q],:) = G'*A([p q],:); 
A(:,[p q]) = A(:,[p q])*G;
\end{lstlisting}
From the code above, we get both the updated matrix $A^{(k+1)}$ and the orthogonal matrix $G_k$. Line 2 requires a matrix multiplication between $\R^{2\times 2}\times \R^{2\times n}$ which needs $\mathcal O(n)$ flops, and the same argument applies to line 3. Hence when updating the matrix $A^{(k)}$, we only require $\mathcal O(n)$ flops. Assemble all above, we have  

\begin{lstlisting}
function [V,A,counter] = jacobi_classical(A,tol)
counter = 0; n = length(A); V = eye(n); done_rot = true;
tol1 = tol * norm(A);
while done_rot
  if isint(counter/(n*n)), A = (A + A')/2; end % maintain symmetry
  done_rot = false; [p,q] = maxoff(A);
  if abs(A(p,q)) >  tol1 * sqrt(abs(A(p,p) * A(q,q)))
    counter = counter + 1; done_rot = true;
    [c,s] = jacobi_pair(A,p,q);
    J = [c,s;-s,c]; 
    A([p,q],:) = J'*A([p,q],:);
    A(:,[p,q]) = A(:,[p,q]) * J;
    V(:,[p,q]) = V(:,[p,q]) * J;
  else
    A = diag(diag(A)); % output diagonal matrix
    break;
  end
end
end
\end{lstlisting}

After several iterations, although in exact arithmetic, the matrix $A$ should always be symmetric, in floating-point arithmetic, we need to maintain its symmetry as shown in line 5. We can then test our code by the following routine 
\begin{lstlisting}
clc;clear; format short e
n = 1e2; A = randn(n); A = A + A'; tol = 2^(-53);
[V,D,iter] = jacobi_classical(A,tol);
disp('norm(AV - VD)'); disp(norm(A *V - V * D))
\end{lstlisting}

Here we should expect the norm $\norm{AV - VD}$ is about $nu_d \norm{A}$ which is $3.0224\times 10^{-13}$ in our test.

\begin{lstlisting}
norm(AV - VD)
  4.3553e-13
\end{lstlisting}

Hence the algorithm works as expected.

\section{Cyclic-by-row Jacobi Algorithm}\label{sec:cyclic-jacobi}
The classical Jacobi method requires at each step the searching of $n(n-1)$ for one maximum entry in modulus. For $n$ is large, this can be extremely expensive. It will be better if we fixed a sequence of choice $(p,q)$, here in cyclic order shown in the following scheme suggested by~\ycite[1953]{Gre1953Computingeigenvalueseigenvectors} and~\ycite[1960, Section~1.2]{FH1960cyclicJacobimethod}

\begin{equation}\notag
    \begin{aligned}
        (p_0,q_0) &= (1,2), \\
        (p_{k+1},q_{k+1}) & = 
        \begin{cases}
            (p_k,q_{k}+1), & \ift p_k < n-1,q_k < n,\\
            (p_k+1,p_k+2), & \ift p_k < n-1, q_k = n,\\
            (1,2),& \ift p_k = n-1,q_k = n.
        \end{cases}
    \end{aligned}
\end{equation}

We keep using $(p,q)$ in this fashion until we meet the required tolerance and this procedure can be described in Algorithm~\ref{alg:cyclic-jacobi}.

\begin{algorithm}\label{alg:cyclic-jacobi}
Given a symmetric matrix $A\in\R\nn$ and a positive tolerance $tol$, this
algorithm overwrites $A$ with $V\tp AV$ where $V$ is orthogonal.
\begin{code}
{$V = I_n,\text{ done\_rot = true}$}  \\
while \= {$\text{done\_rot}$}\\
\> $\text{done\_rot = false}$\\
\> for \= {$p = 1,\dots, n-1$}  \\
\> \> for \= {$ q = p + 1,\dots, n$} \\
\> \> \> if \= {$\abs{a_{pq}} > tol\cdot \sqrt{\abs{a_{pp}a_{qq}}}$} \\
\> \> \> \> $\text{done\_rot = true}$ \\
\> \> \> \> Construct Givens rotation matrix $G$ using
Algorithm~\ref{alg:jacobi-pair} \\
\> \> \> \> Update $G\tp A G \to A$ \\
\> \> \> \> Update $VG \to V$ \\
\> \> \> else  \\
\> \> \> \>  Remove the off-diagonal entries of $A$ \\
\> \> \> \>  continue \\
\> \> \> end \\
\> \> end \\
\> end \\ 
end 
\end{code}
\end{algorithm}

Similarly, the structure of Algorithm~\ref{alg:cyclic-jacobi} adapts
from~\ycite[2001, Algorithm~2.1]{DHT2001AnalysisCholeskyMethod} and the
stopping criterion is adapted from~\ycite[2013,
Section~8.5.5]{van2013mc}, \ycite[2000,
Theorem~1.1]{DV2000Approximateeigenvectorsas} and \ycite[1992,
Section~1]{DV1992Jacobismethodis}.

\subsection{Implementation and Testing}
The MATLAB code is similar to the classical Jacobi algorithm, and the only difference is the way of choosing $(p,q)$:

\begin{lstlisting}
function [V,A,iter] = jacobi_cyclic(A,tol,maxiter)
n = length(A); V = eye(n); iter = 0; done_rot = true;
while done_rot && iter < maxiter
  done_rot = false;
  for p = 1:n-1
    for q = p+1:n
      if abs(A(p,q)) > tol * sqrt(abs(A(p,p)*A(q,q)))
        done_rot = true;
        [c,s] = jacobi_pair(A,p,q);
        J = [c,s;-s,c];
        A([p,q],:) = J'*A([p,q],:);
        A(:,[p,q]) = A(:,[p,q]) * J;
        V(:,[p,q]) = V(:,[p,q]) * J;
      end
    end
  end
  if done_rot 
    A = (A + A')/2; iter = iter + 1; 
  else
    A = diag(diag(A));
    return;
  end
end
end
\end{lstlisting}

Notice that, since the cyclic-by-row Jacobi algorithm does not reduce $\off(A\iter{k})$ in the most optimal way, it may require much more iterations than the classical Jacobi algorithm. Hence we restrict the maximum number of iterations by \inline{maxiter} at line 3.

Using the same routine as presented in Section~\ref{sec:jacobi-imple-and-testing}, we have 
\begin{lstlisting}
norm(AV - VD)
   5.9258e-13
\end{lstlisting}
The required tolerance is $nu_d\norm{A} \approx 3.0260\times 10^{-13}$ and
therefore we indeed have a spectral decomposition of $A$ at double precision.

\section{Comparision}

In the previous two sections, we discussed both the classical Jacobi algorithm and the cyclic-by-row Jacobi algorithm. However, we need to choose between them. In this section, we will examine the performance concerning both the dimension $n$ and the condition number $\kappa(A)$. We can generate Figure~\ref{fig:classical-cyclic-compare} using code in Appendix~\ref{app:code-for-fig1}. Notice that, in the code, we call a function call \inline{jacobi_eig} and this is the assembled version of the Algorithm~\ref{alg:classical Jacobi} and~\ref{alg:cyclic-jacobi} shown in Appendix~\ref{app:code-full-jacobi}.

\begin{figure}[!tbhp]
\centering
\includegraphics[width=0.85\textwidth]{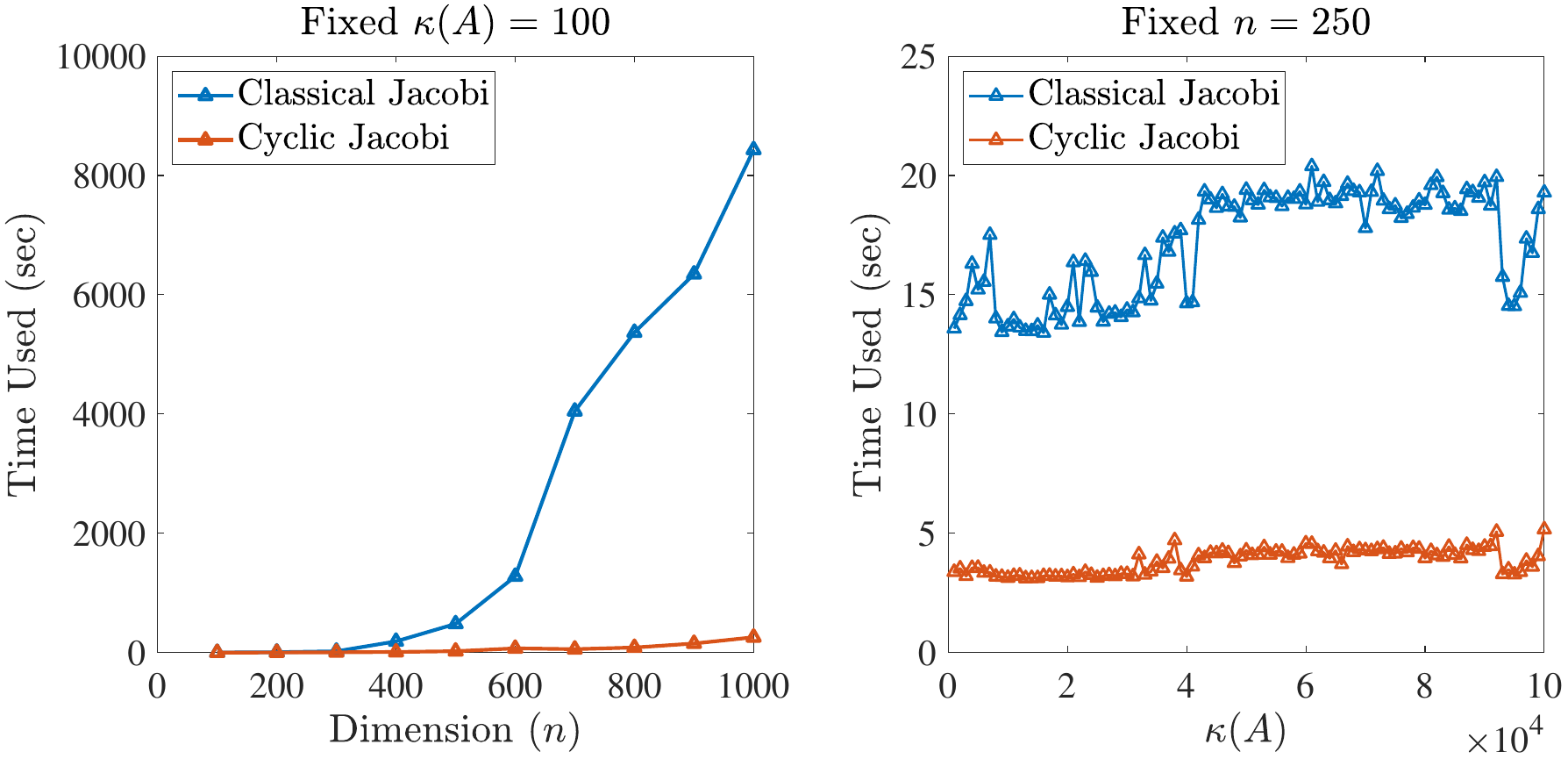}
\caption[The time of applying both the classical and the cyclic-by-row Jacobi algorithm on matrix $A\in\R\nn$ with respect to both the dimension $(n)$ and the condition number $\kappa(A)$.]{The time of applying both the codes \inline{jacobi_classical} and \inline{jacobi_cyclic} on a matrix $A\in\R\nn$ with respect to both dimension $(n)$ and the condition number $\kappa(A)$. The left figure fixes the condition number $\kappa(A) = 100$ and the right figure fixes the dimension $n = 250$.}
\label{fig:classical-cyclic-compare}
\end{figure}

From the right figure, we can see the condition number does not affect the
time too much. However, from the left figure, we observe that the downside
of searching for the maximum off-diagonal entry is significant for large
$n$. The cyclic-by-row Jacobi algorithm uses about $1/30$ of the time used
by the classical Jacobi algorithm for $n = 1000$. Therefore, although the
classical Jacobi reduces the $\off(A)$ in the most optimal way, in
practice, we should always stick with the cyclic-by-row Jacobi algorithm.


\chapter{Orthogonalization}\label{chap:orthogonalisation}

\section{Introduction}\label{sec.approx-eig}

From Chapter~\ref{chap:jacobi_algorithm}, given a symmetric matrix
$A\in\R\nn$, we can use the Jacobi algorithm to find an spectral
decomposition $A = Q\varLambda Q\tp$. However, for large $n$, this
procedure can still be time-consuming and is not generally as fast as the
symmetric QR Algorithm~\ycite[2013, Section~8.5]{van2013mc}. However, the
Jacobi algorithm can exploit the situation that $A$ is almost diagonal
($\off(A)$ is small). In 2000, Drma{\v{c}} and Veseli\'c proposed a way of
speeding up this procedure using approximate eigenvectors as
preconditioner~\ycite[2000, Section~1]{DV2000Approximateeigenvectorsas}.
They proposed the following 
\begin{enumerate}
\item Given a symmetric matrix $A\in\R\nn$ and its approximate orthogonal eigenvector matrix $Q_{\text{app}}$.
\item Computing $A' = Q_{\text{app}}\tp A Q_{\text{app}}$.
\item Diagonalizing $A'$.
\end{enumerate}
Suppose $Q_{\text{app}}$ is orthogonal in double precision and step 2 and 3 are carried out at double precision, then the eigenvalues we computed will have a small relative error and more efficient than the pure Jacobi method~\ycite[2000, Section~2]{DV2000Approximateeigenvectorsas}. In this section, we will provide a way of constructing this approximate vector matrix:
\begin{enumerate}
\item Compute the spectral decomposition at single precision such that 
\begin{equation}\notag
  A = Q_\ell D_\ell Q_\ell\tp, \quad \norm{AQ_\ell - Q_\ell D_\ell} \lesssim nu_\ell \norm{A},\quad \norm{Q_\ell \tp Q_\ell - I} \lesssim nu_\ell,
\end{equation}
where $n$ is the dimension of the real symmetric matrix $A$ and $u_\ell$ is the unite roundoff at single precision. This can be achieved using \inline{eig} function in MATLAB via the following routine 
\begin{lstlisting}
[Q_low,D] = eig(single(A));
Q_low = double(Q_low);
\end{lstlisting}
\item Orthogonalize the matrix $Q_\ell$ to $Q_d$ such that it is orthogonal at double precision,
\begin{equation}\notag
  \norm{Q_d\tp Q_d - I} \lesssim nu_d.
\end{equation}
\end{enumerate}

In this chapter, we will first review two methods to orthogonalize $Q_\ell$: the QR factorization and the polar decomposition. We aim to find the orthogonal (at double precision) matrix $Q_d$ that minimizes the norm $\norm{Q_d - U}$ where $U$ is the exact matrix of the eigenvectors of $A$. However, in practice, we have no access to the matrix $U$, therefore, we instead find the orthogonal matrix $Q_d$ that minimizes $\norm{Q_d - Q_\ell}$. This problem can be transformed to: Given $Q_\ell$, find the orthogonal matrix $Q_d$ such that it satisfies
\begin{equation}\notag
  \min{\norm{Q_\ell - Q_d}},\quad \text{under constraint: $\norm{Q_d\tp Q_d  = I}$.}
\end{equation}

Finally, the code in MATLAB will be presented and we will access them to decide which to use in practice.

\section{The Householder QR Factorization}\label{sec:Householder-QR}
The QR factorization of $A\in\R\nn$ is a factorization
\begin{equation}\notag
  A = QR
\end{equation}
where $Q\in\R\nn$ is orthogonal and $R\in\R\nn$ is upper triangular. Suppose we have a QR factorization at double precision of $Q_\ell$,
\begin{equation}\label{eq.2.2}
  Q_\ell = Q_d R,
\end{equation}
where $Q_d$ is orthogonal at double precision, and this is already our desired preconditioner. The QR factorization can be achieved using Householder transformations in Section~\ref{sec:Householder}. 

\subsection{Theory of Householder QR Factorization}\label{subsec.houseQR}

The idea can be illustrated using a small matrix $A\in\R^{4\times 4}$:
\begin{equation}\notag
  \begin{aligned}
    A &=
        \begin{bmatrix}
          \times & \times & \times & \times\\
          \times & \times & \times & \times\\
          \times & \times & \times & \times\\
          \times & \times & \times & \times
        \end{bmatrix}
        \stackrel{P_1 \in\R^{4\times 4}}{\longrightarrow}
        \begin{bmatrix}
          \times &\mid \times & \times & \times\\\hline
          0 &\mid \times & \times & \times\\
          0 & \mid\times & \times & \times\\
          0 & \mid\times & \times & \times\\
        \end{bmatrix}
        \stackrel{P_2 \in\R^{3\times 3}}{\longrightarrow}
        \begin{bmatrix}
          \times & \times & \times & \times\\
          0 & \times & \times & \times\\ \hline
          0 & 0 &\mid \times & \times\\
          0 & 0 & \mid \times & \times\\
        \end{bmatrix}\\
      & \stackrel{P_3 \in\R^{2\times 2}}{\longrightarrow}
        \begin{bmatrix}
          \times & \times & \times & \times\\
          0 & \times & \times & \times\\
          0 & 0 &\times & \times\\
          0 & 0 & 0 & \times\\
        \end{bmatrix} = R,
  \end{aligned}
\end{equation}
where $P_k$, $k = 1,2,3$, are the Householder transformations. Notice that $P_2\in\R^{3\times 3}$ only applies to the bottom right corner of the matrix, and this can be done by embedding $P_2$ inside a larger $4\times 4$ matrix via
\begin{equation}\notag
  \wt P_2 =
  \begin{bmatrix}
    1 & 0 \\
    0 & P_2
  \end{bmatrix}
\end{equation}

Hence from the example above, we can generalize the Householder QR factorization for $A\in\R\nn$. It will produce a sequence of matrices $\{A_k\}_{k=1}^{n}$ where $A_1$ is the input $A$ and $A_n$ is the upper triangular matrix $R$. Notice, there are $k-1$ stages where the $k$th stage transform $A_k$ to $A_{k+1}$.

At the $k$th stage of the Householder QR factorization, we can carefully partition the matrix $A_k$ into the following form

\begin{equation}\notag
  A_{k} =
  \begin{bmatrix}
    R_{k-1} & w_k & C_k \\
    0 & x_k & B_k
  \end{bmatrix},
  \quad R_{k-1}\in\R^{(k-1)\times (k-1)}, x_k\in\R^{n-k+1},
\end{equation}
where $R_{k-1}$ is upper triangular (achieved at the ($k-1$)th stage) and $x_k$ is the vector that we should focused on the $k$th stage. Choose a Householder transformation $P_k$ such that $P_kx_k = \sign{x_{k,1}}\norm{x_k}e_1$, where $e_1$ is the first column of $I_{n-k+1}$ and $x_{k,1}$ is the first entry of the vector $x_k$. We then embed matrix $P_k$ into a larger matrix $\wt P_k\in\R\nn$ via
\begin{equation}\label{eq.qrQ}
  \wt P_k =
  \begin{bmatrix}
    I_{k-1} & 0 \\
    0 & P_k
  \end{bmatrix},\quad P_k \in \R^{(n-k+1)\times(n-k+1)}.
\end{equation}
Then let $A_{k+1} = \wt P_k A_k$, we have 
\begin{equation}\label{eq.qrstep}
  A_{k+1} = 
  \begin{bmatrix}
    R_{k-1} & w_k & C_k \\
    0 & P_k x_k & P_k B_k
  \end{bmatrix}
  = 
  \begin{bmatrix}
    R_{k-1} & w_k & C_k \\
    0 & \sign{x_{k,1}}\norm{x_k}e_1 & P_k B_k
  \end{bmatrix},
\end{equation}
which is closer to the upper triangular form. After $n-1$ steps, the matrix $A_{n}$ will be upper triangular and we denote as $R$ and we obtain the QR factorization of $A$
\begin{equation}\notag
  A_{n-1} = R = \wt P_{n-1}\wt P_{n-2} \cdots \wt P_1 A =: Q\tp A.
\end{equation}
Since $\wt P_k$ are composed of $P_k$, Householder matrices, and the identity matrix. Hence it is obvious that $\wt P_k$ are also symmetric and orthogonal. Then we can construct $Q$ via $\wt P_1 \cdots\wt P_{n-1}$. This procedure can be summarized in the Algorithm~\ref{alg:householderQR}.

\begin{algorithm}
\label{alg:householderQR}
Given $A\in\R\nn$, this algorithm computes an orthogonal matrix $Q$ and an upper triangular matrix $R$ such that $A = QR$.
\begin{code}
$Q = I_n$\\
for \= {$k=1:n-1$}\\
\> {$[v,\beta] = \mathtt{house}(A(k:n,k))$}\qquad \%{Algorithm~\ref{alg.house}}\\
\> {$A(k:n,k:n) = (I_{n-k+1}-\beta vv\tp )A(k:n,k:n)$} \\
\> {$Q(1:n,k:n) = Q(1:n,k:n)(I_{n-k+1}-\beta vv\tp)$} \\
end
\end{code}
\end{algorithm}

Line $4$ is directly adapted from~\eqref{eq.qrstep}. Line $5$ can be viewed by partitioning the orthogonal matrix $Q_k$ into four parts and multiplying the matrix $\wt P_k$ structured as~\eqref{eq.qrQ},
\begin{equation}\label{eq:update-Q}
  Q_{k+1} = Q_k\wt P_k = 
  \begin{bmatrix}
    Q_1 & Q_2 \\
    Q_3 & Q_4
  \end{bmatrix}
  \begin{bmatrix}
    I_{k-1} & 0 \\
    0 & P_k
  \end{bmatrix} = 
  \begin{bmatrix}
    Q_1 & Q_2 P_k \\
    Q_3 & Q_4 P_k
  \end{bmatrix}.
\end{equation}
Hence, at each step of the algorithm, we only update the final $n-k+1$ columns of the matrix $Q_k$. Also, focused on the line 4, if we construct $P$ and do a matrix-matrix multiplication between $D_k,P_k\in\R^{(n-k+1)\times (n-k+1)}$ where $D_k = [x_k,B_k]$, then the cost will be about $2(n-k+1)^3$. If we utilize the components we computed from Algorithm~\ref{alg.house}, the Householder vector $v$ and the coefficient $\beta$, we can reduce the computational cost by doing the matrix-vector products instead. We can rewrite $P_k D_k$ using $v$ and $\beta$
\begin{equation}\notag
  P_k D_k = (I - \beta vv\tp) D_k = D_k - \beta v v\tp D_k = D_k - (\beta v)\cdot(v\tp D_k).
\end{equation}
By this procedure, we transform a matrix-matrix multiplication to
\begin{enumerate}
\item Vector-matrix multiplication: $v\tp D_k$ requires $2(n-k+1)^2$ flops.
\item Scalar-vector inner product: $\beta v$ requires $O(n-k+1)$ flops.
\item Vector-vector outer product: $(\beta v)\cdot (v\tp D_k)$ requires $(n-k+1)^2$ flops.
\item Matrix-matrix subtraction: $D_k - (\beta v)\cdot(v\tp D_k)$ requires $O(n-k+1)^2$ flops.
\end{enumerate}
Therefore, the overall cost will be $4(n-k+1)^2$ flops. Compare with the matrix-matrix multiplication which requires $O((n-k+1)^3)$ flops, this approach is much more efficient.

\subsection{Implementation and Testing}\label{sec.qrtesting}
Based on these analyses in Section~\ref{subsec.houseQR}, we can implement it into \mat.

\begin{lstlisting}
function [Q,A] = myqr(A)
[m,n] = size(A); Q = eye(n);
if m ~= n, error('Input should be a square matrix.'), end
for k = 1:n-1
[v,b] = house(A(k:n,k)); % compute the components of P
A(k:n,k:n) = A(k:n,k:n) - (b*v)*(v'*A(k:n,k:n)); % update A
Q(:,k:n) = Q(:,k:n) - (Q(:,k:n)*v)*(b*v'); % update Q
end
\end{lstlisting}

Recall from the last section, we discussed that at the $k$th step, updating $A$ requires about $4(n-k+1)^2$ flops. Similarly, we can update $Q$ based on~\eqref{eq:update-Q} which only involves $(n-k+1)$ columns of $Q$. Therefore, line 7 requires $4(n-k+1)n$ flops at the $k$th step. We can then compute the theoretical overall cost,
\begin{equation}\notag
  \text{Cost} = \sum_{k = 1}^{n-1} 4(n-k+1)^2 + 4(n-k+1)n = O(10n^3/3).
\end{equation}

To test our code, we are focused on two quantities, $\norm{Q\tp Q - I}$ and $\norm{QR - A}$. The former evaluates whether our computed $Q$ is orthogonal at double precision and the latter evaluates whether we have a QR factorization. In theory, we should have 
\begin{equation}\notag
  \norm{Q\tp Q - I} \lesssim nu_d, \quad \norm{QR - A}/\norm{A} \lesssim nu_d.
\end{equation}
We can use the following MATLAB routine, notice that we add the accuracy of the MATLAB built-in function \inline{qr} as a reference.
\begin{lstlisting}
clc; clear; format short e; 
n = 1e2; ud = 2^(-53); A = randn(n); 
[Q,R] = myqr(A); [Q1,R1] = qr(A); % qr factorization
orthg = norm(Q'*Q - eye(n)); % check orthogonality
qralg = norm(Q*R - A)/norm(A); % check if we have a qr factorization
fprintf('Orthogonal? %d \nMy QR accuracy is %d\n',orthg, qralg);
fprintf('MATLAB qr function accuracy is %d \n', norm(Q1*R1 - A)/norm(A));
fprintf('n*(machine precision) = %d\n', ud*n);
\end{lstlisting}
And we have the output:
\begin{lstlisting}
Orthogonal? 3.567685e-15 
My QR accuracy is 1.493857e-15
MATLAB qr function accuracy is 1.026458e-15 
n*(machine precision) = 1.110223e-14
\end{lstlisting}

The first and second outputs are all smaller than $nu_d \approx 10^{-13}$ therefore our code works well. In addition, the second and the third outputs are the accuracy of functions \inline{myqr} and \inline{qr} and we can see these are about the same, therefore our code achieves similar accuracy as the MATLAB built-in function.

\begin{figure}[!tbhp]
\centering 
\includegraphics[width=0.5\textwidth]{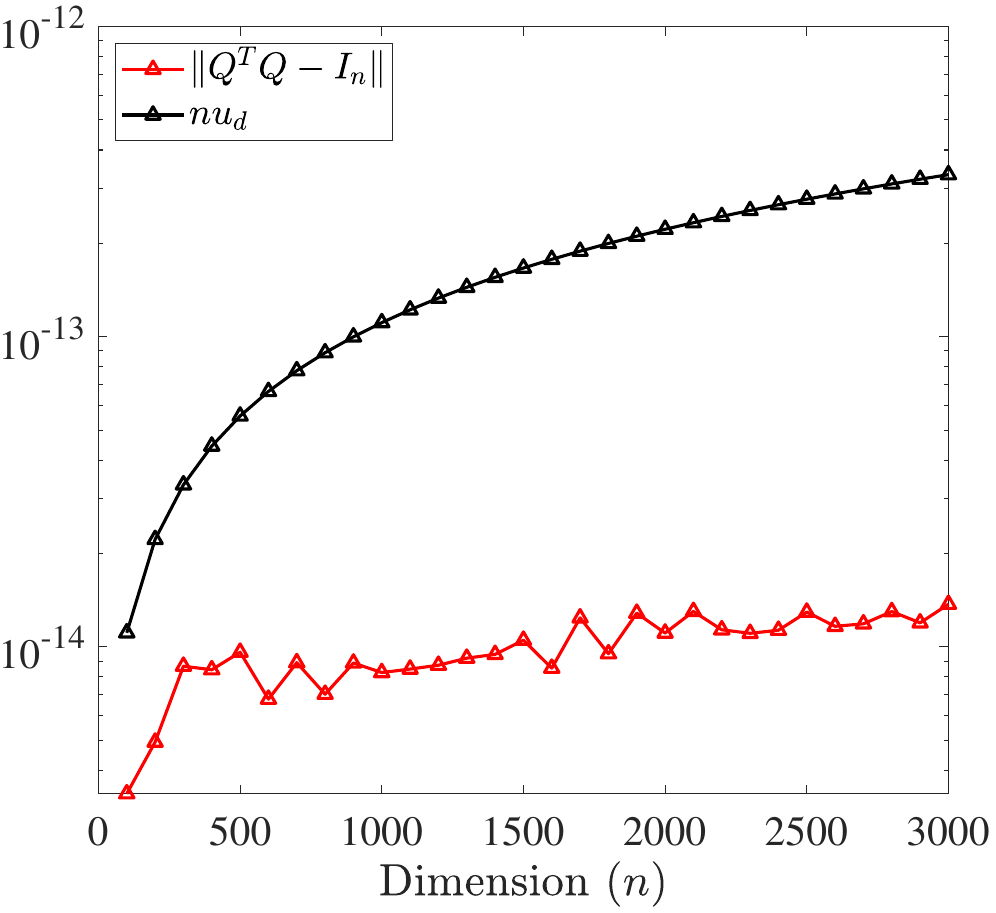}
\caption[Behavior of $\norm{Q\tp Q - I_n}$ with $n$ increases and $\kappa(A)$ fixes using the Householder QR factorization.]{Behavior of $\norm{Q\tp Q - I_n}$ with $n$ increases from $100$ to $3000$ with step size $100$ and $\kappa(A) = 100$ using my own QR code \inline{myqr}. The red line shows the computed results and the black line shows the reference line $nu_d$. The code to regenerate this graph can be found in Appendix~\ref{app:fig:orthog}.}
\label{fig:orthog}
\end{figure}

Moreover, I am interested in how $\norm{Q\tp Q - I_n}$ behaved as $n$ increases. From Figure~\ref{fig:orthog}, we can see $\norm{Q\tp Q - I_n}$ does not increase too much when $n$ increases and is well bounded by the reference line $nu_d$.

By now, we have a way to orthogonalize a matrix. However, this approach does not utilize the property that our input matrix $Q_\ell$ is an almost orthogonal matrix. Namely, applying \inline{myqr} to $Q_\ell$ has no difference from applying \inline{myqr} to a general matrix. In the next section, we will introduce a way to orthogonalize a matrix which exploits the fact that the input is almost orthogonal.

\section{The Polar Decomposition}\label{sec:polar-decomposition}

In complex analysis, it is known that for any $\alpha\in\C$, we can write $\alpha$ in polar form, namely $\alpha = r\eu^{\im \theta}$. The polar decomposition is its matrix analogue. The polar decomposition can be derived from the singular value decomposition discussed in Section~\ref{sec:svd}. Conventionally, we will present our analysis in complex but this can be restricted to real. All the norm $\gnorm{\cdot}$ discussed in this section will be the orthogonally invariant norm.

\begin{theorem}
[Polar decomposition]
\label{thm.polar-decompo}
If $A\in\C\nn$, then there exists a unitary matrix $U\in\C\nn$ and a unique Hermitian positive semidefinite matrix $H\in\C\nn$ such that $A = UH$. We call $U$ the unitary factor.
\end{theorem}

\begin{proof}
By Theorem~\ref{thm.svd}, suppose $A\in\C\nn$ and $\rank(A) = r$, then the SVD of $A$ can be written as 
\begin{equation}\label{eq.polar-decomp-factors}
  A = P\varSigma_r V\ctp = 
  \underbrace{PV\ctp}_{U} 
  \underbrace{
  V \Sigma_r
  V\ctp}_{H}
  \equiv UH,
\end{equation}
where $\Sigma_r = \diag(\sigma_1,\dots,\sigma_r,0,\dots,0)\in\R\nn$.
Unitarity of $U$ can be easily proved by its definition 
\begin{equation}\notag
  UU\ctp = PV\ctp VP\ctp = I_n,\quad U\ctp U = VP\ctp P V\ctp = I_n.
\end{equation}
$H$ is a positive semidefinite due to its eigenvalues are the singular values of $A$ which are all non-negative. By forming $A\ctp A$, we have $A\ctp A = H\ctp U\ctp U H = H\ctp H = H^2$, therefore $H$ can be uniquely determined by $(A\ctp A)^{1/2}$~\ycite[2008, Theorem~1.29]{higham08fm}.
\end{proof}

Notice that, if $A$ is full rank, then $H$ is clearly non-singular by construction, therefore $U$ can be uniquely determined by $U = AH\inv$. Similar to the QR factorization, we have a decomposition of $A$ where one of the factors is unitary, and we expect that we can use this unitary factor as the preconditioner. The next section gives a special property which shows the unitary factor $U$ of the polar decomposition of $A$ is the one that minimizes the distance $\gnorm{A-U}$.

\subsection{Best Approximation Property}

The unitary factor $U$ of the polar decomposition of $A = UH$ satisfies the best approximation property.

\begin{theorem}[{\ycite[1955, Theorem~1]{FH1955Somemetricinequalities},\ycite[2008, Theorem~8.4]{higham08fm}}]
\label{thm.best_approx}
Let $A\in\C\nn$ and $A= UH$ is its polar decomposition. Then 
\begin{equation}\notag
  \gnorm{A-U} = \min\left\{\gnorm{A-Q}\mid Q\ctp Q = I_n\right\}
\end{equation}
holds for every unitarily invariant norm.
\end{theorem}

To prove Theorem~\ref{thm.best_approx}, we need the following two results:

\begin{lemma}[{\ycite[1960, Lemma~4]{Mir1960Symmetricgaugefunctions}}]
\label{lemma.3.2}
Let $X, Y$ be Hermitian matrices, and denote by 
\begin{equation}\notag
  \xi_1\geq \cdots \geq \xi_n,\quad 
  \eta_1 \geq \cdots \geq \eta_n,\quad 
  \zeta_1\geq \cdots \geq \zeta_n
\end{equation}
be the eigenvalues of $X, Y$ and $X-Y$ respectively. Then if we denote by $(\xi-\eta)_1,\dots,(\xi-\eta)_n$ the numbers $\xi_1 -\eta_1,\dots,\xi_n-\eta_n$ arranged in non-ascending order of magnitude, then we have
\begin{equation}\notag
  \sum_{i=1}^k (\xi-\eta)_i \leq \sum_{i=1}^k(\zeta_i), \quad 1\leq k\leq n.
\end{equation}
\end{lemma}

\begin{proof}
See \ycite[1955, Theorem~2]{FH1955Somemetricinequalities}.
\end{proof}

\begin{theorem}[{\ycite[1960, Theorem~5]{Mir1960Symmetricgaugefunctions}, \ycite[2008, Theorem~B.3]{higham08fm}}]
\label{thm.3.3}
Let $\rho_1\geq \cdots\geq\rho_n$ and $\sigma_1\geq \cdots\geq \sigma_n$ be the singular values of the complex matrices $A$ and $B$ respectively. Then 
\begin{equation}\notag
  \gnorm{A-B}\geq \gnorm{\Sigma_A - \Sigma_B},
\end{equation}
where $\Sigma_A$ and $\Sigma_B$ are diagonal matrices with entries the singular values of $A$ and $B$ arranged in non-ascending order respectively.
\end{theorem}

\begin{proof}
We have the fact by H. Wielandt mentioned in~\ycite[1955, p.~113]{FH1955Somemetricinequalities}: For any matrix $M$ of order $n\times n$, the Hermitian matrix 
\begin{equation}\notag
  \wt M = 
  \begin{bmatrix}
    0 & M\\
    M\ctp & 0
  \end{bmatrix}
\end{equation}
of order $2n\times 2n$. The eigenvalues of $\wt M$ are precisely the singular values of $M$ and their negatives. We denote $\tau_1 \geq \cdots \geq \tau_n$ the singular values of $A-B$ and we have 
\begin{equation}\notag
  \begin{bmatrix}
    0 & A \\
    A\ctp & 0
  \end{bmatrix}
  -
  \begin{bmatrix}
    0 & B\\
    B\ctp & 0
  \end{bmatrix}
  =
  \begin{bmatrix}
    0 & A-B\\
    A\ctp - B\ctp & 0
  \end{bmatrix}
\end{equation}
and the eigenvalues of these three Hermitian matrices are 
\begin{equation}\notag
  \begin{aligned}
    \rho_1\geq\cdots\geq\rho_n &\geq -\rho_n\geq \cdots\geq -\rho_1, \\
    \sigma_1\geq\cdots\geq\sigma_n&\geq -\sigma_n\geq \cdots\geq -\sigma_1, \\
    \tau_1\geq\cdots\geq\tau_n&\geq -\tau_n\geq \cdots\geq -\tau_1.
  \end{aligned}
\end{equation}
By Lemma~\ref{lemma.3.2}, if we denote $(\rho-\sigma)_1,\dots,(\rho-\sigma)_n$ the numbers $\rho_1 -\sigma_1,\dots,\rho_n-\sigma_n$ arranged in non-ascending order of magnitude, then we have 
\begin{equation}\notag
  \sum_{i=1}^k \tau_i \geq \sum_{i=1}^k (\rho - \sigma)_i,\quad 1\leq k \leq n.
\end{equation}
Here we focus on the proof for $2$-norm, to see the proof considering all the orthogonally invariant norms, refer to~\ycite[1960, Theorem~1]{Mir1960Symmetricgaugefunctions}. By Corollary~\ref{col.2norm}
\begin{equation}\notag
  \begin{aligned}
    \norm{A-B} &  = \norm{\diag\{\tau_1,\dots,\tau_n\}} = \tau_1 \\
               & \geq (\rho - \sigma)_1 = \max_k{\abs{\rho_k - \sigma_k}}\\
               & = \norm{\diag\{{\rho_1-\sigma_1},\dots,{\rho_n-\sigma_n}\}} \\
               & = \norm{\Sigma_A - \Sigma_B} \qedhere
  \end{aligned} 
\end{equation}
\end{proof}

\begin{proof}
[Proof of Theorem~\ref{thm.best_approx}]
For any unitary matrix $Q$, it has singular values all equal to $1$. Hence by Theorem~\ref{thm.3.3}, if $A$ has the singular value decomposition $A=P\Sigma_A V\ctp$, then $$\gnorm{A-Q}\geq\gnorm{\Sigma_A - I_n}.$$
Remain to prove that $\gnorm{A-U} = \gnorm{\Sigma_A - I_n}$. From Theorem~\ref{thm.polar-decompo}, the Hermitian part $H$ can be written as $V\Sigma_A V\ctp$, hence 
\begin{equation}\notag
  \gnorm{A-U} = \gnorm{U(H-I_n)} = \gnorm{H-I_n} = \gnorm{V \Sigma_A V\ctp - I_n} = \gnorm{\Sigma_A - I_n},
\end{equation}
which completes the proof.
\end{proof}

This property indicates that the unitary factor of the input $Q_\ell$ will be the desired orthogonal matrix and then we can use it as the preconditioner. The rest of this section will focus on how to compute this factor.

\subsection{SVD approach}\label{sec:svdapproach}

The most obvious approach will be the singular value decomposition approach.

For $A\in\C\nn$, by Theorem~\ref{thm.polar-decompo}, we can use the singular value decomposition to compute the polar decomposition using the following procedure:
\begin{enumerate}
\item Compute the singular value decomposition of $A\in\C\nn$ where $A = P\Sigma V\ctp$.
\item Form $U$ via~\eqref{eq.polar-decomp-factors}, $U = P V\ctp$.
\end{enumerate}
The computational cost of computing the singular value decomposition is approximately $21n^3$ flops \ycite[2013, Figure~8.6.1]{van2013mc} with $O(n^3)$ flops to form $U$. This is good generally, but it does not utilize the property that, in our scenario, $A$ is almost orthogonal ($\norm{A\ctp A - I_n}$ is small). The next two sections will discuss a more sophisticated way of computing the unitary factor.

\subsection{Newton's Method Approach}

Consider the Newton iteration~\ycite[2008, Section~8.3]{higham08fm}
\begin{equation}\label{eq.iter}
  \begin{aligned}
    X_0 &= A\in\C\nn, \text{nonsingular}\\
    X_{k+1} &= \frac{1}{2}(X_k + X_k^{-*}),\quad k = 0,1,\dots,
  \end{aligned}
\end{equation}
where $(X_k^{-*})$ denotes the conjugate transpose of the inverse of $X$. The following theorem shows the iteration~\eqref{eq.iter} converges to the unitary factor of $A$.
\begin{theorem}
[Newton iteration for the unitary polar factor]
\label{thm.newton-iteration}
Suppose $A\in\C\nn$ is nonsingular, then the iteration~\eqref{eq.iter} produces a sequence $\{X_k\}_k$ such that 
\begin{equation}\notag
  \lim_{k\to\infty} X_k = U,
\end{equation}
where $U$ is the unitary factor of the polar decomposition of $A$.
\end{theorem}

\begin{proof}
See \ycite[1986, Section~3.2]{Hig1986Computingpolardecomposition}.

Let the singular value decomposition of $A$ be $P \Sigma Q\ctp$, where $P$ and $Q$ are unitary. Then by Theorem~\ref{thm.polar-decompo}, we have $A = UH$ where $U = PQ\ctp$ and $H = Q\Sigma Q\ctp$. The trick of the proof is to define a new sequence ${D_k}$ where 
\begin{equation}\label{eq.XtoD}
  D_k = P\ctp X_k Q.
\end{equation}
From iteration~\eqref{eq.iter}, we have 
\begin{equation}\notag
  \begin{aligned}
    D_0 &= P\ctp X_0 Q = P\ctp AQ = \Sigma,\\
    D_{k+1} &= P\ctp X_{k+1} Q = P\ctp \frac 12 (X_k + X_k^{-*})Q \\
        & = \frac 12 (P\ctp X_k Q + P\ctp X_k^{-*}Q).
  \end{aligned}
\end{equation}

Notice that $P\ctp X_k Q = D_k$ and $P\ctp X_k^{-*}Q = (P\ctp X_k Q)^{-*}$, hence the iteration becomes 
\begin{equation}\label{eq.iterD}
  \begin{aligned}
    D_0 = \Sigma, \quad D_{k+1} = \frac 1 2 (D_k + D_k^{-*}).
  \end{aligned}
\end{equation}
Since $D_0$ is diagonal with positive diagonal entries (since $A$ is full rank), hence the iteration~\eqref{eq.iterD} is well-defined and will produce a sequence of diagonal matrices $\{D_k\}$ and we can write $D_k = \diag(d_i^{(k)})$, where $d_i^{(k)} > 0$. Using this notation, we can rewrite the iteration~\eqref{eq.iterD} elementwise:
\begin{equation}\notag
  d_i^{(0)} = \sigma_i,\quad d_i^{(k+1)} = \frac 12 \left(d_i^{(k)} + \frac{1}{d_i^{(k)}}\right)
\end{equation}
where $i = 1,\dots,n$. By an argument discussed in \ycite[1964, p.~84]{Hen1964ElementsNumericalAnalysis}, we can write 
\begin{equation}\label{thm.conv-proof1}
  \begin{aligned}
    d_i^{(k+1)} - 1 = \frac{1}{2d_i^{(k)}}\left({d_i^{(k)}}^2 -2d_i^{(k)} + 1 \right) = \frac{1}{2d_i^{(k)}} \left(d_i^{(k)}-1\right)^2
  \end{aligned}
\end{equation}
This implies that the diagonal entries of $D_k$ converge to $1$ quadratically as $k\to\infty$. Since the argument holds for any $i  = 1,\dots,n$, we conclude that $\lim_{k\to\infty} D_k = I_n$ and use the definition of $D_k$ we have 
\begin{equation}\notag
  \lim_{k\to\infty} X_k = PD_k Q\ctp = PQ\ctp = U,
\end{equation}
as $k\to\infty$, where $U$ is the unitary factor of the polar decomposition of $A$.
\end{proof}

The Newton iteration for the unitary factor of a square matrix $A$ can also be derived by applying the Newton's method to the equation $X\ctp X = I_n$. Consider a function $F(X) = X\ctp X - I_n$, then $X$ is unitary if $X$ is a zero of the function $F(X)$. 

Suppose $Y$ is an approximate solution to the equation $F(\wt X) = 0$, we can write $\wt X = Y + E$ and substitute $\wt X$ into the equation $F(\wt X) = 0$
\begin{equation}\notag
  (Y+E)\ctp (Y+E) - I_n = Y\ctp Y + Y\ctp E + E\ctp Y + E\ctp E - I_n = 0.
\end{equation}
Dropping the second order term, we get Newton iteration $Y\ctp Y + Y\ctp E + E\ctp Y - I_n = 0$ and this is the Sylvester equation for $E$~\ycite[2008, Problem~8.18]{higham08fm}. We aim to solve $E$ and update $Y$ by $Y + E$, which gives the Newton iteration at the $k$th step:
\begin{equation}\label{eq.it1}
  \begin{cases}
    X_k \ctp X_{k} + X_{k}\ctp E_k + E_k\ctp X_k - I_n = 0,\\
    X_{k+1} = X_{k} + E_k.
  \end{cases}
\end{equation}
Assume that $X_k\ctp E_k$ is Hermitian, namely $X_k\ctp E_k = E_k\ctp X_k$, then we can further simplify the iteration by rewritten~\eqref{eq.it1} as $X_k\ctp E_k = (I_n - X_k\ctp X_k)/2$ and the expression we got for $X_k\ctp E_k$ is indeed Hermitian, hence our assumption is valid. Therefore, we have a explicit expression for $E_k$, $E_k = (X_k^{-*} - X_k)/2$, and the iteration~\eqref{eq.it1} becomes 
\begin{equation}\notag
  X_{k+1} = X_{k} + \frac 12 (X_k^{-*}-X_k) = \frac 12(X_k^{-*} + X_k)
\end{equation}
choosing $X_0 = A$, we have our desired Newton iteration~\eqref{eq.iter}.

\subsubsection{Convergence of Newton Iteration}

The convergence of the iteration~\eqref{eq.iter} can be seen from the proof of Theorem~\ref{thm.newton-iteration}.

\begin{theorem}
[Convergence of the iteration~\eqref{eq.iter}, {\ycite[2008, Theorem~8.12]{higham08fm}}] \label{thm.conv.new}
Let $A\in\C\nn$ be non-singular. Then the Newton iterates $X_k$ in~\eqref{eq.iter} converges quadratically to the unitary polar factor $U$ of $A$ with 
\begin{equation}\notag
  \gnorm{X_{k+1} - U} \leq \frac 12 \gnorm{X_k\inv}\gnorm{X_k - U}^2,
\end{equation}
where $\gnorm{\cdot}$ is any unitarily invariant norm.
\end{theorem}

\begin{proof}
Given $A\in\C\nn$, the singular value decomposition of $A$ is $P\Sigma Q\ctp$ and the unitary factor is $U = PQ\ctp$. From~\eqref{thm.conv-proof1}, we have 
\begin{equation}\notag
  d_i^{(k+1)} - 1 = \frac{1}{2d_i^{(k)}}\left(d_i^{(k)}-1\right)^2.
\end{equation}
This is true for all $1\leq i\leq n$, therefore we can rewrite this in matrix form
\begin{equation}\label{eq.3.28}
  \begin{bmatrix}
    d_1\iter{k+1}-1 & & \\
                    & \ddots & \\
                    & & d_n\iter{k+1}-1
  \end{bmatrix} = 
  \frac{1}{2}
  \begin{bmatrix}
    1/d_1\iter{k} & & \\
                  & \ddots & \\
                  & & 1/d_n\iter{k}
  \end{bmatrix}
  \begin{bmatrix}
    d_1\iter{k}-1 & & \\
                  & \ddots & \\
                  & & d_n\iter{k}-1
  \end{bmatrix}^2
\end{equation}

Recall $D_k = \diag(d_i\iter{k})$, then~\eqref{eq.3.28} is same as $D_{k+1}-I = D_k\inv (D_k - I)^2/2$. Therefore we can conclude that $D_{k}$ converges to the identity $I$ quadratically and therefore $X_{k}$ converges to the unitary factor $U$ quadratically.

From $D_{k+1}-I = D_k\inv (D_k - I)^2/2$, by taking the norm on both sides we have 
\begin{equation}\label{eq.Dineq}
  \begin{aligned}
    \gnorm{D_{k+1}-I} &= \frac 12 \gnorm{D_k\inv (D_k - I)^2} \leq \frac 12 \gnorm{D_k\inv} \gnorm{D_k - I}^2
  \end{aligned}
\end{equation}
The final inequality comes from the fact that every unitarily invariant norm is subordinate~\ycite[1997, Proposition~IV~2.4]{Bha1997MatrixAnalysis}. Also, the unitary invariance implies the following three equalities,
\begin{itemize}
\item $\gnorm{D_{k+1}- I} = \gnorm{PD_{k+1} Q\ctp - PQ\ctp} = \gnorm{X_{k+1} - U}$. 
\item $\gnorm{D_k\inv} = \gnorm{QD_k\inv P\ctp} =\gnorm{X_k\inv}$.
\item Using the same argument as 1, we have $\gnorm{D_k - I} = \gnorm{X_k - U}$.
\end{itemize}
As a result, \eqref{eq.Dineq} can be rewritten as 
\begin{equation}\notag
  \gnorm{X_{k+1} - U} \leq \frac 12 \gnorm{X_k\inv}\gnorm{X_k - U}^2.
\end{equation}
\end{proof}

We can test the iteration~\eqref{eq.iter} by the following code 
\begin{lstlisting}
A = my_testmatrices(10,100); X = A;
for i = 1:10
Y = inv(X);
X = (Y' + X)/2;
end
[P,~,V] = svd(A); U = P*V'; 
disp(norm(X-U));
\end{lstlisting}
With only $10$ iterations, the distance between our computed polar factor and the theoretical polar factor (formed by singular value decomposition) is as small as the unit roundoff. However, the disadvantage of this iteration is that it requires the explicit form of the inverse of a matrix at each step, therefore we can introduce another method using one step of the Schulz iteration~\ycite[1965]{Ben1965iterativemethodcomputing} which will remove the requirement of the matrix inverse.

\subsection{Newton--Schulz Iteration Approach}

Given $A\in\C\nn$, an iterative method for computing the inverse $A\inv$, proposed by Schulz \ycite[1933, p.~58]{Sch1933IterativeBerechungder}, is defined by 
\begin{equation}
  \label{eq.schuiter}
  X_{k+1} = X_{k}(2I - AX_k), \quad k = 0,1,2,\dots,
\end{equation}
where $X_0$ is an approximation to $A\inv$. Based on the \eqref{eq.iter}, we can use a one-step Schulz iteration to approximate the inverse of $X_k\ctp$ which can be written as $Y_1 = Y_0(2I - X_k\ctp Y_0)$, where $Y_0$ is an approximation to $X_k^{-*}$. In our situation, the input $X_0$ is an almost unitary matrix ($\norm{X_0\ctp X_0 - I}$ is small) and the output is $U$ which is unitary. Hence the sequence generated by the Newton iteration, $\{X_k\}$, should always be almost unitary, therefore $X_k$ should be a good approximate to $(X_k\ctp)\inv$. The idea of Newton--Schulz iteration can be formulated by the Algorithm~\ref{alg:ns1}.
\begin{algorithm}
\label{alg:ns1}
Given $A\in\C\nn$ such that $\norm{A\ctp A - I} \approx u_\ell$, the algorithm computes the unitary factor $U$ of $A$ defined by Theorem~\ref{thm.polar-decompo}.
\begin{code}
{$X_0 = A$} \\
for \= {$k = 0,1,2,\dots$} \\
\> {$X_k^{-*} = X_k(2I - X_k\ctp X_k)$} \\
\> {$X_{k+1} = \frac{1}{2}\left(X_k^{-*} + X_k\right)$} \\
end
\end{code}
\end{algorithm}

Notice that, the one-step Schulz iteration can be embedded into the Newton iteration by substituting the $X_k^{-*}$ from line $3$ into line $4$ and we have 
\begin{equation}\label{eq.itNS}
  X_{k+1} = \frac{1}{2}\left(X_k(2I - X_k\ctp X_k) + X_k\right) = \frac 12 X_k\left(3I - X_k\ctp X_k\right), \quad X_0 = A,
\end{equation}
and this is known as the Newton--Schulz iteration~\ycite[2008, Section~8.3]{higham08fm}.

\subsubsection{Convergence of Newton--Schulz Iteration}

We would like to know, the iteration~\eqref{eq.itNS} converges to $1$ under which condition. Using the similar approach as proof of Theorem~\ref{thm.newton-iteration}, by writing $D_k = P\ctp X_k Q$ where $P$ and $Q$ are the unitary parts of the singular value decomposition of $A = P\Sigma Q\ctp$, we can consider the iteration~\eqref{eq.itNS} as the scalar form:
\begin{equation}\label{NSscalar}
  d_i\iter{k+1} = \frac 12 d_i\iter{k}\left(3-{d_i\iter{k}}^2\right),\quad d_i\iter{0} = \sigma_i,\quad \text{for $k=0,1,\dots,$}
\end{equation}
where $\sigma_i$ is the $i$th singular value of $A$. To see the convergence of~\eqref{eq.itNS}, we first find the interval of convergence of~\eqref{NSscalar}.

\begin{lemma}\label{lemma.scalarNS}
If $\sigma_i\in (0,\sqrt{3})$, then the iteration~\eqref{NSscalar} converges quadratically to $1$.
\end{lemma}

\begin{proof}
If $\sigma_i = 1$, then the iteration converges immediately. To simplify the notation, we rewrite the iteration~\eqref{NSscalar} as 
\begin{equation}\label{eq.3.34}
  d_{k+1} = f(d_k) = \frac 12 d_k(3-d_k^2),\quad k = 0,1,\dots.
\end{equation}
\quad \textit{Case 1:} If $d_k\in (0,1)$, then $(3-d_k^2)/2 > 1$ which gives $d_k(3-d_k^2)/2 > d_k$. Differentiate $f(d_k)$ once we have $f'(d_k) = -d_k^2 +(3-d_k^2)/2$ and this will always be positive for $d_k\in(0,1)$. Therefore, if we start with $d_k\in (0,1)$, then $d_k < f(d_k) = d_{k+1} < f(1) = 1$, namely, start from $d_0 \in (0,1)$, then the iteration~\eqref{eq.3.34} will converges to $1$.

\quad \textit{Case 2:} If $d_k\in (1,\sqrt{3})$, then by similar approach we have the following
\begin{equation}\notag
  f(1) = 1,\quad f(\sqrt{3}) = 0,\quad f'(d_k) <1 \quad \text{for $d_k\in(1,\sqrt{3})$.}
\end{equation}

Therefore, the function $f$ will monotonically decreasing from $1$ to $0$ within the interval $(1,\sqrt{3})$ and $\sup_{d_k\in(1,\sqrt{3})}f(d_k) = 1$. Hence, if $d_k \in (0,\sqrt{3})$, then $f(d_k)$ will belongs to $(0,1)$, then the iteration~\eqref{eq.3.34} converges to $1$ as discussed in case $1$. The properties of the function $f(d_k)$ in both case $1$ and case $2$ can be seen from the Figure~\ref{fig:11} which shows that for $d_{k}\in (0,1)$, the function $f$ monotonically increasing from $0$ to $1$ and for $d_k\in(1,\sqrt{3})$, the function $f$ monotonically decreasing from $1$ to $0$.

\begin{figure}[!tbhp]
\centering 
\includegraphics[width=0.5\textwidth]{%
  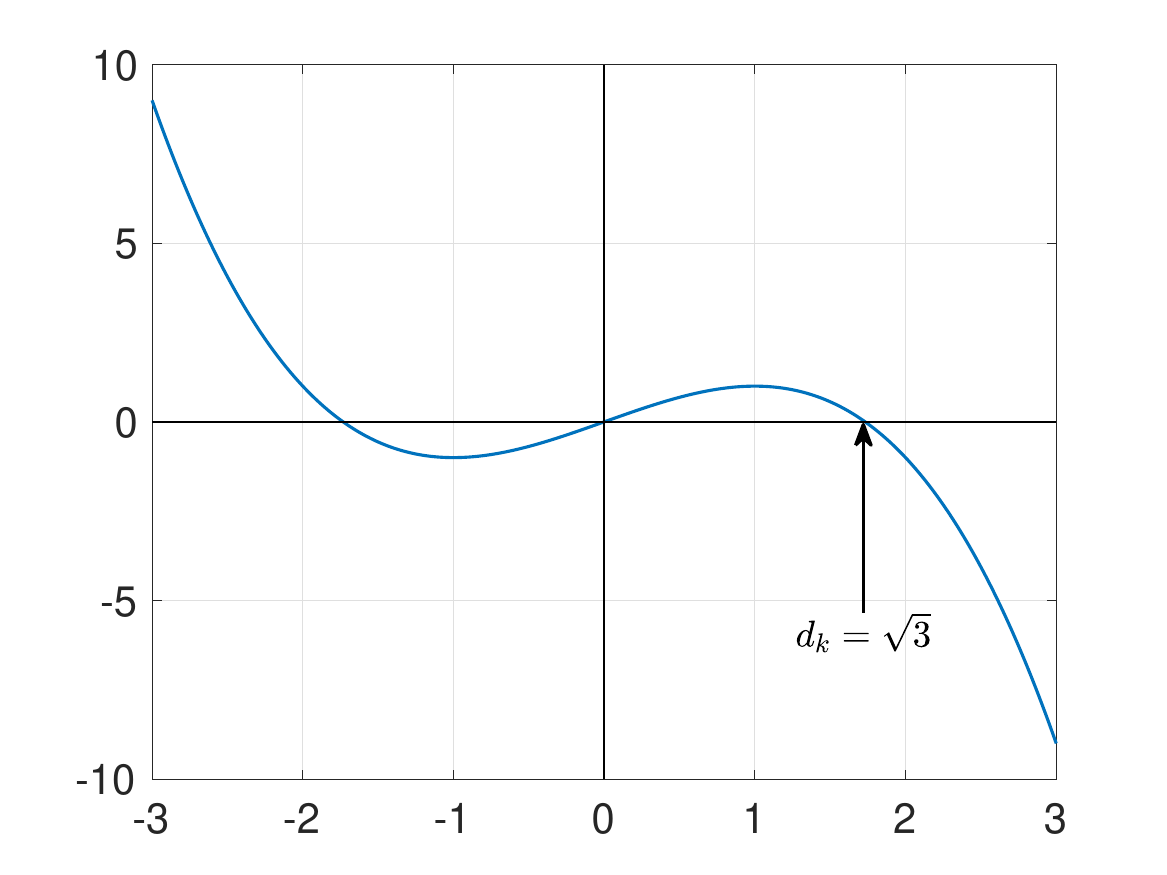}
\caption[The plot of $f(d_k)$ defined in~\eqref{eq.3.34}.]{The plot of $f(d_k)$ defined in~\eqref{eq.3.34} over $[-3,3]$. Three zeros of $f(d_k)$ are $d_k = -\sqrt{3},0$ and $\sqrt{3}$.}
\label{fig:11}
\end{figure}

Other parts of the real line may still converges to $1$ but the interval of converges are too short to consider and control (for example, if we start with $d_0 = 2.2$, then $d_1 = f(d_0) = -2.0240$ and $d_2 = 1.1097 \in (0,\sqrt{3})$ which will converges to $1$).

Remain to show that the convergence is of order $2$. By consider the difference $d_{k+1} - 1$ (similar approach as~\eqref{thm.conv-proof1}), we have 
\begin{equation}\notag
  d_{k+1} - 1 = -\frac 12 (d_k-1)^2 (d_k + 2)
\end{equation}
which implies quadratic convergence.
\end{proof}

From the quadratic convergence of the scalar Newton--Schulz iteration, we have the following theorem.

\begin{theorem}
[Convergence of Newton--Schulz iteration~\eqref{eq.itNS}]
\label{thm.conv.nsiteration}
For $A\in\C\nn$ of rank $n$ which has the polar decomposition $A = UH$, if $\sigma_i(A) \in (0,\sqrt{3})$ or equivalently $\norm{A} < \sqrt{3}$, then the Newton--Schulz iteration 
\begin{equation}\notag
  X_{k+1} = \frac 12 X_k(3I-X_k\ctp X_k),\quad X_0 = A,
\end{equation}
converges to the unitary factor $U$ quadratically as $k\to\infty$ and 
\begin{equation}\label{eq.NSconv}
  \norm{X_{k+1} - U} \leq \frac 12 \norm{X_k + 2U} \norm{X_k - U}^2.  
\end{equation}
\end{theorem}

\begin{proof}
The result is obvious follows the proof of Theorem~\ref{thm.newton-iteration} and~\ref{thm.conv.new} and Lemma~\ref{lemma.scalarNS}.
\end{proof}

\subsubsection{Implementation and Testing}

The Newton--Schulz iteration can be utilized on the almost orthogonal matrix $Q_\ell$ as discussed in Section~\ref{sec.approx-eig} since the eigenvalues of $Q_\ell$ are all near $1$ which definitely satisfies the convergence condition for the Newton--Schulz iteration. For $A\in\C\nn$, iteration~\eqref{eq.itNS} can be formalized by Algorithm~\ref{alg:NSiter} with restrictions on the singular values of $A$ based on Theorem~\ref{thm.conv.nsiteration}.

\begin{algorithm}
\label{alg:NSiter}
Given a matrix $A\in\C\nn$ with $\sigma(A)\subset (0,\sqrt{3})$, this
algorithm computes the unitary factor $U$ of the polar decomposition of $A
= UH$ using the Newton--Schulz iteration~\eqref{eq.itNS}.
\begin{code}
{$X_0 = A$}\\
for \= {$k = 0,1,\dots$}\\
\> {Compute $X_k\ctp X_k$}\\
\> if \= {$\gnorm{X_k\ctp X_k - I}_F \lesssim nu_d$}\\
\> \> {\bfseries Break}, output $X_k$\\
\> end\\
\> {$X_{k+1} = \frac{1}{2}X_k(3I-X_k\ctp X_k)$}\\
end
\end{code}
\end{algorithm}

Notice that, in line 3, we use the Frobenius norm to examine the orthogonality instead of the 2-norm, this is due to the Frobenius norm being easier and cheaper to compute. We can further simplify the algorithm based on the following analysis.

From now on, the explanation will focus on the 2-norm, but it can be generalized to all unitarily invariant norms. Suppose $Q_\ell$ has a singular value decomposition $Q_\ell = P\Sigma Q\ctp$, then we can rewrite $\norm{Q_\ell\ctp Q_\ell - I}$ as 
\begin{equation}\notag
  \begin{aligned}
    \norm{Q_\ell\ctp Q_\ell - I} & = \norm{Q\Sigma\ctp P\ctp P\Sigma Q\ctp -I} = \norm{\Sigma^2 - I} \approx nu_\ell.
  \end{aligned}
\end{equation}

Since $\Sigma = \diag(\sigma_1,\dots,\sigma_n)$ where $\sigma_i$ are sorted in non-ascending order. Using Corollary~\ref{col.2norm} and the fact that the singular values of the almost unitary matrix are around one, we have $\norm{\Sigma^2 - I} \approx {\sigma_1^2 - 1} \approx n \cdot 10^{-8}$, therefore 
\begin{equation}\notag
  \norm{\Sigma - I} \approx \sigma_1 - 1 \approx \sqrt{1 + n\cdot 10^{-8}} - 1.
\end{equation} 

From relation~\eqref{eq.NSconv} and the fact that $\norm{Q_\ell - U} = \norm{\Sigma - I}$, using the notations in~\eqref{eq.itNS}, we have 
\begin{equation}\notag
  \norm{X_1 - U} \leq \norm{X_0 - U}^2 = \norm{Q_\ell -U}^2 = \norm{\Sigma - I}^2 \approx (\sqrt{1+n\cdot 10^{-8}} - 1)^2.
\end{equation}
Similarly, after the second iteration, we have 
\begin{equation}\notag
  \norm{X_2 - U} \leq \norm{X_1 - U}^2 \approx (\sqrt{1+n\cdot 10^{-8}} - 1)^4.
\end{equation}

\begin{figure}[!tbhp]
\centering
\includegraphics[width=0.55\textwidth]{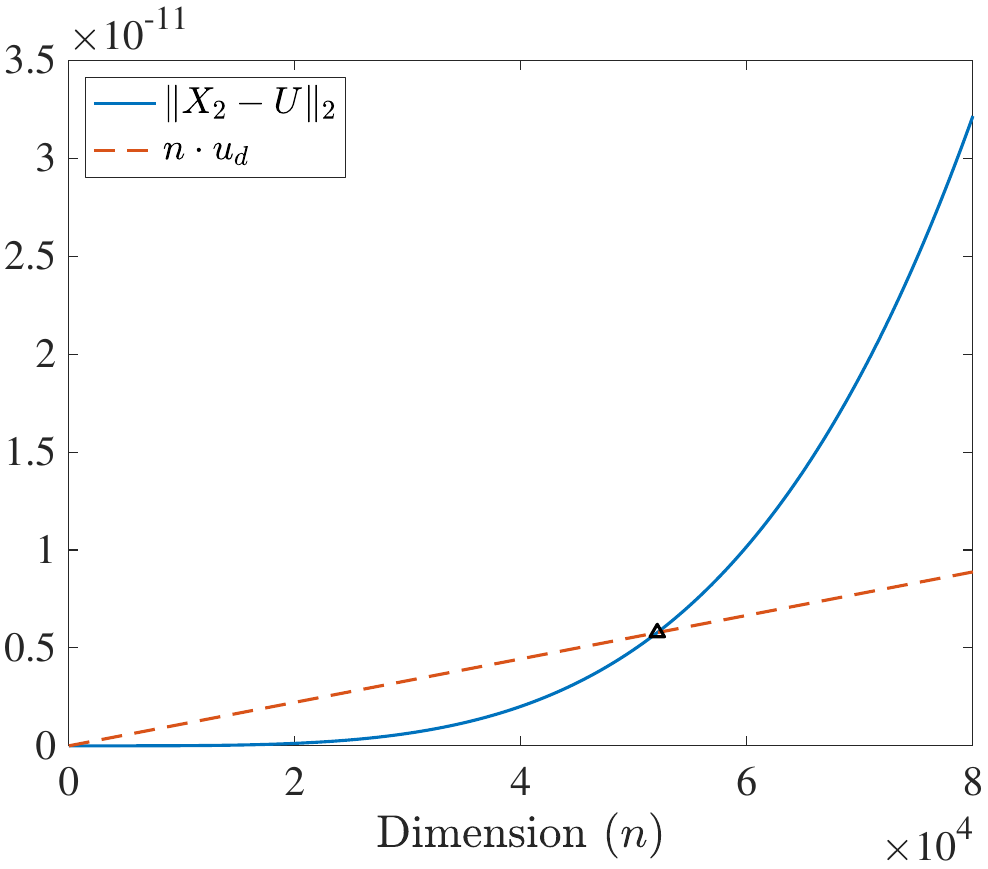}
\caption[The difference between the iterator $X_2$ and the unitary factor $U$ against the dimension $n$.]{The figure shows, after two iterations, the difference between the matrix $X_2$ (from Algorithm~\ref{alg:NSiter}) and the unitary factor $U$ with respect to the dimension. The black marked point is the intersection at $n = 52071$.}
\label{fig:seconditer}
\end{figure}

We can compare this quantity with the desired tolerance $nu_d$ and we can produce Figure~\ref{fig:seconditer} using Appendix~\ref{app:seconditer}.

From Figure~\ref{fig:seconditer}, we observe that for $n$ less than $52000$, Algorithm~\ref{alg:NSiter} will converge in $2$ steps. If the dimension exceeds $52400$, by the same analysis, for any $Q_\ell\in\R\nn$, if $n < 2172548$, then the Newton--Schulz iteration will converge in three steps. In practice, MATLAB is not able to generate such a big matrix. Therefore, we can remove the terminate condition in Algorithm~\ref{alg:NSiter} as shown in Algorithm~\ref{alg:NSiter2},

\begin{algorithm}
\label{alg:NSiter2}
Given a matrix $Q_\ell\in\C\nn$ such that $\norm{Q_\ell\ctp Q_\ell - I} \approx n\cdot 10^{-8}$ and $n < 2\times 10^6$, this algorithm computes the unitary factor $U$ of the polar decomposition of $Q_\ell = UH$ using the Newton--Schulz iteration~\eqref{eq.itNS}.
\begin{code}
{$X_0 = Q_\ell$}\\
if \= {$n \geq 52400$} \\
\> {$\alpha = 2$}\\
else \\
\> {$\alpha = 1$}\\
end\\
for \= {$k = 0:\alpha$}\\
\> {$X_{k+1} = \frac{1}{2}X_k(3I-X_k\ctp X_k)$}\\
end\\
{Output $X_{\alpha + 1}$} 
\end{code}
\end{algorithm}

Implementation of Algorithm~\ref{alg:NSiter2} in \mat~gives 
\begin{lstlisting}
function A = newton_schulz(A)
[~,n] = size(A); 
if n >= 52000, iter = 3;
else, iter = 2; end
for k = 1:iter
A = 0.5*A*(3*eye(n) - A'*A); 
end
\end{lstlisting}
For each iteration, there are only 2 matrix-matrix multiplications which cost $4n^3$ flops and hence the total cost when applying \inline{newton_schulz} to $A\in\C\nn$ will be about $8n^3 + O(n^2)$ flops. Compare to the singular value decomposition approach discussed in Section~\ref{sec:svdapproach} which costs about $21n^3$, the Newton--Schulz iteration is usually the preferred method when the input is almost unitary.

\begin{figure}[!tbhp]
\centering
\includegraphics[width=0.5\textwidth]{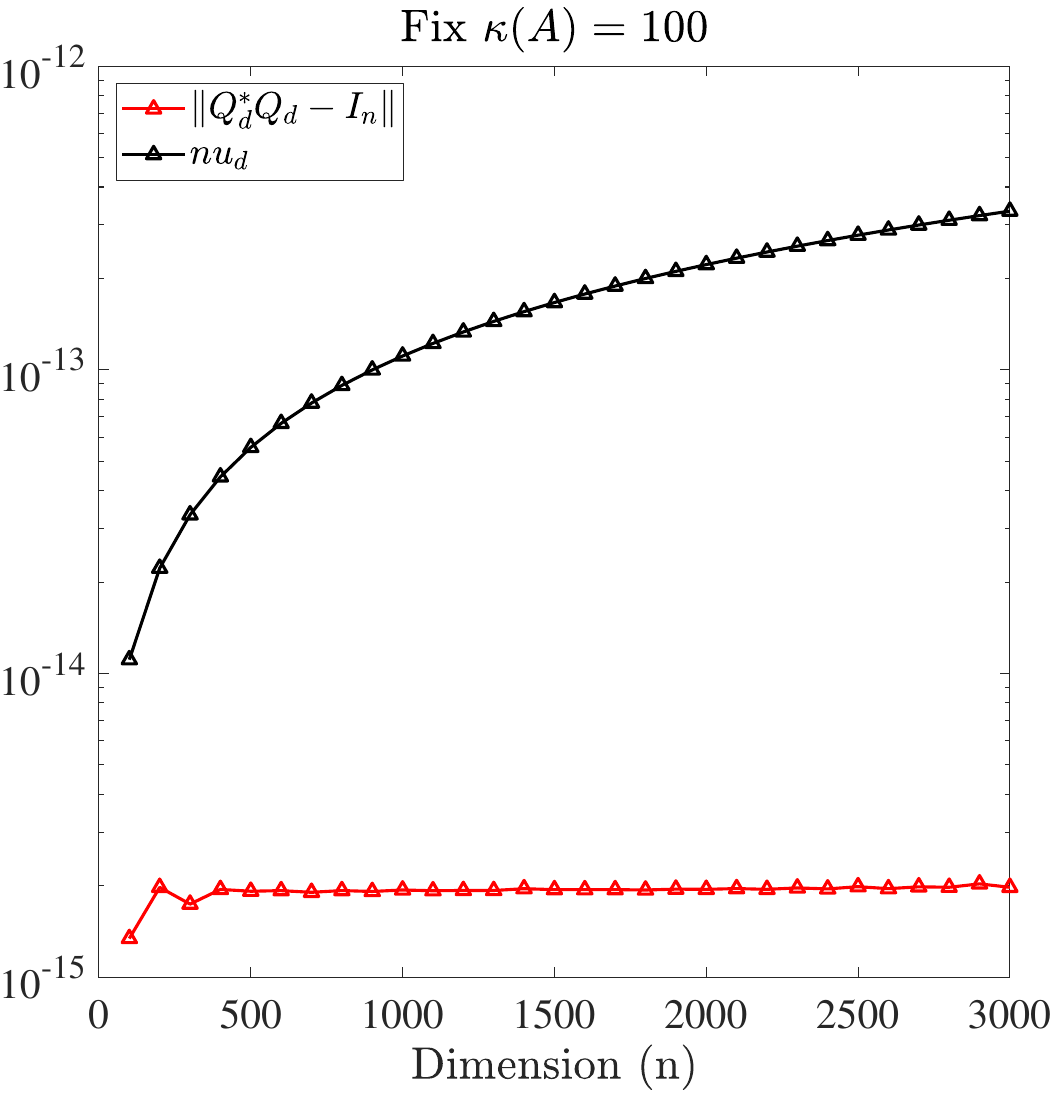}
\caption[Behavior of $\norm{Q_d\ctp Q_d - I}$ with the dimension of the input $Q$ increases and $\kappa(A)$ fixes using the Newton--Schulz iteration.]{This figure shows how the quantity $\norm{Q_d\ctp Q_d - I}$ increases as the dimension of the input matrix $Q$ increases while we fix $\kappa(A)$. The red line shows the quantity $\norm{Q_d\ctp Q_d - I}$ after orthogonalization using Newton--iteration and the black line is the reference line $nu_d$.}
\label{fig:nsorthog}
\end{figure}

We will perform the following test: Given a matrix $A\in\R\nn$ with fixed
condition number, here $\kappa(A) = 100$, we compute its spectral
decomposition at single precision $A = QDQ\ctp$. Then we use this $Q$ as
our input of \inline{newton_schulz} and see how the quantity
$\norm{Q_d\ctp  Q_d - I}$,
where $Q_d$ is the output of the function,
changes as the dimension $n$ changes.
Using the code from Appendix~\ref{app:fig:nsorthog},
we produce Figure~\ref{fig:nsorthog}. 

From Figure~\ref{fig:nsorthog}, the quantity $\norm{Q_d\ctp Q_d - I}$ is well bounded by the tolerance $nu_d$ and the result is consistent with our expectations and analyses.

\section{Numerical Comparison}

In this section, we will examine the performance and accuracy of both the QR factorization and the polar decomposition using the Newton--Schulz iteration. Theoretically, for our almost orthogonal matrix $Q_\ell$, the QR factorization requires about $10/3n^3$ flops and the Newton--Schulz iteration method requires about $8n^3$ flops. However, MATLAB computes matrix-matrix products faster, hence Newton--Schulz iteration has the potential to be faster than the QR factorization. We set up our numerical comparison via the following 
\begin{enumerate}
\item Generate the symmetric matrix $A\in\R\nn$ and its approximate
spectral decomposition $A = Q_\ell D_\ell Q_\ell\ctp$. Here, we fixed $\kappa(A) = 100$. Even though no matter how we change the condition number of $A$, the eigenvector matrix at single precision, $Q_\ell$, will always have condition number near $1$, recall this is also the reason why we can always apply the Newton--Schulz iteration to the matrix $Q_\ell$.
\item For different $n$, applying \inline{myqr} and \inline{newton_schulz} to $Q_\ell$ and output $Q_d$ which has been numerically proved that it is orthogonal at double precision, then we examine these functions via
\begin{enumerate}
\item Accuracy: $\norm{Q\ctp Q - I}$.
\item Performance : time used measured by \inline{tic,toc}.
\end{enumerate}
\end{enumerate}

\subsection{Accuracy}\label{sec:Accuracy}

\begin{figure}[!tbhp]
\centering
\includegraphics[width=0.5\textwidth]{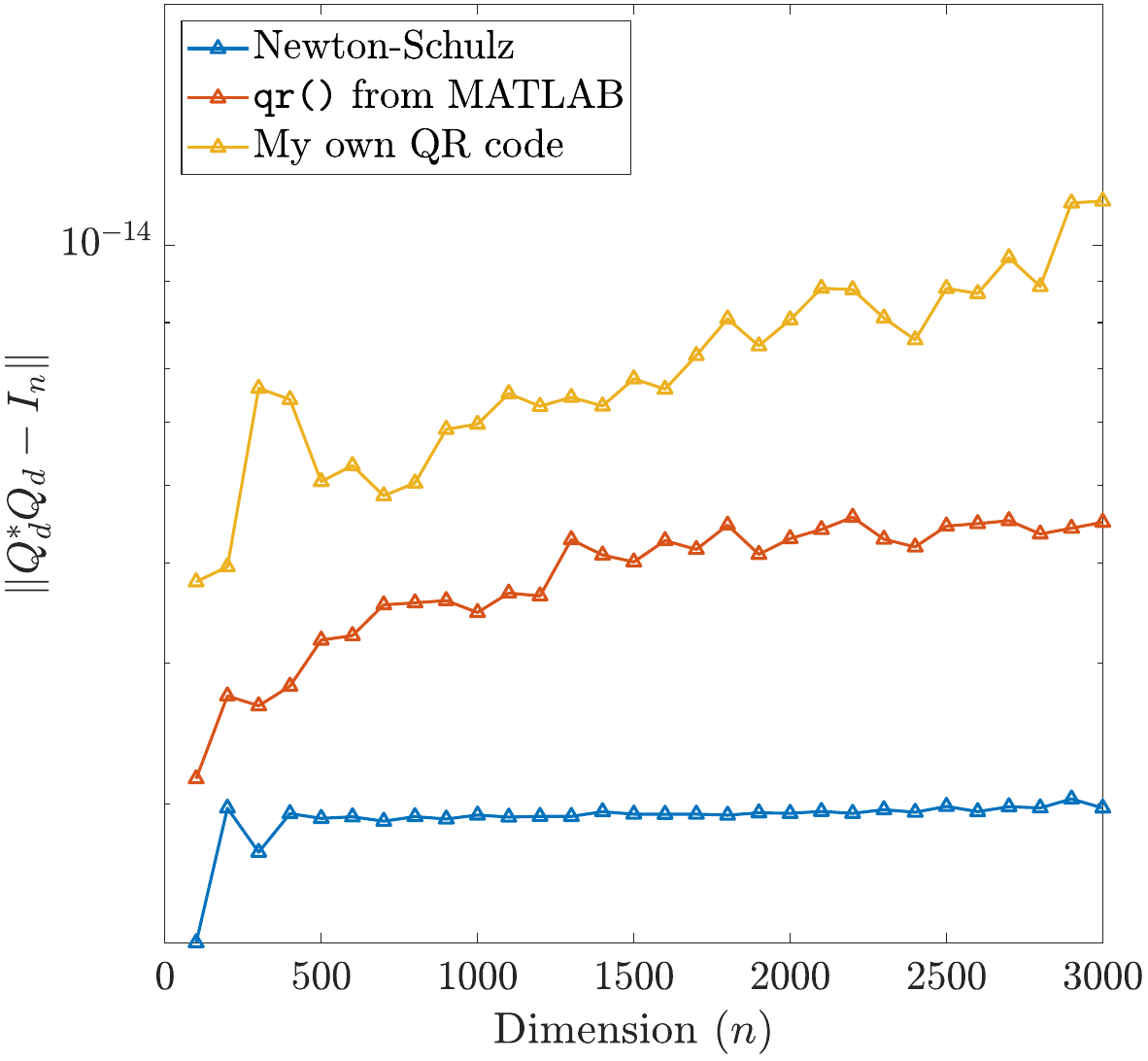}
\caption[Behavior of the quantity $\norm{Q_d\ctp Q_d - I_n}$ with the dimension $n$ increases using both the QR factorization and the Newton--Schulz iteration.]{Behavior of $\norm{Q_d\ctp Q_d - I_n}$ changes as the dimension $n$ increases using both the QR factorization and the Newton Schulz iteration. The blue line shows the results using the Newton--Schulz iteration, the red line shows the results using the MATLAB built-in \inline{qr()} function and the black line shows the results using the my own QR factorization code \inline{myqr()} from Section~\ref{sec.qrtesting}. The code to regenerate this graph can be found in Appendix~\ref{app:fig:acccompare}.}
\label{fig:acccompare}
\end{figure}

The accuracy is measured by $\norm{Q_d\ctp Q_d - I}$ where $Q_d$ is the output of \inline{myqr} and \inline{newton_schulz}. From previous testings, both QR factorization and Newton--Schulz iteration are accurate enough to generate an orthogonal matrix at double precision. Therefore, here we are comparing the accuracy between these two methods shown in Figure~\ref{fig:acccompare}. In addition, we add the accuracy of the MATLAB built-in function \inline{qr()} served as a reference line.

From the figure, we can see that the Newton--Schulz iteration is generally more accurate than the QR factorization method. 

\subsection{Performance}\label{sec:qr-ns-performance}

\begin{figure}[!tbhp]
\centering
\includegraphics[width=0.5\textwidth]{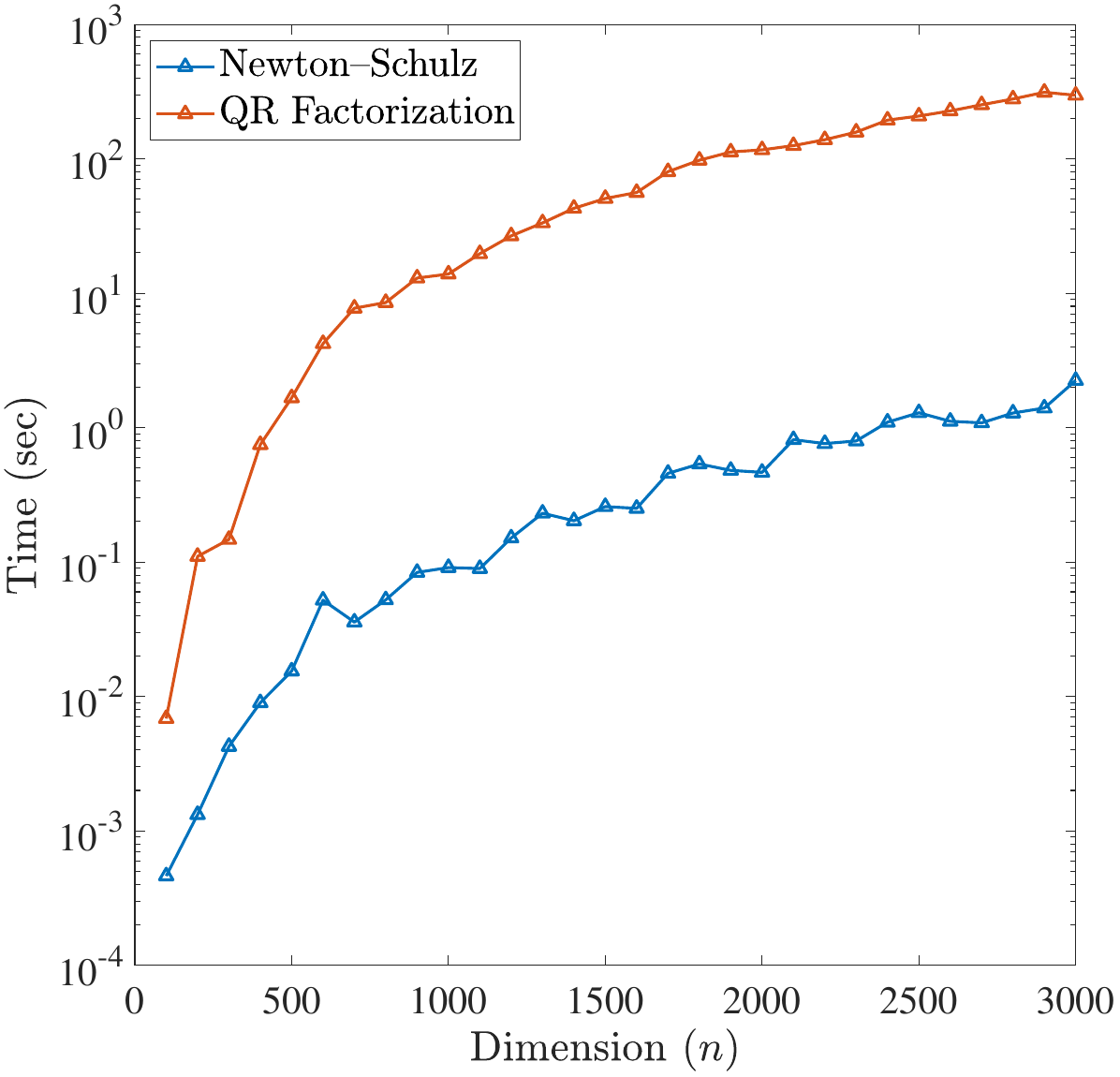}
\caption[The time of applying the QR factorization and the Newton--Schulz iteration to the same almost orthogonal matrix $Q$ with respect to the dimension $(n)$.]{The time of applying the QR factorization and the Newton--Schulz iteration to the same almost orthogonal matrix $Q$. The red line shows the results using the Newton--Schulz iteration and the black line shows the results using the QR factorization. The code to regenerate this graph can be found in Appendix~\ref{app:fig:timecompare}.}
\label{fig:timecompare}
\end{figure}

The performance can be measured by the MATLAB built-in function \inline{tic} and \inline{toc}. Using similar code from Section~\ref{sec:Accuracy}, we are able to compare the time used by \inline{myqr} and \inline{newton_schulz} and produce Figure~\ref{fig:timecompare}.

From the figure, the Newton--Schulz iteration is generally faster than the QR factorization approach. 

Consequently, based on the accuracy and performance, we shall stick with the Newton--Schulz iteration approach as the method to orthogonalize the matrix $Q_\ell$.


\chapter{Mixed Precision Jacobi Algorithm}\label{chap:mixed-precision}

In this chapter, a mixed precision algorithm will be presented with intensive testings on its convergence, speed and accuracy. As mentioned in the Section~\ref{sec:test-matrices}, we will perform the testings for both mode \inline{'geo'} and mode \inline{'ari'}. 

\section{Algorithm}\label{sec:precondition-algorithm}

As discussed from~\ycite[2000]{DV2000Approximateeigenvectorsas}, we can improve the Jacobi algorithm by using a preconditioner. From the discussions in Chapter~\ref{chap:jacobi_algorithm} and~\ref{chap:orthogonalisation}, we assemble them all and collapse into the following algorithm,
\begin{algorithm}
\label{alg:jacobi-preconditioned}
Given a symmetric matrix $A\in\R\nn$ and a positive tolerance $tol$, this
algorithm computes the spectral decomposition of $A$ using the
cyclic-by-row Jacobi algorithm with preconditioning. 
\begin{code}
{Compute the spectral decomposition in single precision, $A = Q_\ell D_\ell
  Q_\ell\tp$.}\\
Orthogonalize $Q_\ell$ by using the Newton--Schulz iteration as shown in
  the Algorithm~\ref{alg:NSiter2}, \nonumberbreak and we have the output
  $Q_d$.\\
Precondition $A$ via $A_{\text{cond}} = Q_d\tp A Q_d$ and apply the
  cyclic-by-row Jacobi algorithm
  \nonumberbreak on $A_{\text{cond}}$ to get the orthogonal
  matrix $V$ and the diagonal matrix $\varLambda$ such that $V\varLambda
  V\tp = A_{\text{cond}}$.\\
{Output the diagonal matrix $\varLambda$ and the orthogonal matrix $Q = Q_d
  V$ such that $A = Q\varLambda Q\tp$.}
\end{code}
\end{algorithm}

\subsection{Implementation and Testing}\label{sec:ref111}

It is straightforward to implement Algorithm~\ref{alg:jacobi-preconditioned} by utilizing the existing codes. MATLAB code can be found in Appendix~\ref{app:mixed-precision-jacobi}. 

\begin{figure}[!tbhp]
  \centering
  \includegraphics[width=0.85\textwidth]{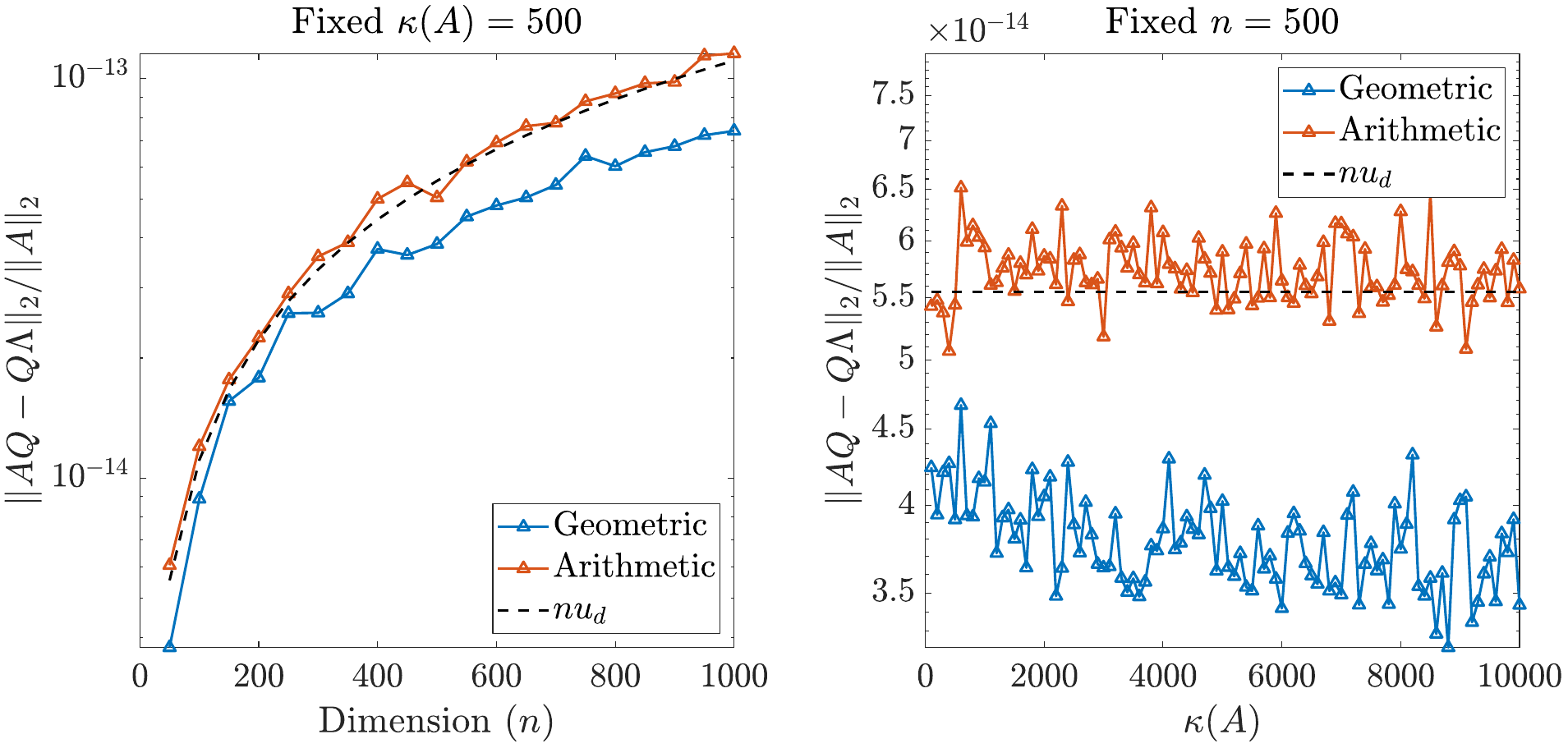}
  \caption[Behavior of the relative error $\norm{AQ-Q\varLambda}/\norm{A}$ where $A\in\R\nn$ is symmetric and $Q$ and $\varLambda$ are the results computed using Algorithm~\ref{alg:jacobi-preconditioned}.]{Behavior of the relative error $\norm{AQ-Q\varLambda}/\norm{A}$ with the dimension $n$ and the condition number $\kappa(A)$ increases. The matrices for the testings on the blue line have geometrically distributed singular values and those on the red line have arithmetically distributed singular values. The black line is the reference tolerance $nu_d$.}
  \label{fig:precondition-accuracy}
  \end{figure}

We would like to first access its accuracy concerning both the dimension and the condition number. Using code in Appendix~\ref{app:precondition-accuracy} we can produce Figure~\ref{fig:precondition-accuracy}.

As shown from the figure, the Algorithm~\ref{alg:jacobi-preconditioned} gives more accurate results for the matrices with geometrically distributed singular values. But regardless of the distributions, the computed $Q$ and $\varLambda$ satisfies the condition 
\begin{equation}\notag
  \norm{AQ - Q\varLambda} \lesssim nu_d\norm{A},
\end{equation}
and therefore this algorithm can be used to
find the spectral decomposition of a real symmetric matrix in practice.
The remaining of this chapter is to find out two things  
\begin{enumerate}
  \item How many sweeps of Jacobi will be necessary for convergence for a preconditioned real symmetric matrix $A_\text{cond}$ as described in Algorithm~\ref{alg:jacobi-preconditioned}?
  \item How much time can be saved by using preconditioning?
\end{enumerate}

\section{Convergence and Speed}

In this section, we will analyze how many sweeps will be sufficient for the cyclic-by-row Jacobi algorithm converges for a preconditioned symmetric matrix $A_{\text{cond}}$.

\subsection{Preconditioning}\label{p3:sec:precondition}

Using the notations from Algorithm~\ref{alg:jacobi-preconditioned}, $Q_d$ is the nearest orthogonal matrix compare to the approximate eigenvector matrix $Q_\ell$, hence $A_{\text{cond}}$ is expected to be almost diagonal. Consequently, we expect this eigenvalue problem should be able to be solved more quickly by the Jacobi algorithm \ycite[2000, Sec.~1, Corollary~3.2]{DV2000Approximateeigenvectorsas}. Quantitatively, after preconditioning, the magnitude of the off-diagonal entries of $A_\text{cond}$ should be no more than $nu_\ell \norm{A}$, which is equivalently saying $\off(A_\text{cond}) \lesssim n^2(n-1)u_\ell \norm{A}$. In practice, this bound is sharper than the theoretical bound and we can have 
\begin{equation}\label{p3:eq:offAbound}
  \off(A_\text{cond}) \lesssim nu_\ell \norm{A}.
\end{equation}

This upper bound can be numerically tested by the following procedure: First, generate a symmetric matrix $A\in\R\nn$ with controlled dimension and condition number, then we perform steps 1,2 and 3 in the Algorithm~\ref{alg:jacobi-preconditioned} to see how much $\off(A)$ is reduced based on this precondition technique. Using the code in Appendix~\ref{p3:sec:figure1}, we generate Figure~\ref{p3:fig:approx-eig-test}.

\begin{figure}[!tbhp]
\centering
\includegraphics[width=0.8\textwidth]{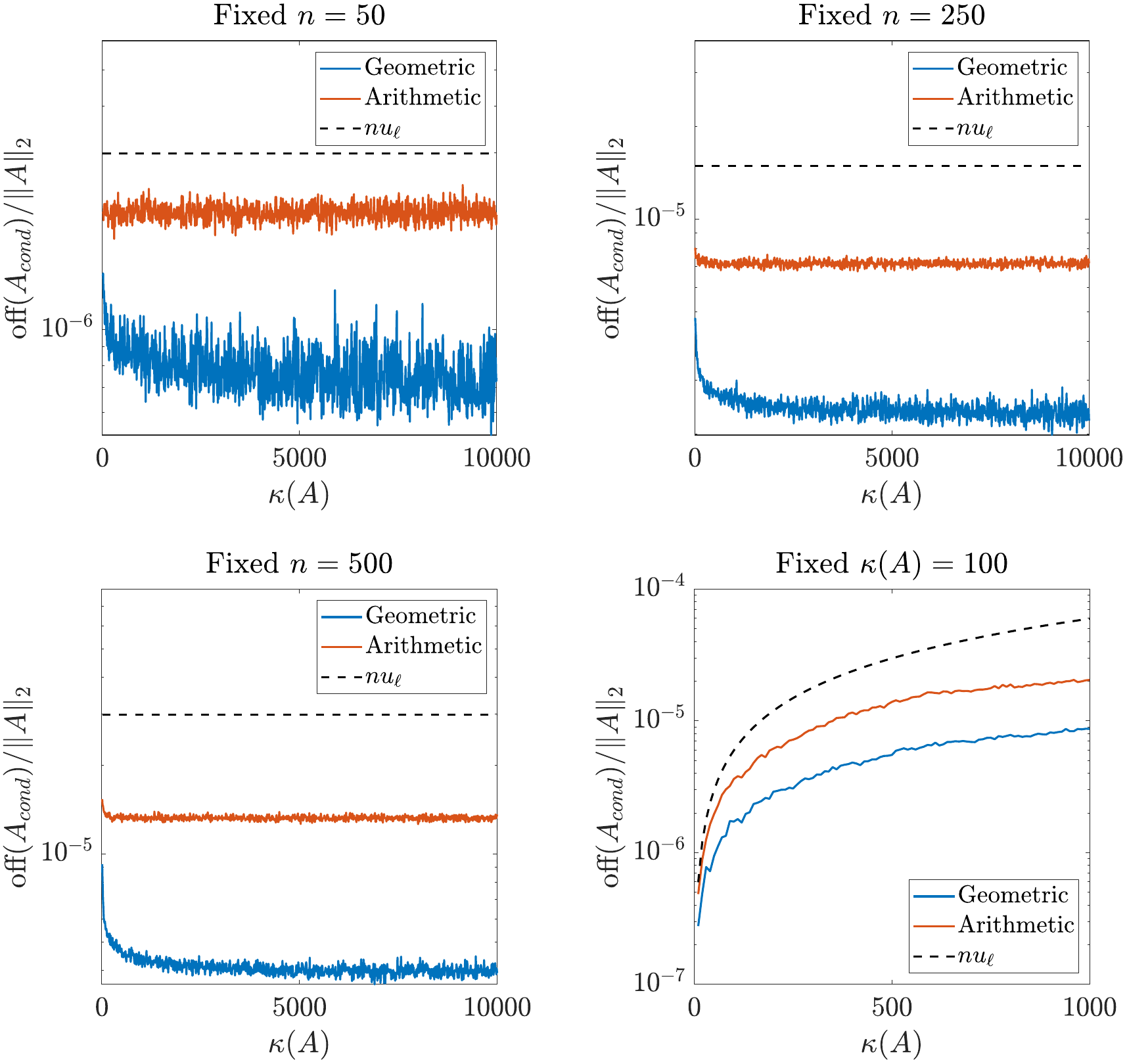}
\caption[Four comparisons of $\off(A_\text{cond})/\norm{A}$ with the reference line $nu_\ell$ where $A_\text{cond}$ is the preconditioned matrix described in the Algorithm~\ref{alg:jacobi-preconditioned}.]{Top-left, Top-right, Bottom-left: These figures show the $\off(A_\textup{cond})/\norm{A}$ for $A\in\R\nn$ with different condition number while fixing $n = 50, 250$ and $500$ respectively and varying $\kappa(A)$ from $10$ to $10^3$. Bottom-right: This figure shows the $\off(A_\textup{cond})/\norm{A}$ for $A\in\R\nn$ with different $n$ while fixing $\kappa(A) = 100$ and varying $n$ from $10$ to $1000$. The blue and red line represent the testings on the real symmetric matrices with geometrically and arithmetically distributed singular values respectively.}
\label{p3:fig:approx-eig-test}
\end{figure}

From the top-left, top-right and bottom-left parts of the Figure~\ref{p3:fig:approx-eig-test}, we see that for $n=50,250$ and $500$, the quantity $\off(A_\text{cond})/\norm{A}$ is always smaller or similar to $nu_\ell$ regardless of its singular values distribution. Notice that, as $n$ changes from $50$ to $500$, the perturbation of the quantity $\off(A_\text{cond})/\norm{A}$ becomes more and more negligible compare to $nu_\ell$. Together with the bottom-right figure, as the size of $A$ passes $100$, the quantity $\off(A_\text{cond})/\norm{A}$ is always well bounded by $nu_\ell$. Based on these observations, we verify the idea of preconditioning and prove numerically that after preconditioning $A$ by $Q_d\tp A Q_d$, the off-diagonal entries are reduced significantly and this property can be captured by the Jacobi algorithm.

\subsection{Quadratic Convergence}\label{sec:quadratic-conv}

Based on the preconditioning process, in this section, we will discuss the convergence of the cyclic Jacobi algorithm applied to a preconditioned real symmetric matrix. Then we will deliver some numerical testings concerning different dimensions and condition numbers. Finally, we will test the improvements in speed by using the preconditioners.

Recall that the number of iterations in one sweep is equal to half of the number of off-diagonal entries,
\begin{equation}\notag
  \text{one sweep of $A$} = \frac{n(n-1)}{2},\quad A \in\R\nn.
\end{equation}

From the previous derivation, the Jacobi algorithm is proved to converge linearly. However, an improvement by Sch{\"o}nhage in 1961~\ycite[1961]{1961-first-quadratic-convergence-Schu} stated that the Jacobi algorithm for symmetric matrices with distinct eigenvalues is converging quadratically. Let us denote $A\iter{k}\in\R\nn$ be the matrix produced by the $k$th iteration of the Jacobi algorithm and therefore $A\iter{0}\in\R\nn$ is the original input matrix. Denote $\lambda_1,\dots,\lambda_n$ be the eigenvalues of $A\iter{0}$ and suppose we have 
\begin{equation}
  \label{p3:eq:eig-distance}
  \abs{\lambda_i - \lambda_j} < 2\delta,\quad i\neq j.
\end{equation}

Suppose we reach the stage that the $\off(A\iter{k}) < {\delta}/{8}$ and the angles of rotation generate by each iteration are smaller than $\pi/4$, as controlled in Section~\ref{sec:2.3}, Sch{\"o}nhage's result said 
\begin{equation}\notag
  \off(A\iter{N+k}) \leq \frac{n\sqrt{(n-2)/2}}{\delta}\off(A\iter{k})^2,\quad N = \frac{n(n-1)}{2},
\end{equation}
which implies quadratic convergence. Later in 1962, Wilkinson provided a sharper bound~\ycite[1962, Section~3]{Wil1962Notequadraticconvergence}. Under the same condition as described by Sch{\"o}nhage, Wilkinson proposed
\begin{equation}\notag
  \off(A\iter{N+k}) \leq \off(A\iter{k})^2 /\delta.
\end{equation}

In 1966, van Kempen proved the quadratic convergence without the assumption of distinct eigenvalues. In his paper~\ycite[1966, Section~2]{vKem1966quadraticconvergencespecial}, instead of define $\delta$ as in \eqref{p3:eq:eig-distance}, he denoted $\delta$ as 
\begin{equation}
  \label{p3:eq:eig-value-dist-legit}
  \min_{\l_i \neq \l_j} \abs{\l_i - \l_j} \geq 2\delta.
\end{equation}
Suppose after $k$ iterations, $\off(A\iter{k}) < \delta/8$, then 
\begin{equation}\notag
  \off(A\iter{k+N}) \leq \frac{\sqrt{\frac{17}{9}}\off(A\iter{k})^2}{\delta}
\end{equation}
Note that although the above bound is correct, the proof in~\ycite[1966]{vKem1966quadraticconvergencespecial} is not correct. This mistake was unveiled by Hari who proposed a sharper bound~\ycite[1991, Section~2]{Har1991sharpquadraticconvergence}.

Back to our situation, as described in Section~\ref{p3:sec:precondition}, after preconditioning, $\off(A) \lesssim nu_\ell \norm{A}$. Assuming our eigenproblem is well conditioned and the minimum distance between the distinct eigenvalues is large enough such that 
$$\off(A) \lesssim 0.1 < \delta / 8$$
where $\delta$ is defined by~\eqref{p3:eq:eig-value-dist-legit}. Then we can expect quadratic convergence and the number of iterations should not exceed five since it will only take $4$ sweeps for the quantity $\off(A)$ reduces from $0.1$ to $10^{-16}$ ($0.1^4 = 10^{-16}$).

\subsection{Numerical Testing}
\begin{figure}[!tbhp]
\centering
\includegraphics[width=0.8\textwidth]{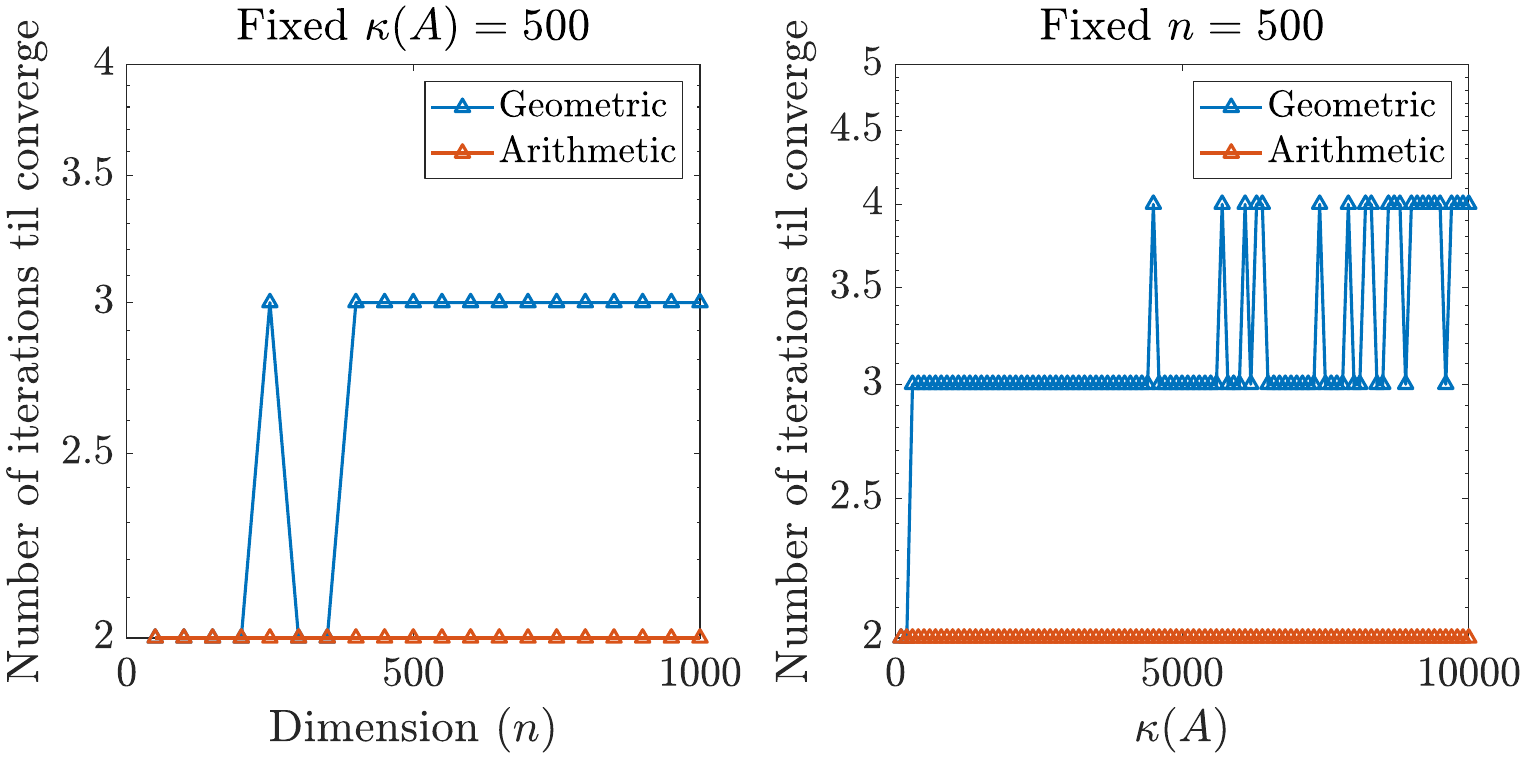}
\caption[Behavior of the number of iterations with the dimension and the condition number of the input matrix vary.]{Behavior of the number of iterations with the dimension and the condition number of the input matrix vary. Matrices in the left and right figures fix its condition number and dimension respectively. The blue and red triangular dots represent the number of iterations for each generated matrices with geometrically and arithmetically distributed singular values respectively.}
\label{fig:typical-sweep}
\end{figure}

By referring to the function \inline{jacobi_precondi} in Appendix~\ref{app:mixed-precision-jacobi}, we can output the number of iterations required. Therefore, we can investigate the number of sweeps required using the routine from Appendix~\ref{app:typical-sweep}. The procedure is exactly the same as the testing in Section~\ref{sec:ref111} and we can produce Figure~\ref{fig:typical-sweep}.

Suppose our testing matrices are not ill-conditioned, then for the matrices with arithmetically distributed singular values, it usually requires two iterations for the cyclic-by-row Jacobi algorithm to converge since the eigenvalue are relatively far apart and this can result in quadratic convergence as discussed in Section~\ref{sec:quadratic-conv}. However, for the matrices with geometrically distributed singular values, the algorithm can take four iterations to converge. By construction, the geometrically distributed singular values has smaller $\delta$ in \eqref{p3:eq:eig-value-dist-legit} and this can result in slower convergence since it will require more sweep for $\off(A)$ to satisfies $\off(A) < \delta/8$. We can exaggerate this observation by generating a matrix with very close eigenvalues. 
\begin{lstlisting}
close all; clear; clc; rng(1,'twister');
A_geo = my_randsvd(1e3,1e16,'geo');
A_ari = my_randsvd(1e3,1e16,'ari');
\end{lstlisting}
The matrix \inline{A_geo} has $\min_{\l_i \neq \l_j}\abs{\l_i - \l_j} \approx 10^{-18}$, whereas the matrix \inline{A_ari} has $\min_{\l_i \neq \l_j}\abs{\l_i - \l_j} \approx 10^{-3}$. Clearly, the Algorithm~\ref{alg:jacobi-preconditioned} can hardly attained the quadratic rate of convergence on \inline{A_geo}.

After applying Algorithm~\ref{alg:jacobi-preconditioned} on both \inline{A_geo} and \inline{A_ari}, the former testing requires 25 sweeps to converge and the latter one only need three sweeps. Therefore, the Jacobi algorithm performs better when the distance between distinct eigenvalues are large since the algorithm can take advantage of quadratic convergence.

\begin{figure}[!tbhp]
\centering
\includegraphics[width=1\textwidth]{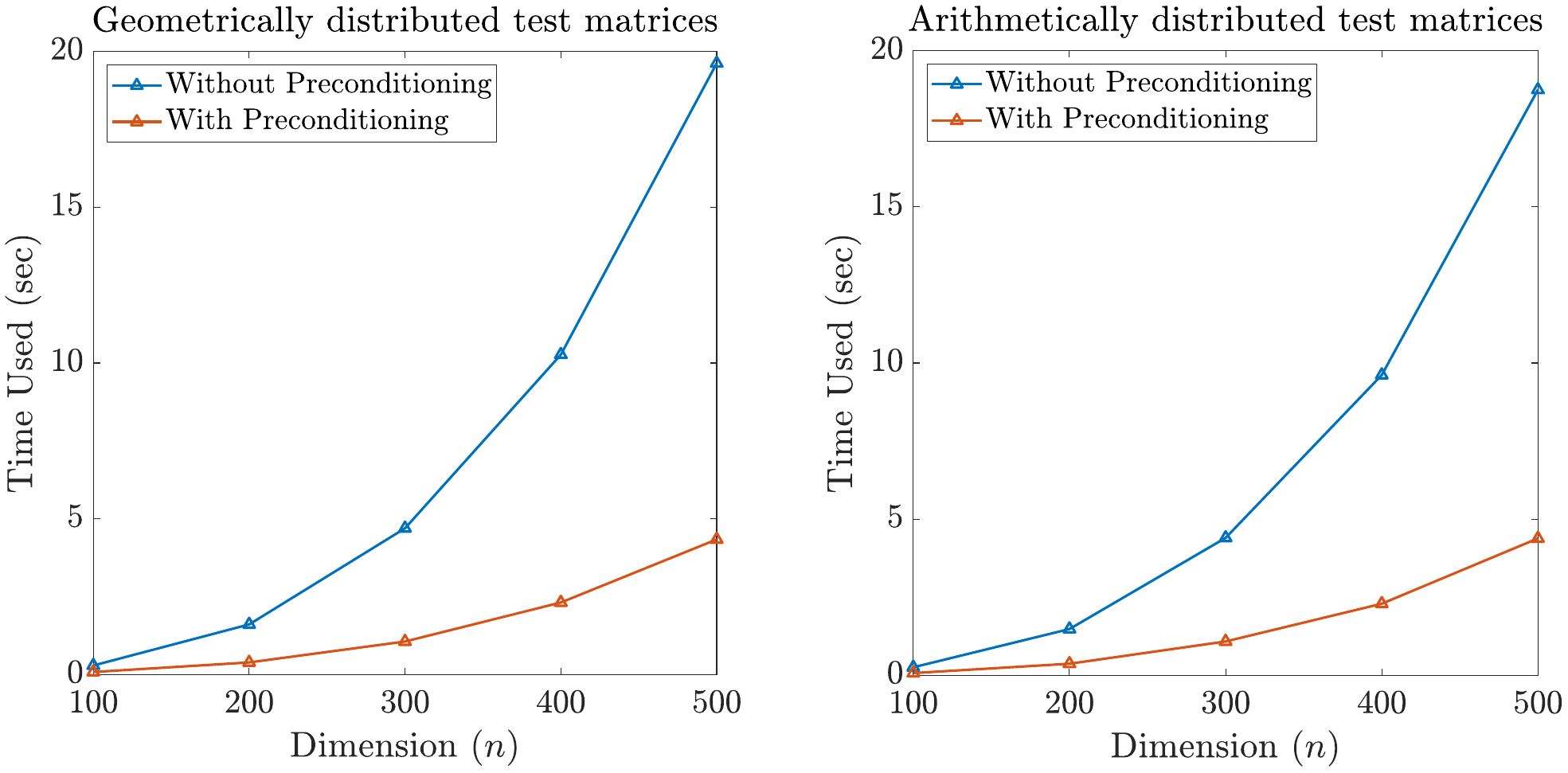}
\caption[Time used by the cyclic-by-row Jacobi algorithm with and without preconditioning.]{Time used by the cyclic-by-row Jacobi algorithm with and without preconditioning on a real symmetric matrix $A$ with respect to the dimension $n$. Here we fixed the condition number $\kappa(A) = 100$. The left and right figures show the testings on the real symmetric matrices with geometrically and arithmetically distributed singular values respectively. The code to regenerate this figure can be found in Appendix~\ref{app:time-improvement}.}
\label{fig:time-improvement}
\end{figure}

The final testing is to see how much time the mixed precision algorithm can save, we can use the methodology in Section~\ref{sec:qr-ns-performance} to produce Figure~\ref{fig:time-improvement}. For large $n$, the benefit of using preconditioner is significant for both distributions of the singular values.

Also, besides using the figure to visualize the difference in time, we can calculate how much the algorithm is improved for each $n$ via
\begin{equation}\notag
  \text{Improvement} = \frac{t_{\text{normal}} - t_{\text{precondition}}}{t_{\text{normal}}},
\end{equation}
where $t_{\text{normal}}$ and $t_{\text{precondition}}$ are the time used by the cyclic-by-row Jacobi algorithm without and with preconditioning.

\begin{lstlisting}
>> format short
>> (t_normal_g - t_precond_g)./t_normal_g % geometrically distributed SVs
ans =
  0.7150    0.7574    0.7739    0.7747    0.7789
>> (t_normal_a - t_precond_a)./t_normal_a % Arithmetically distributed SVs
ans =
  0.7294    0.7513    0.7543    0.7619    0.7663
\end{lstlisting}
The entries are corresponding to $n = 100,200,\dots,500$. Based on the
data, we can conclude that by using the preconditioning technique stated in
Algorithm~\ref{alg:jacobi-preconditioned}, we saved roughly $75\%$ of time.


\chapter{Conclusion and Further Work}

We studied the classical Jacobi algorithm and the cyclic-by-row Jacobi algorithm and concluded that the cyclic-by-row Jacobi algorithm is much more efficient than the other one. Then, after performing implementation and testing on the QR factorization and the Newton--Schulz iteration, in order to orthogonalize an almost orthogonal matrix, we shall always use the Newton--Schulz iteration. Finally, we followed the preconditioning technique in~\ycite[2000]{DV2000Approximateeigenvectorsas} and successfully propose the mixed precision algorithm~\ref{alg:jacobi-preconditioned} for the real symmetric eigenproblem. By utilizing the \inline{eig} function at single precision and preconditioning, we can save roughly 75\% of the time compared to applying the Jacobi algorithm alone.

Further works of this thesis include:
\begin{enumerate}
  \item Rounding error analysis hasn't been conducted.
  \item There exist more ways to orthogonalize an almost orthogonal matrix. For example, the Cholesky QR factorization.
  \item Only one low precision level has been studied. These theories can apply to half precision as well. However, it requires more work on exploring a half precision version of \inline{eig} function in order to become faster than our Algorithm~\ref{alg:jacobi-preconditioned} which uses MATLAB built-in single precision. 
  \item There are other ways that can be used to speed up the Jacobi algorithm. For example, the threshold Jacobi algorithm~\ycite[1965, Section~5.12]{Wil1965AlgebraicEigenvalueProblem} and the block Jacobi procedure~\ycite[2013, Section~8.5.6]{van2013mc}.
\end{enumerate}

\appendix
\chapter{Supplementary Explanation}
\section{Algorithm~\ref{alg:jacobi-pair}}
\label{app:algorithm-jacobi-pair}
We have two different expressions for computing $t$ as shown in Table~\ref{tab:different-approaches-for-t}.
\begin{table}[ht]
  \centering 
  \caption{Two different approaches for computing $t$ in Algorithm~\ref{alg:jacobi-pair}.}
  \label{tab:different-approaches-for-t}
  \begin{tabular}{ccc}
    \toprule 
    & Approach 1 & Approach 2 \\ \midrule  
    $\tau \geq 0$ & $ t = - \tau + \sqrt{\tau^2 + 1}$ & $\displaystyle t = \frac{1}{\tau + \sqrt{1 + \tau^2}}$ \\ \midrule 
    $\tau < 0 $ &$ t = - \tau - \sqrt{\tau^2 + 1}$ & 
    $\displaystyle t = \frac{1}{\tau - \sqrt{1 + \tau^2}}$\\
    \bottomrule
  \end{tabular}
\end{table}

These two expressions are the same by simple calculation,
\begin{equation}
  t=\frac{1}{\tau+\sqrt{1+\tau^{2}}}=\frac{\tau-\sqrt{1+\tau^{2}}}{\left(\tau+\sqrt{1+\tau^{2}}\right)\left(\tau-\sqrt{1+\tau^{2}}\right)}=-\tau+\sqrt{1+\tau^{2}}.
\end{equation}

However, when we implement these into MATLAB, approach 1 will output $c = 1, s = 0$ for $a_{pq} \approx 10^{-8}$, and approach 2's output gives $s \neq 0$ even for $a_{pq} \approx 10^{-14}$. Namely, by using approach 1, the computed output is less accurate than the output from approach 2.

\chapter{Supplementary Code}

\section{Testing Matrices}\label{app:myrandsvd}
By using the following code, we can generate a random real symmetric matrix with predefined condition number.
\begin{lstlisting}
function A = my_randsvd(n, kappa, mode)
classname = 'double';
if isempty(mode)
  mode = 'Geo';
end
switch mode % Set up vector sigma of singular values.
  case {'Geometric','geometric','Geo','geo'} 
    % geometrically distributed singular values
    factor = kappa^(-1/(n-1));
    sigma = factor.^(0:n-1);
  case {'Arithmatic', 'arithmatic', 'Ari', 'ari'}
    % arithmetically distributed singular values
    sigma = ones(n,1) - (0:n-1)'/(n-1)*(1-1/kappa);
  otherwise
    error(message('MATLAB:randsvd:invalidMode'));
end
sigma = cast(sigma,classname); signs = sign(rand(1,n-2)); 
sigma(2:n-1) = sigma(2:n-1).* signs; 
Q = qmult(n,1,classname);
A = Q*diag(sigma)*Q';
A = (A + A')/2;
\end{lstlisting}

\section{Full Jacobi Algorithm Code}\label{app:code-full-jacobi}
\begin{lstlisting}
function [V,D,counter] = jacobi_eig(Aini,tol,type,maxiter)
%JACOBI_EIG   Jacobi Eigenvalue Algorithm
%   [A,V,counter] = JACOBI_EIG(Aini,tol,type,maxiter) computes the eigen-
%   decomposition of Aini, Aini = VDV', where V is orthogonal.
%Input:
%   Aini -------- Input matrix
%   tol --------- Desired tolerance, usually 1.1e-16 (double precision)
%   type -------- Type of Jacobi algorithm, either
%                   * 'classical': For classical Jacobi algorithm
%                   * 'cyclic': For cyclic-by-row Jacobi algorithm
%   maxiter ----- Maximum number of sweep, specially for cyclic-by-row
%               Jacobi algorithm. If undetermined, the maximum sweep will
%               be set as 10 sweeps.
%Output:
%   V ------------ Orthogonal matrix
%   D ------------ Resulting matrix after apply Jacobi algorithm.
%   counter ------ Number of iteration
if nargin <= 2, error('Not enough input arguments.');
elseif nargin >= 5, error('Too many input arguments.');
else
  % select method, 'classical' or 'cyclic'
  switch type
    case {'classical','Classical'}
      method = 1;
    case {'cyclic','Cyclic'}
      method = 2;
    otherwise
      error("The Jacobi algorithm avaliable are " + ...
        "'classical' and 'cyclic'")
  end
  switch method
    case 1
      [V,D,counter] = jacobi_classical(Aini,tol);
    case 2
      switch nargin
        case 4
          maxiteration = maxiter;
          [V,D,counter] = jacobi_cyclic(Aini,tol,maxiteration);
        case 3
          maxiteration = 10;
          [V,D,counter] = jacobi_cyclic(Aini,tol,maxiteration);
        otherwise
          error('Not enough input arguments.');
      end
  end
end
end

% classical Jacobi algorithm
function [V,A,counter] = jacobi_classical(A,tol)
counter = 0; n = length(A); V = eye(n); done_rot = true;
tol1 = tol * norm(A);
while done_rot
  if isint(counter/(n*n)), A = (A + A')/2; end
  done_rot = false; [p,q] = maxoff(A);
  if abs(A(p,q)) >  tol1 * sqrt(abs(A(p,p) * A(q,q)))
    counter = counter + 1; done_rot = true;
    [c,s] = jacobi_pair(A,p,q);
    J = [c,s;-s,c]; 
    A([p,q],:) = J'*A([p,q],:);
    A(:,[p,q]) = A(:,[p,q]) * J;
    V(:,[p,q]) = V(:,[p,q]) * J;
  else
    A = diag(diag(A));
    break;
  end
end
end

% Cyclic-by-row Jacobi algorithm
function [V,A,iter] = jacobi_cyclic(A,tol,maxiter)
n = length(A); V = eye(n); iter = 0; done_rot = true;
while done_rot && iter < maxiter
  done_rot = false;
  for p = 1:n-1
    for q = p+1:n
      if abs(A(p,q)) > tol * sqrt(abs(A(p,p)*A(q,q)))
        done_rot = true;
        [c,s] = jacobi_pair(A,p,q);
        J = [c,s;-s,c];
        A([p,q],:) = J'*A([p,q],:);
        A(:,[p,q]) = A(:,[p,q]) * J;
        V(:,[p,q]) = V(:,[p,q]) * J;
      end
    end
  end
  if done_rot 
    A = (A + A')/2; iter = iter + 1; 
  else
    A = diag(diag(A));
    return;
  end
%   fprintf('off(A) = %d, sweep = %d\n', off(A), iter);
end
end

% compute the Frobenius norm of the off-diagonal entries of the input matrix
function num = off(A)
n = length(A); A(1:n+1:n*n) = 0;
num = norm(A,'fro');
end

% find the index of maximum off-diagonal entry of the input matrix A
function [p,q] = maxoff(A)
n = length(A); % dimension of matrix A
A(1:n+1:n*n) = 0; A = abs(A); % clear the diagonal entries
[val, idx1] = max(A);
[~, q] = max(val); p = idx1(q);
end

% calculate the Givens rotation matrix's entries (c,s).
function [c,s] = jacobi_pair(A,p,q)
if A(p,q) == 0
  c = 1; s = 0;
else
  tau = (A(q,q)-A(p,p))/(2*A(p,q));
  if tau >= 0
    t = 1/(tau + sqrt(1+tau*tau));
  else
    t = 1/(tau - sqrt(1+tau*tau));
  end
  c = 1/sqrt(1+t*t);
  s = t*c;
end
end
\end{lstlisting}

\section{Code for Figure~\ref{fig:classical-cyclic-compare}}\label{app:code-for-fig1}
Figure~\ref{fig:classical-cyclic-compare} can be regenerated using the following routine 
\begin{lstlisting}
clc; clear; close all; format short e;
condi1 = 100; tol = 2^(-53);
rng(1,'twister');
for n = 1e2:1e2:1e3
  A = my_randsvd(n,condi1,'geo');
  fprintf('iteration of n %d/%d\n', n/1e2, 1e3/1e2);
  tic 
  [V1,D1,iter1] = jacobi_eig(A,tol,'Classical');
  t_classical_n(n/1e2) = toc;
  tic
  [V2,D2,iter2] = jacobi_eig(A,tol,'Cyclic',20);
  t_cyclic_n(n/1e2) = toc;
end
nplot = 1e2:1e2:1e3;
n = 2.5e2; 
rng(1,'twister');
for condi = 1e3:1e3:1e5
  A = my_randsvd(n,condi,'geo');
  fprintf('iteration of cond %d/%d...',condi/1e3,1e5/1e3);
  tic
  [V1,D1,iter1] = jacobi_eig(A,tol,'Classical');
  t_classical_condi(condi/1e3) = toc;
  tic
  [V2,D2,iter2] = jacobi_eig(A,tol,'Cyclic',20);
  t_cyclic_condi(condi/1e3) = toc;
  fprintf('complete \n');
end
condiplot = 1e3:1e3:1e5;
close all;
subplot(1,2,1); hold off;
plot(nplot,t_classical_n,'-^','LineWidth',2);
hold on;
plot(nplot,t_cyclic_n,'-^','LineWidth',2);
xlabel('Dimension $(n)$');
ylabel('Time Used (sec)');
title('Fixed $\kappa(A) = 100$');
legend('Classical Jacobi','Cyclic Jacobi','Location','northwest','FontSize',20);
axis square
subplot(1,2,2); hold off;
plot(condiplot,t_classical_condi,'-^','LineWidth',1.5);
hold on;
plot(condiplot,t_cyclic_condi,'-^','LineWidth',1.5);
xlabel('$\kappa(A)$');
ylabel('Time Used (sec)');
title('Fixed $n = 250$');
legend('Classical Jacobi','Cyclic Jacobi','Location','northwest','FontSize',20);
axis square; axis([0,1e5,0,25])
\end{lstlisting}

\section{Code for Figure~\ref{fig:orthog}}\label{app:fig:orthog}
Figure~\ref{fig:orthog} can be regenerated using the following routine 
\begin{lstlisting}
clc; clear; close all; format short e; 
ud = 2^(-53); n = 1e2:1e2:3e3; condi = 100; rng(1,'twister');
for i = 1:length(n)
  fprintf('iteration %d/%d\n', i, length(n));
  A = gallery('randsvd', n(i), -condi, 3);
  [Q,~] = myqr(A); 
  orthog(i) = norm(Q'*Q - eye(n(i)));
end
semilogy(n,orthog,'-^r',LineWidth=1.5); 
hold on; semilogy(n, n.*ud, '-^k',LineWidth=1.5);
legend("$\|Q^TQ - I_n\|$", '$nu_d$','Location','northwest');
xlabel('Dimension $(n)$')
\end{lstlisting}

\section{Code for Figure~\ref{fig:seconditer}}\label{app:seconditer}
Figure~\ref{fig:seconditer} can be regenerated using the following routine 
\begin{lstlisting}
clc; clear; close all;
n = 10:10:8e4;
sq = (sqrt(1+n .* eps(single(1/2))) - 1).^4;
sq_ref = n .* eps(1/2);
plot(n,sq); hold on;
plot(n,sq_ref,'--');
plot(52400,5.81757e-12,'^k')
legend('$\|X_2 - U\|_2$','$n\cdot u_d$','Location','northwest');
xlabel('Dimension $(n)$')
\end{lstlisting}

\section{Code for Figure~\ref{fig:nsorthog}}\label{app:fig:nsorthog}
Figure~\ref{fig:nsorthog} can be regenerated using the following routine 
\begin{lstlisting}
clc; clear; close all; format short e; rng(1,'twister');
ud = 2^(-53); 
n = 1e2:1e2:3e3; condi = 1e2;
for i = 1:length(n)
  fprintf('Start iteration %d/%d,\n',i,length(n));
  A = my_randsvd( n(i), condi, 'geo');
  [Q,~] = eig(single(A)); Q = double(Q);
  Q1 = newton_schulz(Q); 
  orthog(i) = norm(Q1'*Q1 - eye(n(i)));
end
semilogy(n,orthog,'-^r',LineWidth=1.5); 
hold on; 
semilogy(n, n.*ud, '-^k',LineWidth=1.5);
legend("$\|Q_d^*Q_d - I_n\|$",'$nu_d$','Location','northwest'); 
xlabel('Dimension (n)'); title('Fix $\kappa(A) = 100$'); axis square
\end{lstlisting}

\section{Code for Figure~\ref{fig:acccompare}}\label{app:fig:acccompare}
Figure~\ref{fig:acccompare} can be regenerated using the following routine 
\begin{lstlisting}
clc; clear; close all; format short e; rng(1,'twister');
ud = 2^(-53); n = 1e2:1e2:3e3; condi = 1e2;
for i = 1:length(n)
  fprintf('iteration %d/%d\n',i,length(n));
  A = my_randsvd(n(i),condi,'geo');  
  [Q,~] = eig(single(A)); Q = double(Q);
  Q1 = newton_schulz(Q); 
  [Q2,~] = qr(Q);
  [Q3,~] = myqr(Q);
  orthogns(i) = norm(Q1'*Q1 - eye(n(i)));
  orthogqr(i) = norm(Q2'*Q2 - eye(n(i)));
  orthogmyqr(i) = norm(Q3'*Q3 - eye(n(i)));
end
semilogy(n,orthogns,'-^',LineWidth=1.5); hold on; 
semilogy(n,orthogqr,'-^',LineWidth=1.5); 
semilogy(n,orthogmyqr,'-^',LineWidth=1.5); 
legend("Newton-Schulz", ...
  "\texttt{qr()} from MATLAB", ...
  "My own QR code",...
  'interpreter','latex', ...
  'Location','northwest', ...
  'FontSize',20); 
xlabel('Dimension $(n)$');
ylabel('$\|Q_d^*Q_d - I_n\|$')
axis([0,3000,0,2e-14])
axis square
\end{lstlisting}

\section{Code for Figure~\ref{fig:timecompare}}
\label{app:fig:timecompare}
Figure~\ref{fig:timecompare} can be regenerated using the following routine 
\begin{lstlisting}
clc; clear; close all; format short e; 
ud = 2^(-53); n = 1e2:1e2:3e3; condi = 100;
time_ns = zeros(length(n),1); time_qr = time_ns;
rng(1,'twister'); % RNG
for i = 1:length(n)
  fprintf('iteration %d/%d\n',i,length(n));
  A = my_randsvd(n(i),condi,'geo');
  [Q,~] = eig(single(A)); Q = double(Q);
  tic; [Q1] = newton_schulz(Q); time_ns(i) = toc;
  tic; [Q2,~] = myqr(Q); time_qr(i) = toc;
end
semilogy(n,time_ns,'-^',LineWidth=1.5); hold on; 
semilogy(n,time_qr,'-^',LineWidth=1.5); 
legend('Newton--Schulz','QR Factorization','Location','northwest','interpreter','latex','FontSize',20); 
xlabel('Dimension $(n)$')
ylabel('Time (sec)')
axis square  
\end{lstlisting}

\section{Mixed Precision Jacobi Algorithm Code}\label{app:mixed-precision-jacobi}
\begin{lstlisting}
function [Q,D,iter] = jacobi_precondi(A,tol,maxiter)
%JACOBI_PRECONDI(A,tol,maxiter) compute the spectral decomposition of A using cyclic-by-row Jacobi algorithm with preconditioning. 
%% Preliminaries
switch nargin
  case 2, max_it = 10; case 3, max_it = maxiter;
  otherwise, error('Please check the number of inputs');
end
[n,m] = size(A); 
if m ~= n, error('Input matrix should be square matrix'), end
if issymmetric(A) == 0, error('Input matrix should be symmetric'); end
%% Spectral decomposition at single precision
[Q_low,~] = eig(single(A)); Q_low = double(Q_low);
%% Newton--Schulz iteration to orthogonalize
Q_d = newton_schulz(Q_low);
%% Cyclic-by-row Jacobi algorithm on preconditioned A
A_cond = Q_d' * A * Q_d;
[V,D,iter] = jacobi_eig(A_cond,tol,'Cyclic',max_it);
%% Output
Q = Q_d * V;
\end{lstlisting}

\section{Code for Figure~\ref{fig:precondition-accuracy} }\label{app:precondition-accuracy}
Figure~\ref{fig:precondition-accuracy} can be regenerated using the following routine 
\begin{lstlisting}
clc; clear; format short e;
n = 50:50:1e3; condi = 500; nplot = n;
rng(1,'twister');
for i = 1:length(n)
  fprintf('iteration %d/%d...',i,length(n))
  A_geo = my_randsvd(n(i),condi,'geo');
  A_ari = my_randsvd(n(i),condi,'ari');
  [V1,D1,iter1] = jacobi_precondi(A_geo,eps(1/2));
  [V2,D2,iter2] = jacobi_precondi(A_ari,eps(1/2));
  accuracy_n_geo(i) = norm(A_geo * V1 - V1 * D1);
  accuracy_n_ari(i) = norm(A_ari * V2 - V2 * D2);
  reference_n(i) = n(i) * norm(A_geo) * eps(1/2);
  fprintf('complete\n');
end
%%
subplot(1,2,1); hold off;
semilogy(nplot,accuracy_n_geo,'-^'); hold on;
semilogy(nplot,accuracy_n_ari,'-^');
semilogy(nplot,nplot.*eps(1/2),'--k');
xlabel('Dimension $(n)$'); 
ylabel('$\|AQ - Q\Lambda\|_2/\|A\|_2$');
legend('Geometric','Arithmetic','$nu_d\|A\|_2$','Location','southeast','FontSize',20);
title('Fixed $\kappa(A) = 500$')
axis square
%%
condiplot = 1e2:1e2:1e4; n = 5e2;
for i = 1:length(condiplot)
  fprintf('iteration %d/%d...',i,length(condiplot))
  A_geo = my_randsvd(n,condiplot(i),'geo');
  A_ari = my_randsvd(n,condiplot(i),'ari');
  [V1,D1,iter1] = jacobi_precondi(A_geo,eps(1/2));
  [V2,D2,iter2] = jacobi_precondi(A_ari,eps(1/2));
  accuracy_condi_geo(i) = norm(A_geo * V1 - V1 * D1);
  accuracy_condi_ari(i) = norm(A_ari * V2 - V2 * D2);
  reference_n(i) = n * norm(A_geo) * eps(1/2);
  fprintf('complete\n');
end
%%
close all;
subplot(1,2,1); hold off;
semilogy(nplot,accuracy_n_geo,'-^'); hold on;
semilogy(nplot,accuracy_n_ari,'-^');
semilogy(nplot,nplot.*eps(1/2),'--k');
xlabel('Dimension $(n)$'); 
ylabel('$\|AQ - Q\Lambda\|_2/\|A\|_2$');
legend('Geometric','Arithmetic','$nu_d$','Location','southeast','FontSize',20);
title('Fixed $\kappa(A) = 500$')
axis square
subplot(1,2,2); hold off;
semilogy(condiplot,accuracy_condi_geo,'-^'); hold on;
semilogy(condiplot,accuracy_condi_ari,'-^');
semilogy(condiplot,reference_n,'--k'); 
xlabel('$\kappa(A)$'); 
ylabel('$\|AQ - Q\Lambda\|_2/\|A\|_2$');
legend('Geometric','Arithmetic','$nu_d$','Location','northeast','FontSize',20);
title('Fixed $n = 500$')
axis square; axis([0,1e4,0,8e-14]);
\end{lstlisting}

\section{Code for Figure~\ref{p3:fig:approx-eig-test}}
\label{p3:sec:figure1}
\begin{lstlisting}
clc; clear all; close all; figure(1);
n = [50,250,500];
for i = 1:3
  subplot(2,2,i)
  testing_1(n(i));
  if i == 1
    axis([0,1e4,0,6e-6]);
    title('Fixed $n=50$');
  elseif i == 2
    axis([0,1e4,0,3.8e-5]);
    title('Fixed $n = 250$');
  elseif i == 3
    axis([0,1e4,0,8e-5]);
    title('Fixed $n = 500$');
  end
  pause(0.1);
end
subplot(2,2,4); hold off
condi = 1e2; n1 = 10:10:1e3; offA = zeros(1,length(n1)); rng(1,'twister');
for i = 1:length(n1)
  fprintf('start %d/%d\n',i,length(n1));
  A_geo = my_randsvd(n1(i),condi,'geo');
  A_ari = my_randsvd(n1(i),condi,'ari');
  [Q1,~] = eig(single(A_geo)); Q1 = double(Q1);
  Q1_d = newton_schulz(Q1);
  offA_geo(i) = off(Q1_d' * A_geo * Q1_d)/norm(A_geo);
  [Q2,~] = eig(single(A_ari)); Q2 = double(Q2);
  Q2_d = newton_schulz(Q2);
  offA_ari(i) = off(Q2_d' * A_ari * Q2_d)/norm(A_ari);
end
semilogy(n1,offA_geo); hold on;
semilogy(n1,offA_ari);
semilogy(n1,n1 .* double(eps(single(1/2))),'--k');
legend('Geometric','Arithmetic','$nu_\ell$','Location','southeast');
xlabel('$\kappa(A)$'); ylabel('off$(A_{cond})/\|A\|_2$','Interpreter','latex');
set(gca,'FontSize',20);
axis square;
title('Fixed $\kappa(A) = 100$')
function testing_1(n)
condi = 10:10:1e4; offA_geo = zeros(1,length(condi)); offA_ari = offA_geo;
rng(1,'twister');
for i = 1:length(condi)
  fprintf('start %d/%d\n',i,length(condi));
  A_geo = my_randsvd(n,condi(i),'geo');
  [Q1,~] = eig(single(A_geo)); Q1 = double(Q1);
  Q1_d = newton_schulz(Q1);
  offA_geo(i) = off(Q1_d' * A_geo * Q1_d)/norm(A_geo);
  A_ari = my_randsvd(n,condi(i),'ari');
  [Q2,~] = eig(single(A_ari)); Q2 = double(Q2);
  Q2_d = newton_schulz(Q2);
  offA_ari(i) = off(Q2_d' * A_ari * Q2_d)/norm(A_ari);
end
semilogy(condi,offA_geo,'LineWidth',1.5); hold on;
semilogy(condi,offA_ari);
semilogy(condi,n * double(eps(single(1/2))) * ones(1,length(condi)),'--k','LineWidth',1.5);
legend('Geometric','Arithmetic','$nu_\ell$');
xlabel('$\kappa(A)$'); ylabel('off$(A_{cond})/\|A\|_2$','Interpreter','latex');
set(gca,'FontSize',20);
axis square;
end
\end{lstlisting}

\section{Code for Figure~\ref{fig:typical-sweep}}\label{app:typical-sweep}
Figure~\ref{fig:typical-sweep} can be regenerated using the following routine
\begin{lstlisting}
clear; clc; rng(1,'twister');
n = 50:50:1e3; condi = 500; nplot = n;
for i = 1:length(n)
  A_geo = my_randsvd(n(i),condi,'geo');
  [V1,D1,iter1] = jacobi_precondi(A_geo,eps(1/2));
  iteration_n_geo(i) = iter1;
  A_ari = my_randsvd(n(i),condi,'ari');
  [V2,D2,iter2] = jacobi_precondi(A_ari,eps(1/2));
  iteration_n_ari(i) = iter2;
  fprintf('iteration %d/%d complete\n',i,length(n));
end
condiplot = 1e2:1e2:1e4; n = 5e2;
for i = 1:length(condiplot)
  A_geo = my_randsvd(n,condiplot(i),'geo');
  [V1,D1,iter1] = jacobi_precondi(A_geo,eps(1/2));
  iteration_condi_geo(i) = iter1;
  A_ari = my_randsvd(n,condiplot(i),'ari');
  [V2,D2,iter2] = jacobi_precondi(A_ari,eps(1/2));
  iteration_condi_ari(i) = iter2;
  fprintf('iteration %d/%d complete\n',i,length(condiplot));
end
close all
subplot(1,2,1);
semilogy(nplot,iteration_n_geo,'-^'); hold on;
semilogy(nplot,iteration_n_ari,'-^');
xlabel('Dimension $(n)$');
ylabel('Number of iterations til converge');
title('Fixed $\kappa(A) = 500$')
legend('Geometric','Arithmetic','Location','northeast');
axis square; axis([0,1000,0,4]);
subplot(1,2,2);
semilogy(condiplot,iteration_condi_geo,'-^'); hold on;
semilogy(condiplot,iteration_condi_ari,'-^');
legend('Geometric','Arithmetic','Location','northeast');
xlabel('$\kappa(A)$'); 
ylabel('Number of iterations til converge');
title('Fixed $n = 500$')
axis square; axis([0,1e4,0,5]);
\end{lstlisting}

\section{Code for Figure~\ref{fig:time-improvement}}\label{app:time-improvement}
Figure~\ref{fig:time-improvement} can be regenerated using the following routine
\begin{lstlisting}
clc; clear; close all; rng(1,'twister');
n = 1e2:1e2:5e2; condi = 500; tol = 2^(-53);
for i = 1:length(n)
  Ag = my_randsvd(n(i),condi,'geo');
  Aa = my_randsvd(n(i),condi,'ari');
  tic; [V1,D1,iter1] = jacobi_precondi(Ag,tol); 
  t_precond_g(i) = toc;
  tic; [V2,D2,iter2] = jacobi_eig(Ag,tol,'Cyclic'); 
  t_normal_g(i) = toc;
  tic; [Va1,Da1,iter1] = jacobi_precondi(Aa,tol); 
  t_precond_a(i) = toc;
  tic; [Va2,Da2,iter2] = jacobi_eig(Aa,tol,'Cyclic'); 
  t_normal_a(i) = toc;
  fprintf('Iteration %d/%d complete\n',i,length(n));
end
subplot(1,2,1)
plot(n,t_normal_g,'-^');
hold on;
plot(n,t_precond_g,'-^');
xlabel('Dimension $(n)$')
ylabel('Time Used (sec)')
legend('Without Preconditioning', 'With Preconditioning','Location','northwest');
title('Geometrically distributed test matrices')
axis square  
subplot(1,2,2)
plot(n,t_normal_a,'-^');
hold on;
plot(n,t_precond_a,'-^');
xlabel('Dimension $(n)$')
ylabel('Time Used (sec)')
legend('Without Preconditioning', 'With Preconditioning','Location','northwest');
title('Arithmetically distributed test matrices')
axis square  
\end{lstlisting}


\bibliographystyle{nj-plain}
\bibliography{strings, bib}

\def\noopsort#1{}\def\hbk{hardback}\def\pbk{paperback}
\begin{thebibliography}{10}

\bibitem{Ben1965iterativemethodcomputing}
Adi Ben-Israel.
\newblock \href{https://doi.org/10.1090/s0025-5718-1965-0179915-5}{An iterative
  method for computing the generalized inverse of an arbitrary matrix}.
\newblock {\em Math. Comp.}, 19\penalty0 (91):\penalty0 452--455, 1965.

\bibitem{Bha1997MatrixAnalysis}
Rajendra Bhatia.
\newblock \href{http://doi.org/10.1007/978-1-4612-0653-8}{{\em Matrix
  {Analysis}}}.
\newblock Spring{\-}er-Ver{\-}lag, New York, 1997.
\newblock xi+347 pp.
\newblock ISBN 978-1-4612-6857-4.

\bibitem{DHT2001AnalysisCholeskyMethod}
Philip~I. Davies, Nicholas~J. Higham, and Fran\c{c}oise Tisseur.
\newblock \href{https://doi.org/10.1137/S0895479800373498}{Analysis of the
  {Cholesky} method with iterative refinement for solving the symmetric
  definite generalized eigenproblem}.
\newblock {\em SIAM J. Matrix Anal. Appl.}, 23\penalty0 (2):\penalty0 472--493,
  2001.

\bibitem{DV1992Jacobismethodis}
James Demmel and Kre\v{s}imir Veseli\'{c}.
\newblock \href{https://doi.org/10.1137/0613074}{{Jacobi}'s method is more
  accurate than {QR}}.
\newblock {\em SIAM J. Matrix Anal. Appl.}, 13\penalty0 (4):\penalty0
  1204--1245, 1992.

\bibitem{DV2000Approximateeigenvectorsas}
Zlatko Drma\v{c} and Kre\v{s}imir Veseli\'{c}.
\newblock \href{https://doi.org/10.1016/s0024-3795(00)00046-x}{{Approximate}
  eigenvectors as preconditioner}.
\newblock {\em Linear Algebra Appl.}, 309\penalty0 (1):\penalty0 191--215,
  2000.

\bibitem{FH1955Somemetricinequalities}
Ky~Fan and A.~J. Hoffman.
\newblock \href{http://www.jstor.org/stable/2032662}{Some metric inequalities
  in the space of matrices}.
\newblock {\em Proc. Amer. Math. Soc.}, 6\penalty0 (1):\penalty0 111--116,
  1955.

\bibitem{FH1960cyclicJacobimethod}
George~E. Forsythe and Peter Henrici.
\newblock \href{https://www.jstor.org/stable/1993275}{The cyclic {Jacobi}
  method for computing the principal values of a complex matrix}.
\newblock {\em Trans. Amer. Math. Soc.}, 94\penalty0 (1):\penalty0 1--23, 1960.

\bibitem{GvdV2000Eigenvaluecomputation20th}
Gene~H. Golub and Henk~A. van~der Vorst.
\newblock \href{https://doi.org/10.1016/b978-0-444-50617-7.50010-0}{Eigenvalue
  computation in the 20th century}.
\newblock {\em J. Comput. Appl. Math.}, 123\penalty0 (1):\penalty0 209--239,
  2001.

\bibitem{van2013mc}
Gene~H. Golub and Charles~F. Van~Loan.
\newblock
  \href{https://www.press.jhu.edu/books/title/10678/matrix-computations}{{\em
  Matrix Computations}}.
\newblock Johns {Hopkins} studies in the mathematical sciences. 4th edition,
  Johns Hopkins University Press, Baltimore, MD, USA, 2013.
\newblock ISBN 978-1-4214-0794-4.

\bibitem{Gre1953Computingeigenvalueseigenvectors}
Robert~T. Gregory.
\newblock \href{https://doi.org/10.1090/s0025-5718-1953-0057643-6}{{Computing}
  eigenvalues and eigenvectors of a symmetric matrix on the {ILLIAC}}.
\newblock {\em Math. Comp.}, 7\penalty0 (44):\penalty0 215--220, 1953.

\bibitem{Har1991sharpquadraticconvergence}
Vjeran Hari.
\newblock \href{https://doi.org/10.1007/BF01385728}{{On} sharp quadratic
  convergence bounds for the serial {Jacobi} methods}.
\newblock {\em Numer. Math.}, 60\penalty0 (1):\penalty0 375--406, 1991.

\bibitem{Hen1964ElementsNumericalAnalysis}
Peter Henrici.
\newblock {\em {Elements} of {Numerical} {Analysis}}.
\newblock Wiley, New York, USA, 1964.
\newblock ISBN 978-0471372417.

\bibitem{Hig1986Computingpolardecomposition}
Nicholas~J. Higham.
\newblock \href{https://doi.org/10.1137/0907079}{Computing the polar
  decomposition{\textemdash}with applications}.
\newblock {\em SIAM J. Sci. Statist. Comput.}, 7\penalty0 (4):\penalty0
  1160--1174, 1986.

\bibitem{high:ASNA2}
Nicholas~J. Higham.
\newblock \href{http://doi.org/10.1137/1.9780898718027}{{\em Accuracy and
  Stability of Numerical Algorithms}}.
\newblock 2nd edition, Society for Industrial and Applied Mathematics,
  Philadelphia, PA, USA, January 2002.
\newblock xxx+680 pp.
\newblock ISBN 0-89871-521-0.

\bibitem{higham08fm}
Nicholas~J. Higham.
\newblock \href{http://doi.org/10.1137/1.9780898717778}{{\em Functions of
  Matrices: {Theory} and Computation}}.
\newblock Society for Industrial and Applied Mathematics, Philadelphia, PA,
  USA, 2008.
\newblock xx+425 pp.
\newblock ISBN 978-0-898716-46-7.

\bibitem{Jacobi-original-paper-1846}
Carl G.~J. Jacobi.
\newblock \href{https://doi.org/10.1515/crll.1846.30.51}{{\"{U}}ber ein
  leichtes verfahren die in der {Theorie} der {S{\"a}cularst{\"o}rungen}
  vorkommenden {Gleichungen} numerisch aufzul{\"o}sen}.
\newblock 1846\penalty0 (30):\penalty0 51--94, 1846.

\bibitem{MATLAB:2022}
MATLAB.
\newblock version 9.12.0.2009381 (r2022a), 2022.

\bibitem{Mir1960Symmetricgaugefunctions}
Leon Mirsky.
\newblock \href{https://doi.org/10.1093/qmath/11.1.50}{Symmetric gauge
  functions and unitarily invariant norms}.
\newblock {\em Q. J. Math.}, 11\penalty0 (1):\penalty0 50--59, 1960.

\bibitem{IEEE2019}
\rule{0.35in}{1pt}.
\newblock \href{https://doi.org/10.1109/IEEESTD.2019.8766229}{{IEEE} standard
  for floating--point arithmetic}.
\newblock {\em {IEEE} {Std} 754-2019 ({Revision} of {IEEE} 754--2008)}, pages
  1--84, 2019.

\bibitem{1961-first-quadratic-convergence-Schu}
Arnold Sch{\"o}nhage.
\newblock \href{https://doi.org/10.1007/bf01386036}{Zur konvergenz des
  {Jacobi}-verfahrens}.
\newblock {\em Numer. Math.}, 3\penalty0 (1):\penalty0 374--380, 1961.

\bibitem{Sch1933IterativeBerechungder}
G{\"{u}}nther Schulz.
\newblock \href{https://doi.org/10.1002/zamm.19330130111}{{Iterative Berechung}
  der reziproken {Matrix}}.
\newblock {\em Z. Angew. Math. Mech.}, 13\penalty0 (1):\penalty0 57--59, 1933.

\bibitem{vKem1966quadraticconvergencespecial}
H.~P.~M. van Kempen.
\newblock \href{https://doi.org/10.1007/bf02165225}{On the quadratic
  convergence of the special cyclic {Jacobi} method}.
\newblock {\em Numer. Math.}, 9\penalty0 (1):\penalty0 19--22, 1966.

\bibitem{Wil1962Notequadraticconvergence}
James~H. Wilkinson.
\newblock \href{https://doi.org/10.1007/bf01386321}{Note on the quadratic
  convergence of the cyclic {Jacobi} process}.
\newblock {\em Numer. Math.}, 4\penalty0 (1):\penalty0 296--300, 1962.

\bibitem{Wil1965AlgebraicEigenvalueProblem}
James~H. Wilkinson.
\newblock {\em The Algebraic Eigenvalue Problem}.
\newblock Oxford University Press, Oxford, UK, 1965.
\newblock xviii+662 pp.
\newblock ISBN 978-0198534181.

\end{thebibliography}

\end{document}
